\documentclass{article}
\usepackage{rotating}
\usepackage{graphics}
\usepackage[all]{xy}
\xyoption{v2}
\xyoption{knot}

\author {Stefan Forcey}

\title {Higher Dimensional Enrichment}

\newtheorem{theorem}{Theorem}
\newtheorem{definition}{Definition}

\newcommand{\MySection}[1]
{\section{ #1}}

\newcommand{\bcal}[1]{\mbox{\boldmath${\cal {#1}}$}}

\newenvironment{proof}%!!!!!!
{\noindent{\bf Proof}}%!!!!!!

\begin{document}
\SloppyCurves{

\maketitle
\begin{abstract}
 Lyubashenko has described enriched 2--categories as categories enriched over 
     ${\cal V}$--Cat, the 2--category of categories enriched over a symmetric monoidal ${\cal V}$. Here 
     I generalize this to a $k$--fold monoidal ${\cal V}$. The latter is defined 
     as by Balteanu, Fiedorowicz, Schw${\rm \ddot a}$nzl and Vogt but with the addition of making 
     visible the coherent associators  $\alpha^i$.
     The symmetric case can easily be recovered. The introduction of this paper
     proposes a recursive definition of ${\cal V}$--$n$--categories and their morphisms.
     Then I consider the special case of  
     ${\cal V}$--2--categories  and give the details of the proof that with their 
     morphisms these form the
    structure of a 3--category. 
%50 pages; complete definition of V-n-categories and morphisms between them; journal style and notation--note tensor^(1) replaces box^(2); updated references
%Lyubashenko has described enriched 2-categories as categories enriched over V-Cat, the 2-category of categories enriched over a symmetric monoidal V. I have generalized this to the k-fold monoidal V. The symmetric case can easily be recovered. The introduction of this paper proposes a recursive definition of V-n-categories and their morphisms. Then I consider the special case of V-2-categories and give the details of the proof that with their morphisms these form the structure of a 3-category. 
\end{abstract}

  %{intro2}
      
      %\clearpage
                %\newpage
            \MySection{Introduction}
      
    There is an ongoing massive effort to link category theory and geometry, just a part 
    of the broad undertaking known as categorification as described by Baez and Dolan in \cite{Baez1}.
    This effort has as a partial goal that
    of understanding
    the categories and functors that correspond to loop spaces and their
    associated topological functors. 
    Progress towards this goal has been advanced greatly by the recent work of Balteanu,
    Fiedorowicz, Schw${\rm \ddot a}$nzl, and Vogt in \cite{Balt}
    where they show a direct correspondence between $k$-fold monoidal
    categories and  $k$-fold loop spaces through the
    categorical nerve. 
    
    As I pursued part of a plan to relate the enrichment functor to topology, I
    noticed that the concept of higher dimensional
    enrichment would be important in its relationship to double, triple and
    further iterations of delooping.
    The concept of enrichment over a monoidal category is well known, and 
    enriching over the category of categories enriched over 
    a monoidal category is  
    defined, for the case of symmetric categories, in the paper on
    $A_{\infty}$--categories
    by Lyubashenko,  \cite{Lyub}. It seems that it is a good idea to generalize his definition first
    to the case of an iterated monoidal
    base category
    and then to 
    define ${\cal V}$--$(n+1)$--categories as
    categories enriched over ${\cal V}$--$n$--Cat, the $(k-n)$--fold monoidal
    strict 
    $(n+1)$--category of ${\cal V}$--$n$--categories where $k<n \in {\mathbf N}$. Of course the facts
    implicit in this last statement must be verified.
    At each stage of successive enrichments, the number of monoidal products
    should decrease and the categorical dimension
    should increase, both by one. This is motivated by topology. When we
    consider the loop space of a topological space,
    we see that paths (or  1--cells) in the original 
    are now points (or objects) in the derived space. There is also now
    automatically a product structure
    on the points in the derived space, where multiplication is given by
    concatenation of loops. Delooping is the inverse 
    functor here, and thus involves shifting objects to the status of
    1--cells and decreasing the number of ways to multiply.

    The concept of a $k$--fold monoidal strict $n$--category  
    is easy enough to define as a tensor object in a category of $(k-1)$--fold
    monoidal $n$--categories with cartesian product. Thus the products and
    accompanying associator and interchange transformations
    are strict bi--$n$--functors and $n$--natural transformations respectively. That
    this sort of structure 
    ($(k-n)$--fold monoidal strict 
    $n+1$ category) is possessed by  ${\cal V}$--$n$--Cat for ${\cal V}$ $k$--fold monoidal is shown 
    for $n = 1$ and all $k$ in my paper \cite{forcey1}. 
    A full inductive proof covering all $n,k$ is a work in progress, and this
    paper fills in one of the gaps; specifically 
    showing how the categorical dimension is increased when the base
    2--category is ${\cal V}$--Cat. I hope to see
    eventually how to take the long proof contained here and turn it into
    part (1) of the induction step in the full proof.
    Part (2), showing how the number of products is decreasing, should
    actually be easier since rather than involving
    new structure at each step it only seems to require a repeating of the
    axiom checking that is done in the 
    initial case. This is described for the case of ${\cal V}$--$2$--Cat
    being  $(k-2)$--fold monoidal in my previously
    mentioned paper. In general the decrease is engineered by a shift in index--we define 
    new products ${\cal V}$--$n$--Cat$\times {\cal V}$--$n$--Cat $\to {\cal V}$--$n$--Cat by 
    using cartesian products of object sets and letting hom--objects of the
    $i$th product of enriched $n$--categories
    be the $(i+1)$th product of hom--objects of the component categories.
    Symbolically,
    $$({\cal A} \otimes^{(n)}_{i} {\cal B})((A,B),(A',B')) = {\cal A}(A,A')\otimes^{(n-1)}_{i+1} {\cal B}(B,B').$$
    The superscript $(n)$ is not necessary since the product is defined by
    context, 
            but I insert it to make clear at what
            level of enrichment the product is occurring. 
    Defining the necessary natural transformations for this new product as
    ``based upon'' the old ones, and the checking of 
    the axioms that define their structure is briefly mentioned later on in
    this paper and more fully  described in \cite{forcey1}  
     for certain cases.
      
    The definition of a category enriched over ${\cal V}$--$n$--Cat is simply
    stated by describing the process as enriching over 
    ${\cal V}$--$n$--Cat with the first of the $k-n$ ordered products. In  detail this means that
    
    \begin{definition}  \label{super} A (small, strict) ${\cal V}$--{\it (n+1)--category}
    ${\bcal U}$ consists of
            
            \begin{enumerate}
                \item A set of objects $\left|{\bcal U}\right|$
                \item For each pair of objects $A,B \in \left|{\bcal U}\right|$ a
                ${\cal V}$--$n$--category ${\bcal U}(A,B).$
                \item For each triple of objects $A,B,C \in 
            \left|{\bcal U}\right|$ a ${\cal V}$--$n$--functor 
            $${\bcal M}_{ABC}:{\bcal U}(B,C) \otimes^{(n)}_{1} {\bcal U}(A,B) \to {\bcal U}(A,C)$$
                \item For each object $A \in \left|{\bcal U}\right|$ a ${\cal V}$--$n$--functor
            $${\bcal J}_A: {\bcal I}^{(n)} \to {\bcal U}(A,A)$$ 
                \item Axioms: The ${\cal V}$--$n$--functors that play the
    role of composition and identity
                obey commutativity of a pentagonal diagram (associativity
    axiom) and of two triangular diagrams (unit axioms).
               This amounts to saying that the 
            functors given by the two legs of each diagram are equal. 
            $$
  	          
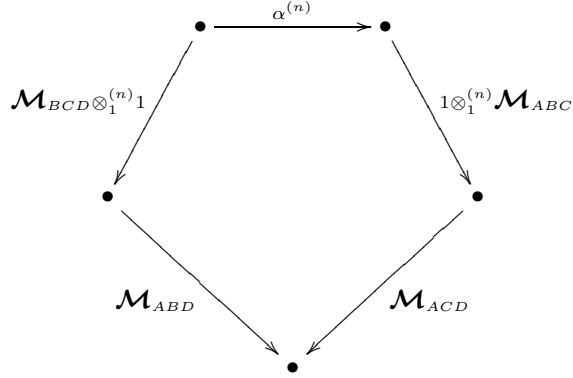
\xymatrix{
  	           &\bullet
  	           \ar[rr]^{\alpha^{(n)}}
  	           \ar[ddl]_{{\bcal M}_{BCD} \otimes^{(n)}_{1} 1}
  	           &&\bullet
  	           \ar[ddr]^{1 \otimes^{(n)}_{1} {\bcal M}_{ABC}}&\\\\
  	           \bullet
  	           \ar[ddrr]_{{\bcal M}_{ABD}}
  	           &&&&\bullet
  	             \ar[ddll]^{{\bcal M}_{ACD}}
  	             \\\\&&\bullet&&&
             }
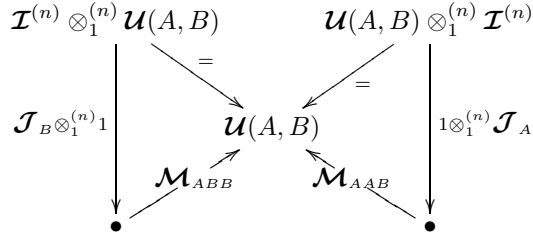
$$
             $$
  	              \xymatrix@C=-5pt{
  	              {\bcal I}^{(n)}\otimes^{(n)}_1 {\bcal U}(A,B)
  	              \ar[rrd]^{=}
  	              \ar[dd]_{{\bcal J}_{B}\otimes^{(n)}_1 1}
  	              &&&&{\bcal U}(A,B)\otimes^{(n)}_1 {\bcal I}^{(n)}
  	              \ar[dd]^{{1}\otimes^{(n)}_1 {\bcal J}_{A}}
  	              \ar[lld]^{=}\\
  	              &&{\bcal U}(A,B)\\
  	              \bullet
  	              \ar[rru]|{{\bcal M}_{ABB}}
  	              &&&&\bullet
  	              \ar[llu]|{{\bcal M}_{AAB}}
  	              }
          $$
             \end{enumerate}
   \end{definition}      
    
    This definition requires that there be a definitions of the unit ${\bcal I}^{(n)}$ and of ${\cal
    V}$--$n$--functors in place. Since the axioms talk
    about equality, the ${\cal V}$--$n$--natural transformations and various
    higher morphisms are not directly required.
    They will be utilized however, as soon as the word strict is removed from
     Definition~\ref{super}. I will sketch definitions
    of morphisms in 
    a moment. First, since we hope to realize monoidal structure on ${\cal
    V}$--$n$--Cat, it needs a recursively defined 
    unit ${\cal V}$--$n$--category.
    
    \begin{definition} The unit object in ${\cal V}$--$n$--Cat is the ${\cal
    V}$--$n$--category 
      ${\bcal I}^{(n)}$ with one object $\mathbf{0}$ and with ${\bcal
    I}^{(n)}(\mathbf{0},\mathbf{0}) = {\bcal I}^{(n-1)},$
      where ${\bcal I}^{(n-1)}$ is the unit object in ${\cal
    V}$--$(n-1)$--Cat. Of course we let ${\bcal I}^{(0)}$ be $I$ in
      ${\cal V}.$ Also ${\bcal M}_{000} = {\bcal J}_{0} = 1_{{\bcal I^{(n)}}}.$
      \end{definition}
      
    Now we can define the morphisms:  
    
              \begin{definition} \label{onecell}
              For two ${\cal V}$--$n$--categories ${\bcal U}$ and ${\bcal W}$
    a ${\cal V}${\it --$n$--functor} $T:{\bcal U}\to{\bcal W}$ is
    a function on objects 
              $\left|{\bcal U}\right|\to\left|{\bcal W}\right|$ and a
    family of ${\cal V}$--$(n-1)$--functors
              $T_{UU'}:{\bcal U}(U,U')\to{\bcal W}(TU,TU').$ These latter
    obey commutativity of the usual diagrams.
              \begin{enumerate}
               \item For $U,U',U'' \in \left|{\bcal U}\right|$
               $$
                 \xymatrix@C=65pt@R=65pt{
                 &\bullet
                 \ar[rr]^{{\bcal M}_{UU'U''}}
                 \ar[d]^{T_{U'U''} \otimes^{(n-1)}_1 T_{UU'}}
                 &&\bullet
                 \ar[d]^{T_{UU''}}&\\
                 &\bullet
                 \ar[rr]_{{\bcal M}_{(TU)(TU')(TU'')}}
                 &&\bullet
                 }
               $$
               \item
               $$
                 \xymatrix{
                 &&\bullet
                 \ar[dd]^{T_{UU}}\\
                 {\bcal I}^{(n-1)}
                 \ar[rru]^{{\bcal J}_U}
                 \ar[rrd]_{{\bcal J}_{TU}}\\
                 &&\bullet
                 }
               $$
              Here a ${\cal V}$--$0$--functor is just a morphism in ${\cal V}.$
             \end{enumerate}
              \end{definition}
              
             The 1--cells we have just defined play a special role in the  
    definition of a general $k$--cell for $k\ge 2$. 
    
              \begin{definition} \label{kcell}
	                  A ${\cal V}${\it --n--k--cell} $\alpha$ between $(k-1)$--cells
	        $\psi^{k-1}$ and $\phi^{k-1}$, written 
	                   $$\alpha:\psi^{k-1}\to \phi^{k-1}:\psi^{k-2}\to \phi^{k-2}: ... :\psi^{2}\to \phi^{2}:F\to G:{\bcal U}\to {\bcal W}$$ 
	                  where $F$ and $G$ are  ${\cal V}$--$n$--functors,
	                  is a function sending
	                  each $U \in  \left|{\bcal U}\right|$ to a 
	                  ${\cal V}$--$((n-k)+1)$--functor 
	                 $$\alpha_{U}: {\bcal I}^{((n-k)+1)}\to {\bcal W}(FU,GU)(\psi^{2}_U(\mathbf{0}),\phi^{2}_U(\mathbf{0}) ... (\psi^{k-1}_U(\mathbf{0}),\phi^{k-1}_U(\mathbf{0}))$$
	                  in such a way that we have commutativity of the following
	        diagram. Note that the final (curved) equal sign is implied recursively by the diagram for the $(k-1)$--cells. 
                  \end{definition}

                  \noindent
		  	          \begin{center}
		  	          \resizebox{6.47in}{!}{ %!!!!
		  	          %%%%\resizebox{7in}{!}{
		  	          \begin{small}
		  	                     $$
		  	                      \xymatrix@C=-105pt@R=35pt{ %!!!!
		  	                      %%%%\xymatrix@C=-275pt@R=35pt{
		  	                      &
		  	                      &*\txt{${\bcal W}(FU',GU')(\psi^{2}_{U'}(\mathbf{0}),\phi^{2}_{U'}(\mathbf{0})) ... (\psi^{k-1}_{U'}(\mathbf{0}),\phi^{k-1}_{U'}(\mathbf{0}))$\\$\otimes^{((n-k)+1)}_{k-1} {\bcal W}(FU,FU')(F(x_2),F(y_2))...(F(x_{k-1}),F(y_{k-1}))  $}
		  	                      \ar[dr]^-{{\bcal M}}
		  	                       \\
		  	                    &
		  	                    {\bcal I}^{((n-k)+1)} \otimes^{((n-k)+1)}_{k-1} {\bcal U}(U,U')(x_2,y_2)...(x_{k-1},y_{k-1})\text{ }
		  	                      \ar[ru]^<<{\alpha_{U'} \otimes^{((n-k)+1)}_{k-1} F\text{ }}
		  	                    &&{\bcal W}(FU,GU')(\psi^{2}_{U'}(\mathbf{0})F(x_2),\phi^{2}_{U'}(\mathbf{0})F(y_2)) ... (\psi^{k-1}_{U'}(\mathbf{0})F(x_{k-1}),\phi^{k-1}_{U'}(\mathbf{0})F(y_{k-1})) 
		  	                    \ar@{=}@/^2pc/[dd]
		  	                    \\
		  	                    \text{ }\text{ }\text{ }{\bcal U}(U,U')(x_2,y_2)...(x_{k-1},y_{k-1})\text{ }\text{ }
		  	                      \ar[ru]^{=}
		  	                    \ar[rd]_{=}
		  	                    \\
		  	                    &
		  	                    {\bcal U}(U,U')(x_2,y_2)...(x_{k-1},y_{k-1}) \otimes^{((n-k)+1)}_{k-1} {\bcal I}^{((n-k)+1)}\text{ }
		  	                      \ar[rd]_<<{G \otimes^{((n-k)+1)}_{k-1} \alpha_{U}\text{ }}
		  	                    &&{\bcal W}(FU,GU')(G(x_2)\psi^{2}_U(\mathbf{0}),G(y_2)\phi^{2}_U(\mathbf{0})...(G(x_{k-1})\psi^{k-1}_U(\mathbf{0}),G(y_{k-1})\phi^{k-1}_U(\mathbf{0}))   
		  	                    \\
		  	                      &
		  	                      &*\txt{${\bcal W}(GU,GU')(G(x_2),G(y_2))...(G(x_{k-1}),G(y_{k-1}))$\\$  \otimes^{((n-k)+1)}_{k-1} {\bcal W}(FU,GU)(\psi^{2}_U(\mathbf{0}),\phi^{2}_U(\mathbf{0}) ... (\psi^{k-1}_U(\mathbf{0}),\phi^{k-1}_U(\mathbf{0}))$}
		  	                      \ar[ur]_-{{\bcal M}}
		  	                      }
		  	                    $$
		  	                    \end{small}
		  	                    }
	                    \end{center}

           Thus for a given value of $n$ there are $k$--cells up to $k = n+1$,  making   ${\cal V}$--$n$--Cat a potential $(n+1)$--category.
           This last definition is best grasped by looking at examples. The cases for $n = 1, 2$ are given in detail in the following section.

   % {Def}
    \clearpage
          \newpage
        \MySection{Review of Definitions}
        
     In this section I briefly review the definitions of a category enriched over a monoidal category ${\cal V}$, a
     category enriched over an iterated monoidal category, and an enriched 2--category.
     I begin with the basic definitions of enrichment, included due to how often they 
      are referred to and followed as models in the rest of the paper. This first set of definitions
      can be found with more detail in \cite{Kelly} and \cite{EK1}.
      
      \begin{definition} For our purposes a {\it monoidal category} is a category ${\cal V}$
        together with a functor 
        $\otimes: {\cal V}\times{\cal V}\to{\cal V}$  and an object $I$ such that
        \begin{enumerate}
        \item $\otimes$ is  associative up to the coherent natural transformations $\alpha$. The coherence
        axiom is given by the commuting pentagon
        
        \noindent
	\begin{center}
	\resizebox{5.5in}{!}{
        $$
        \xymatrix@C=-25pt{
        &((U\otimes V)\otimes W)\otimes X \text{ }\text{ }
        \ar[rr]^{ \alpha_{UVW}\otimes 1_{X}}
        \ar[ddl]^{ \alpha_{(U\otimes V)WX}}
        &&\text{ }\text{ }(U\otimes (V\otimes W))\otimes X
        \ar[ddr]^{ \alpha_{U(V\otimes W)X}}&\\\\
        (U\otimes V)\otimes (W\otimes X)
        \ar[ddrr]|{ \alpha_{UV(W\otimes X)}}
        &&&&U\otimes ((V\otimes W)\otimes X)
        \ar[ddll]|{ 1_{U}\otimes \alpha_{VWX}}
        \\\\&&U\otimes (V\otimes (W\otimes X))&&&
       }
       $$
       }
	                    \end{center}
	                    
      \item $I$ is a strict $2$-sided unit for $\otimes$.
      \end{enumerate}
      \end{definition}
      \begin {definition} \label{V:Cat} A (small) ${\cal V}$ {\it--Category} ${\cal A}$ is a set $\left|{\cal A}\right|$ of 
      {\it objects}, 
      a {\it hom-object} ${\cal A}(A,B) \in \left|{\cal V}\right|$ for
      each pair of objects of ${\cal A}$, a family of {\it composition morphisms} $M_{ABC}:{\cal A}(B,C)
      \otimes{\cal A}(A,B)\to{\cal A}(A,C)$ for each triple of objects, and an {\it identity element} $j_{A}:I\to{\cal A}(A,A)$ for each object.
      The composition morphisms are subject to the associativity axiom which states that the following pentagon commutes
      
      \noindent
      	          \begin{center}
	          \resizebox{5.5in}{!}{
      $$
      \xymatrix@C=-15pt{
      &({\cal A}(C,D)\otimes {\cal A}(B,C))\otimes {\cal A}(A,B)\text{ }\text{ }
      \ar[rr]^{\scriptstyle \alpha}
      \ar[dl]^{\scriptstyle M \otimes 1}
      &&\text{ }\text{ }{\cal A}(C,D)\otimes ({\cal A}(B,C)\otimes {\cal A}(A,B))
      \ar[dr]^{\scriptstyle 1 \otimes M}&\\
      {\cal A}(B,D)\otimes {\cal A}(A,B)
      \ar[drr]^{\scriptstyle M}
      &&&&{\cal A}(C,D)\otimes {\cal A}(A,C)
      \ar[dll]^{\scriptstyle M}
      \\&&{\cal A}(A,D))&&&
      }$$
      }
     	                    \end{center}
     	                    
      and to the unit axioms which state that both the triangles in the following diagram commute
      
      $$
        \xymatrix{
        I\otimes {\cal A}(A,B)
        \ar[rrd]^{=}
        \ar[dd]_{j_{B}\otimes 1}
        &&&&{\cal A}(A,B)\otimes I 
        \ar[dd]^{1\otimes j_{A}}
        \ar[lld]^{=}\\
        &&{\cal A}(A,B)\\
        {\cal A}(B,B)\otimes {\cal A}(A,B)
        \ar[rru]^{M_{ABB}}
        &&&&{\cal A}(A,B)\otimes {\cal A}(A,A)
        \ar[llu]^{M_{AAB}}
        }
     $$

     \end{definition}
    
     \begin{definition} \label{enriched:funct} For ${\cal V}$--categories 
     ${\cal A}$ and ${\cal B}$, a ${\cal V}$--$functor$ $T:{\cal A}\to{\cal B}$ is a function
      $T:\left| {\cal A} \right| \to \left| {\cal B} \right|$ and a family of 
      $hom-object$ morphisms $T_{AB}:{\cal A}(A,B) \to {\cal B}(TA,TB)$ in ${\cal V}$ indexed by 
      pairs $A,B \in \left| {\cal A} \right|$.
      The usual rules for a functor that state $T(f \circ g) = Tf \circ Tg$ 
      and $T1_{A} = 1_{TA}$ become in the enriched setting, respectively, the commuting diagrams
     
     $$
      \xymatrix{
      &{\cal A}(B,C)\otimes {\cal A}(A,B)
      \ar[rr]^{\scriptstyle M}
      \ar[d]^{\scriptstyle T \otimes T}
      &&{\cal A}(A,C)
      \ar[d]^{\scriptstyle T}&\\
      &{\cal B}(TB,TC)\otimes {\cal B}(TA,TB)
      \ar[rr]^{\scriptstyle M}
      &&{\cal B}(TA,TC)
      }
     $$
    and
     $$
      \xymatrix{
      &&{\cal A}(A,A)
      \ar[dd]^{\scriptstyle T_{AA}}\\
      I
      \ar[rru]^{\scriptstyle j_{A}}
      \ar[rrd]_{\scriptstyle j_{TA}}\\
      &&{\cal B}(TA,TA).
      }
     $$
    ${\cal V}$--functors can be composed to form a category called ${\cal V}$--Cat. This category
    is actually enriched over $\mathbf{Cat}$, the category of (small) categories with cartesian product. 
       \end{definition}
       
    \begin{definition} \label{enr:nat:trans}
    For ${\cal V}$--functors $T,S:{\cal A}\to{\cal B}$ a ${\cal V}$--{\it natural  
     transformation} $\alpha:T \to S:{\cal A} \to {\cal B}$
    is an $\left| {\cal A} \right|$--indexed family of 
    morphisms $\alpha_{A}:I \to {\cal B}(TA,SA)$ satisfying the ${\cal V}$--naturality
    condition expressed by the commutativity of 
    
    $$
      \xymatrix{
      &I \otimes {\cal A}(A,B)
      \ar[rr]^-{\scriptstyle \alpha_{B} \otimes T_{AB}}
      &&{\cal B}(TB,SB) \otimes {\cal B}(TA,TB)
      \ar[rd]^-{\scriptstyle M}
    \\
      {\cal A}(A,B)
      \ar[ru]^{=}
    \ar[rd]_{=}
      &&&&{\cal B}(TA,SB)
    \\
      &{\cal A}(A,B) \otimes I
      \ar[rr]_-{\scriptstyle S_{AB} \otimes \alpha_{A}}
      &&{\cal B}(SA,SB) \otimes {\cal B}(TA,SA)
      \ar[ru]^-{\scriptstyle M}
      }
     $$
    
     \end{definition}
     
     For two ${\cal V}$--functors  $T,S$ to be equal is to say $TA = SA$ for all $A$ 
     and for the ${\cal V}$--natural isomorphism $\alpha$ between them to have 
     components $\alpha_{A} = j_{TA}$. This latter implies equality of the hom--object morphisms: 
     $T_{AB} = S_{AB}$ for all pairs of objects. The implication is seen by combining the second diagram in 
     Definition~\ref{V:Cat} with all the diagrams
    in Definitions~\ref{enriched:funct} and \ref{enr:nat:trans}.

    The fact that ${\cal V}$--Cat has the structure of a 2--category is demonstrated in \cite{Kelly}. Now
    we review the transfer to enriching over a $k$--fold monoidal category. The latter sort of category was
    developed and defined in \cite{Balt}. The authors describe its structure as arising from its description
    as a monoid in the
        category of $(k-1)$--fold monoidal categories. Here is that definition altered only slightly 
    to make visible the coherent associators as in \cite{forcey1}. In that paper I describe its structure as 
    arising from its description
    as a tensor object in the
        category of $(k-1)$--fold monoidal categories
        
        \begin{definition} An $n${\it -fold monoidal category} is a category ${\cal V}$
        with the following structure. 
        \begin{enumerate}
        \item There are $n$ distinct multiplications
        $$\otimes_1,\otimes_2,\dots, \otimes_n:{\cal V}\times{\cal V}\to{\cal V}$$
        for each of which the associativity pentagon commutes
        
        \noindent
		          \begin{center}
	          \resizebox{5.5in}{!}{
        $$
        \xymatrix@C-=2pt@C=-30pt{
        &((U\otimes_i V)\otimes_i W)\otimes_i X \text{ }\text{ }
        \ar[rr]^{ \alpha^{i}_{UVW}\otimes_i 1_{X}}
        \ar[ddl]^{ \alpha^{i}_{(U\otimes_i V)WX}}
        &&\text{ }\text{ }(U\otimes_i (V\otimes_i W))\otimes_i X
        \ar[ddr]^{ \alpha^{i}_{U(V\otimes_i W)X}}&\\\\
        (U\otimes_i V)\otimes_i (W\otimes_i X)
        \ar[ddrr]^{ \alpha^{i}_{UV(W\otimes_i X)}}
        &&&&U\otimes_i ((V\otimes_i W)\otimes_i X)
        \ar[ddll]^{ 1_{U}\otimes_i \alpha^{i}_{VWX}}
        \\\\&&U\otimes_i (V\otimes_i (W\otimes_i X))&&&
        }
        $$
        }
		                    \end{center}

        ${\cal V}$ has an object $I$ which is a strict unit
        for all the multiplications.
        \item For each pair $(i,j)$ such that $1\le i<j\le n$ there is a natural
        transformation
        $$\eta^{ij}_{ABCD}: (A\otimes_j B)\otimes_i(C\otimes_j D)\to
        (A\otimes_i C)\otimes_j(B\otimes_i D).$$
        \end{enumerate}
        These natural transformations $\eta^{ij}$ are subject to the following conditions:
        \begin{enumerate}
        \item[(a)] Internal unit condition: 
        $\eta^{ij}_{ABII}=\eta^{ij}_{IIAB}=1_{A\otimes_j B}$
        \item[(b)] External unit condition:
        $\eta^{ij}_{AIBI}=\eta^{ij}_{IAIB}=1_{A\otimes_i B}$
        \item[(c)] Internal associativity condition: The following diagram commutes
        
        $$
          \diagram
          ((U\otimes_j V)\otimes_i (W\otimes_j X))\otimes_i (Y\otimes_j Z)
          \xto[rrr]^{\eta^{ij}_{UVWX}\otimes_i 1_{Y\otimes_j Z}}
          \ar[d]^{\alpha^i}
          &&&\bigl((U\otimes_i W)\otimes_j(V\otimes_i X)\bigr)\otimes_i (Y\otimes_j Z)
          \dto^{\eta^{ij}_{(U\otimes_i W)(V\otimes_i X)YZ}}\\
          (U\otimes_j V)\otimes_i ((W\otimes_j X)\otimes_i (Y\otimes_j Z))
          \dto^{1_{U\otimes_j V}\otimes_i \eta^{ij}_{WXYZ}}
          &&&((U\otimes_i W)\otimes_i Y)\otimes_j((V\otimes_i X)\otimes_i Z)
          \ar[d]^{\alpha^i \otimes_j \alpha^i}
          \\
          (U\otimes_j V)\otimes_i \bigl((W\otimes_i Y)\otimes_j(X\otimes_i Z)\bigr)
          \xto[rrr]^{\eta^{ij}_{UV(W\otimes_i Y)(X\otimes_i Z)}}
          &&& (U\otimes_i (W\otimes_i Y))\otimes_j(V\otimes_i (X\otimes_i Z))
          \enddiagram
          $$
         \item[(d)] External associativity condition: The following diagram commutes
          $$
          \diagram
          ((U\otimes_j V)\otimes_j W)\otimes_i ((X\otimes_j Y)\otimes_j Z)
          \xto[rrr]^{\eta^{ij}_{(U\otimes_j V)W(X\otimes_j Y)Z}}
          \ar[d]^{\alpha^j \otimes_i \alpha^j}
          &&& \bigl((U\otimes_j V)\otimes_i (X\otimes_j Y)\bigr)\otimes_j(W\otimes_i Z)
          \dto^{\eta^{ij}_{UVXY}\otimes_j 1_{W\otimes_i Z}}\\
          (U\otimes_j (V\otimes_j W))\otimes_i (X\otimes_j (Y\otimes_j Z))
          \dto^{\eta^{ij}_{U(V\otimes_j W)X(Y\otimes_j Z)}}
          &&&((U\otimes_i X)\otimes_j(V\otimes_i Y))\otimes_j(W\otimes_i Z)
          \ar[d]^{\alpha^j}
          \\
          (U\otimes_i X)\otimes_j\bigl((V\otimes_j W)\otimes_i (Y\otimes_j Z)\bigr)
          \xto[rrr]^{1_{U\otimes_i X}\otimes_j\eta^{ij}_{VWYZ}}
          &&& (U\otimes_i X)\otimes_j((V\otimes_i Y)\otimes_j(W\otimes_i Z))
          \enddiagram
          $$
        
        \item[(e)] Finally it is required that for each triple $(i,j,k)$ satisfying
        $1\le i<j<k\le n$ 
        the giant hexagonal interchange diagram commutes.
        \end{enumerate}
        \end{definition}
        \begin{small}
        
        \noindent
		          \begin{center}
	          \resizebox{6.5in}{!}{
        $$
        \xymatrix@C=-103pt{
        &((A\otimes_k A')\otimes_j (B\otimes_k B'))\otimes_i((C\otimes_k C')\otimes_j (D\otimes_k D'))
        \ar[ddl]|{\eta^{jk}_{AA'BB'}\otimes_i \eta^{jk}_{CC'DD'}}
        \ar[ddr]|{\eta^{ij}_{(A\otimes_k A')(B\otimes_k B')(C\otimes_k C')(D\otimes_k D')}}
        \\\\
        ((A\otimes_j B)\otimes_k (A'\otimes_j B'))\otimes_i((C\otimes_j D)\otimes_k (C'\otimes_j D'))
        \ar[dd]|{\eta^{ik}_{(A\otimes_j B)(A'\otimes_j B')(C\otimes_j D)(C'\otimes_j D')}}
        &&((A\otimes_k A')\otimes_i (C\otimes_k C'))\otimes_j((B\otimes_k B')\otimes_i (D\otimes_k D'))
        \ar[dd]|{\eta^{ik}_{AA'CC'}\otimes_j \eta^{ik}_{BB'DD'}}
        \\\\
        ((A\otimes_j B)\otimes_i (C\otimes_j D))\otimes_k((A'\otimes_j B')\otimes_i (C'\otimes_j D'))
        \ar[ddr]|{\eta^{ij}_{ABCD}\otimes_k \eta^{ij}_{A'B'C'D'}}
        &&((A\otimes_i C)\otimes_k (A'\otimes_i C'))\otimes_j((B\otimes_i D)\otimes_k (B'\otimes_i D'))
        \ar[ddl]|{\eta^{jk}_{(A\otimes_i C)(A'\otimes_i C')(B\otimes_i D)(B'\otimes_i D')}}
        \\\\
        &((A\otimes_i C)\otimes_j (B\otimes_i D))\otimes_k((A'\otimes_i C')\otimes_j (B'\otimes_i D'))
        }
        $$
        }
		                    \end{center}
	                    
        \end{small}

     The authors of \cite{Balt} remark that a symmetric monoidal category is $n$-fold monoidal for all $n$.
     This they demonstrate by letting
        $$\otimes_1=\otimes_2=\dots=\otimes_n=\otimes$$
        and defining (associators added by myself)
        $$\eta^{ij}_{ABCD}=\alpha^{-1}\circ (1_A\otimes \alpha)\circ (1_A\otimes (c_{BC}\otimes 1_D))\circ (1_A\otimes \alpha^{-1})\circ \alpha$$
        for all $i<j$. Here $c_{BC}: B\otimes C \to C\otimes B$ is the symmetry natural transformation.
     This provides the hint that enriching over a $k$--fold monoidal category is precisely a generalization of enriching over
     a symmetric category. In the symmetric case, to define a product in ${\cal V}$--Cat, we need $c_{BC}$ in order to create a
     middle exchange morphism $m$.  
     To describe products in ${\cal V}$--Cat for ${\cal V}$ $k$--fold monoidal we simply use $m=\eta$.
     
     Categories enriched over $k$--fold monoidal ${\cal V}$ are carefully defined in \cite{forcey1}, where they  
     are shown to be the objects of a $(k-1)$--fold monoidal 2--category. Here we need only the definitions. Simply put,
     a category enriched over a $k$--fold monoidal ${\cal V}$ is a category enriched in the usual sense over $({\cal V}, \otimes_1, I, \alpha).$
     The other $k-1$ products in ${\cal V}$ are used up in the structure of ${\cal V}$--Cat. I will always denote the product(s) in ${\cal V}$--Cat
       with a superscript in parentheses that corresponds to the level of enrichment of the components of their domain.
       The product(s) in ${\cal V}$ 
    should logically then have a superscript (0) but I have suppressed this for brevity and to agree with my sources.
     For ${\cal V}$ $k$--fold monoidal we define the $i$th product of ${\cal V}$--categories 
     ${\cal A} \otimes^{(1)}_{i}{\cal B}$
        to have objects $\in \left|{\cal A}\right|\times \left|{\cal B}\right|$
        and to have hom--objects $\in \left|{\cal V}\right|$ given by
        
        $$({\cal A} \otimes^{(1)}_{i} {\cal B})((A,B),(A',B')) = {\cal A}(A,A')\otimes_{i+1} {\cal B}(B,B').$$
        
        Immediately we see that ${\cal V}$--Cat is $(k-1)$--fold monoidal by definition. (The full proof of
        this is in \cite{forcey1}.)
        The composition morphisms are 
        \begin{small}
        $$M_{(A,B)(A',B')(A'',B'')} : ({\cal A} \otimes^{(1)}_{i}{\cal
        B})((A',B'),(A'',B''))\otimes_{1}({\cal A} \otimes^{(1)}_{i}{\cal
        B})((A,B),(A',B'))\to ({\cal A} \otimes^{(1)}_{i}{\cal B})((A,B),(A'',B''))$$
        \end{small}
        %\clearpage
        given by
        $$
        \xymatrix{
        ({\cal A} \otimes^{(1)}_{i}{\cal B})((A',B'),(A'',B''))\otimes_{1}({\cal A} \otimes^{(1)}_{i}{\cal B})((A,B),(A',B'))
        \ar@{=}[d]\\
        ({\cal A}(A',A'')\otimes_{i+1}{\cal B}(B',B''))\otimes_{1}({\cal A}(A,A')\otimes_{i+1}{\cal B}(B,B'))
        \ar[d]_{\eta^{1,i+1}}\\
        ({\cal A}(A',A'')\otimes_{1}{\cal A}(A,A'))\otimes_{i+1}({\cal B}(B',B'')\otimes_{1}{\cal B}(B,B'))
        \ar[d]_{M_{AA'A''}\otimes_{2}M_{BB'B''}}\\
        ({\cal A}(A,A'')\otimes_{i+1}{\cal B}(B,B''))
        \ar@{=}[d]\\
        ({\cal A} \otimes^{(1)}_{i}{\cal B})((A,B),(A'',B''))
       }
    $$
        
        The identity element is given by $j_{(A,B)} = \xymatrix{I = I \otimes_{i+1} I
         							\ar[d]^{j_A \otimes_{i+1} j_B}
         							\\{\cal A}(A,A)\otimes_{i+1} {\cal B}(B,B)
         							\ar@{=}[d]
         							\\({\cal A} \otimes^{(1)}_{i}{\cal B})((A,B),(A,B))
     							}$
      
      The unit object in ${\cal V}$--1--Cat $={\cal V}$--Cat is the enriched category ${\bcal I}^{(1)}={\cal I}$ where $\left|{\cal I}\right| = \{0\}$ and 
       ${\cal I}(0,0) = I$. Of course $M_{000} = 1_I = j_0.$
       
       That each product $\otimes^{(1)}_i$ thus defined is a 2--bi--functor  ${\cal V}$--Cat $\times$ ${\cal V}$--Cat $\to$ ${\cal V}$--Cat
       is seen easily. Its action on enriched functors and natural transformations is to form formal products using $\otimes_{i+1}$ of their
       associated morphisms. That the result is a valid enriched functor or natural transformation always follows from the
       naturality of $\eta.$
       
       Associativity in ${\cal V}$--Cat
          must hold for each $\otimes^{(1)}_{i}$. The components of the 2--natural
          isomorphism $\alpha^{(1)i}$
          $$\alpha^{(1)i}_{{\cal A}{\cal B}{\cal C}}: ({\cal A} \otimes^{(1)}_{i} {\cal B})
           \otimes^{(1)}_{i} {\cal C} \to {\cal A} \otimes^{(1)}_{i} ({\cal B} \otimes^{(1)}_{i}
          {\cal C})$$
          are ${\cal V}$--functors
          that send ((A,B),C) to (A,(B,C)) and whose hom-components 
          \begin{small}
          $$\alpha^{(1)i}_{{{\cal A}{\cal B}{\cal C}}_{((A,B),C)((A',B'),C')}}: [({\cal A} \otimes^{(1)}_{i} {\cal B}) \otimes^{(1)}_{i} {\cal C}](((A,B),C),((A',B'),C'))
          \to [{\cal A} \otimes^{(1)}_{i} ({\cal B} \otimes^{(1)}_{i} {\cal C})]((A,(B,C)),(A',(B',C')))$$
          \end{small}
          are given by 
          $$\alpha^{(1)i}_{{{\cal A}{\cal B}{\cal C}}_{((A,B),C)((A',B'),C')}}
          = \alpha^{i+1}_{{\cal A}(A,A'){\cal B}(B,B'){\cal C}(C,C')}.$$
               
     Now for the interchange 2--natural transformations $\eta^{(1)ij}$
       for $j\ge i+1$.
       We define the component morphisms $\eta^{(1)i,j}_{{\cal A}{\cal B}{\cal C}{\cal D}}$
       that make a 2--natural transformation between 2--functors. Each component must be an enriched functor.
       Their action on objects
       is to send $((A,B),(C,D)) \in \left|({\cal A}\otimes^{(1)}_{j} {\cal B})\otimes^{(1)}_{i} ({\cal C}\otimes^{(1)}_{j} {\cal D})\right|$
        to $((A,C),(B,D)) \in \left|({\cal A}\otimes^{(1)}_{i} {\cal C})\otimes^{(1)}_{j} ({\cal B}\otimes^{(1)}_{i} {\cal D})\right|$.
       The hom--object morphisms are given by
        $$\eta^{(1)i,j}_{{{\cal A}{\cal B}{\cal C}{\cal D}}_{(ABCD)(A'B'C'D')}} =
   \eta^{i+1,j+1}_{{\cal A}(A,A'){\cal B}(B,B'){\cal C}(C,C'){\cal D}(D,D')}.$$
       That the axioms regarding the associators and interchange transformations are all obeyed is established
       in \cite{forcey1}.
       
       We now define categories enriched over 
     ${\cal V}$--Cat. These are defined for the symmetric case in \cite{Lyub}. Here as in \cite{forcey1}
    the definition of 
       ${\cal V}$--2--category is generalized for ${\cal V}$ a $k$--fold monoidal category with $k\ge 2.$ 
       The definition for symmetric monoidal ${\cal V}$ can be easily recovered just by letting $\otimes_1=\otimes_2=\otimes$,
       $\alpha^2=\alpha^1 = \alpha$ and $\eta=m.$ 
       
       \begin{definition} A (small, strict) ${\cal V}$--{\it 2--category} ${\bcal U}$ consists of
       
       \begin{enumerate}
           \item A set of objects $\left|{\bcal U}\right|$
           \item For each pair of objects $A,B \in \left|{\bcal U}\right|$ a ${\cal V}$--category ${\bcal U}(A,B).$
       
       Of course then ${\bcal U}(A,B)$ consists of a set of objects (which play the role of the 1--cells in a 2--category) 
       and for each pair $f,g \in \left|{\bcal U}(A,B)\right|$ an object 
       ${\bcal U}(A,B)(f,g) \in {\cal V}$ (which plays the role 
       of the hom--set of 2--cells in a 2--category.) Thus the vertical composition morphisms of these $hom_{2}$--objects are in ${\cal V}$:
       $$M_{fgh}:{\bcal U}(A,B)(g,h) \otimes_{1} {\bcal U}(A,B)(f,g) \to 
       {\bcal U}(A,B)(f,h)$$
       
       Also, the vertical identity for a 1-cell object $a \in \left|{\bcal U}(A,B)\right|$ is  $j_{a} : I \to {\bcal U}(A,B)(a,a)$.
       The associativity and the units of vertical composition are then those given by the respective axioms of enriched categories.  
           \item For each triple of objects $A,B,C \in 
       \left|{\bcal U}\right|$ a ${\cal V}$--functor
       $${\cal M}_{ABC}:{\bcal U}(B,C) \otimes^{(1)}_{1} {\bcal U}(A,B) \to {\bcal U}(A,C)$$ 
       Often I repress the subscripts. We denote ${\cal M}(h,f)$ as $hf$. 
       
       The family of morphisms indexed by pairs of objects $(g,f),(g',f') \in 
       \left|{\bcal U}(B,C) \otimes^{(1)}_{1} {\bcal U}(A,B)\right|$ furnishes the direct analogue of horizontal composition of 2-cells
       as can be seen by observing their domain and range in ${\cal V}$:
       $${\cal M}_{ABC_{(g,f)(g',f')}}:[{\bcal U}(B,C) \otimes^{(1)}_{1} 
       {\bcal U}(A,B)]((g,f),(g',f')) \to {\bcal U}(A,C)(gf,g'f')$$
       Recall that 
       $$[{\bcal U}(B,C) \otimes^{(1)}_{1} {\bcal U}(A,B)]((g,f),(g',f')) = {\bcal U}(B,C)(g,g') \otimes_{2} {\bcal U}(A,B)(f,f').$$
       
          \item For each object $A \in \left|{\bcal U}\right|$ a ${\cal V}$--functor
       $${\cal J}_A: {\cal I} \to {\bcal U}(A,A)$$ 
       We denote ${\cal J}_A(0)$ as $1_{A}$. 
          \item (Associativity axiom of a strict ${\cal V}$--2--category.) We require a commuting pentagon.
          Since the morphisms are 
       ${\cal V}$--functors this amounts to saying that the 
       functors given by the two legs of the diagram are equal. 
       For objects we have the equality $(fg)h = f(gh).$  
       
       For the  hom--object morphisms we have the following family of commuting diagrams for associativity, where the first bullet represents
       $$[({\bcal U}(C,D)\otimes^{(1)}_{1} {\bcal U}(B,C)) \otimes^{(1)}_{1} {\bcal U}(A,B)](((f,g),h),((f',g'),h'))$$
       and the reader may fill in the others
       %$$
       % \xymatrix@C-=1pt{
       %  &[({\bcal U}(C,D)\otimes^{(1)}_{1} {\bcal U}(B,C)) \otimes^{(1)}_{1} {\bcal U}(A,B)](((f,g),h),((f',g'),h'))
       %  \ar[rr]^{\alpha^{(1)}}
       %  \ar[dl]^{{\cal M}_{(f,g)(f',g')} \otimes_{2} 1}
       %  &&[{\bcal U}(C,D) \otimes^{(1)}_{1} ({\bcal U}(B,C) \otimes^{(1)}_{1} {\bcal U}(A,B))]((f,(g,h)),(f',(g',h')))
       %  \ar[dr]^{1 \otimes_{2} {\cal M}_{(g,h)(g',h')}}&\\
       %  {\bcal U}(B,D)(fg,f'g') \otimes_{2} {\bcal U}(A,B)(h,h')
       %  \ar[drr]^{{\cal M}_{(fg,h)(f'g',h')}}
       %  &&&&{\bcal U}(C,D)(f,f') \otimes_{2} {\bcal U}(A,C)(gh,g'h')
       %    \ar[dll]^{{\cal M}_{(f,gh)(f',g'h')}}
       %    \\&&{\bcal U}(A,D)(fgh,f'g'h')&&&
       %    }$$
       
       $$
        
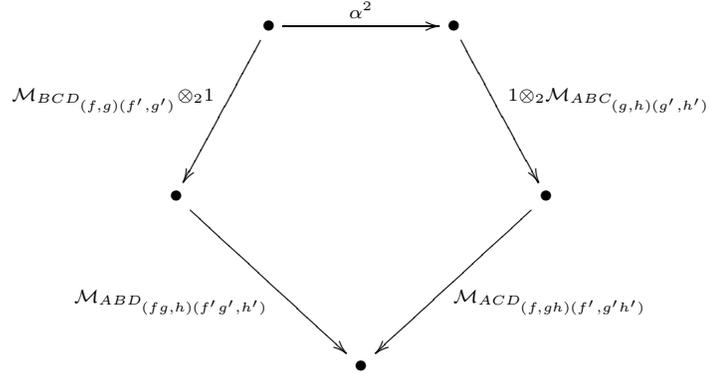
\xymatrix{
         &\bullet
         \ar[rr]^{\alpha^{2}}
         \ar[ddl]_{{\cal M}_{BCD_{(f,g)(f',g')}} \otimes_{2} 1}
         &&\bullet
         \ar[ddr]^{1 \otimes_{2} {\cal M}_{ABC_{(g,h)(g',h')}}}&\\\\
         \bullet
         \ar[ddrr]_{{\cal M}_{ABD_{(fg,h)(f'g',h')}}}
         &&&&\bullet
           \ar[ddll]^{{\cal M}_{ACD_{(f,gh)(f',g'h')}}}
           \\\\&&\bullet&&&
           }$$
      The heuristic diagram for this commutativity is
         $$
         
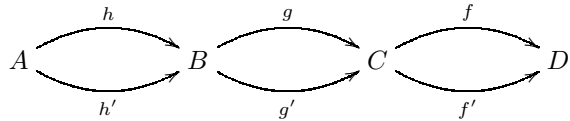
\xymatrix@R-=3pt{
         &&&\\
         A
         \ar@/^1pc/[rr]^h
         \ar@/_1pc/[rr]_{h'}
         &&B
         \ar@/^1pc/[rr]^g
         \ar@/_1pc/[rr]_{g'}
         &&C
         \ar@/^1pc/[rr]^f
         \ar@/_1pc/[rr]_{f'}
         &&D
         \\
         &&&\\
         }
    $$
      \item (Unit axioms of a strict ${\cal V}$--2--category.) We require commuting triangles. 
        For objects we have the equality $f1_{A} = f = 1_{B}f.$
        For the  unit morphisms we have that the triangles in the following diagram commute.

       $$
           \xymatrix@C=-5pt{
           [{\cal I}\otimes^{(1)}_1 {\bcal U}(A,B)]((0,f),(0,g))
           \ar[rrd]^{=}
           \ar[dd]_{{\cal J}_{B_{00}}\otimes_2 1}
           &&&&[{\bcal U}(A,B)\otimes^{(1)}_1 {\cal I}]((f,0),(g,0))
           \ar[dd]^{{1}\otimes_2 {\cal J}_{A_{00}}}
           \ar[lld]^{=}\\
           &&{\bcal U}(A,B)(f,g)\\
           [{\bcal U}(B,B)\otimes^{(1)}_1 {\bcal U}(A,B)]((1_B,f),(1_B,g))
           \ar[rru]|{{\cal M}_{ABB_{(1_B,f)(1_B,g)}}}
           &&&&[{\bcal U}(A,B)\otimes^{(1)}_1 {\bcal U}(A,A)]((f,1_A),(g,1_A))
           \ar[llu]|{{\cal M}_{AAB_{(f,1_A)(g,1_A)}}}
           }
        $$ 
        The heuristic diagrams for this commutativity are
     
     \noindent
     	          \begin{center}
	          \resizebox{6in}{!}{
     $$
      \xymatrix@R-=3pt{
      &\ar@{=>}[dd]^{1_{1_A}}&&\\
      A
      \ar@/^1pc/[rr]^{1_A}
      \ar@/_1pc/[rr]_{1_A}
      &&A
      \ar@/^1pc/[rr]^f
      \ar@/_1pc/[rr]_g
      &&B
      \\
      &&&&\text{ }&\\
      }
      =
      \xymatrix@R-=3pt{
           &&\\
          A
          \ar@/^1pc/[rr]^f
          \ar@/_1pc/[rr]_g
          &&B
         \\
        &&&
        }
      =
      \xymatrix@R-=3pt{
        &&&\ar@{=>}[dd]^{1_{1_B}}\\
        A
        \ar@/^1pc/[rr]^f
        \ar@/_1pc/[rr]_g
        &&B
        \ar@/^1pc/[rr]^{1_B}
        \ar@/_1pc/[rr]_{1_B}
        &&B
        \\
        &&&\\
      }
    $$ 
    }
    \end{center}
    
    \end{enumerate}
         \end{definition}
        Consequences of ${\cal V}$--functoriality of ${\cal M}$ and ${\cal J}$:
       First the ${\cal V}$--functoriality of ${\cal M}$ implies that the following (expanded) diagram commutes
       
       \noindent
       		          \begin{center}
	          \resizebox{6.5in}{!}{
       \begin{small}
       $$
        \xymatrix@C=-130pt{
        &({\bcal U}(B,C)(k,m)\otimes_1 {\bcal U}(B,C)(h,k))\otimes_2 ({\bcal U}(A,B)(g,l)\otimes_1 {\bcal U}(A,B)(f,g))
       \ar[rdd]^{M_{hkm}\otimes_2 M_{fgl}}
       \\\\
        ({\bcal U}(B,C)(k,m)\otimes_2 {\bcal U}(A,B)(g,l))\otimes_1 ({\bcal U}(B,C)(h,k)\otimes_2 {\bcal U}(A,B)(f,g))
       \ar[dd]^{{\cal M}_{ABC_{(k,g)(m,l)}}\otimes_1 {\cal M}_{ABC_{(h,f)(k,g)}}}
       \ar[ruu]^{\eta^{1,2}}
       &&{\bcal U}(B,C)(h,m)\otimes_2 {\bcal U}(A,B)(f,l)
       \ar[dd]^{{\cal M}_{ABC_{(h,f)(m,l)}}}
       \\\\
       {\bcal U}(A,C)(kg,ml)\otimes_1 {\bcal U}(A,C)(hf,kg)
       \ar[rr]^{M_{(hf)(kg)(ml)}}
       &&{\bcal U}(A,C)(hf,ml)
       \\
       }
       $$
       \end{small}
       }
       	                    \end{center}
	                    
       The heuristic diagram is
    $$
    \xymatrix@R-=16pt{
    &&&
    \\
    A
    \ar@/^2pc/[rr]^f
    \ar[rr]^g
    \ar@/_2pc/[rr]^l
    && B
    \ar@/^2pc/[rr]^h
    \ar[rr]^k
    \ar@/_2pc/[rr]^m
    && C\\
    &&&&&
    }
    $$       
       
       Secondly the ${\cal V}$--functoriality of ${\cal M}$ implies that the following (expanded) diagram commutes
       $$
         \xymatrix{
         &&{\bcal U}(B,C)(g,g)\otimes_2 {\bcal U}(A,B)(f,f)
         \ar[dd]^{{\cal M}_{ABC_{(g,f)(g,f)}}}\\
         I
         \ar[rru]^{j_{g}\otimes_2 j_{f}}
         \ar[rrd]_{j_{gf}}\\
         &&{\bcal U}(A,C)(gf,gf)
         }
        $$
     The heuristic diagram here is 
    $$
    \xymatrix@R-=3pt{
    &\ar@{=>}[dd]^{1_f}&&\ar@{=>}[dd]^{1_g}\\
    A
    \ar@/^1pc/[rr]^f
    \ar@/_1pc/[rr]_f
    &&B
    \ar@/^1pc/[rr]^g
    \ar@/_1pc/[rr]_g
    &&C
    \\
    &&&\\
    }
    =
    \xymatrix@R-=3pt{
    &\ar@{=>}[dd]^{1_{gf}}&&\\
    A
    \ar@/^1pc/[rr]^{gf}
    \ar@/_1pc/[rr]_{gf}
    &&C
    \\
    &&&\\
    } 
    $$   
       
       In addition, the ${\cal V}$--functoriality of ${\cal J}$ implies that the following (expanded) diagram commutes
       $$
         \xymatrix{
         &&{\cal I}(0,0)
         \ar[dd]^{{\cal J}_{A_{00}}}\\
         I
         \ar[rru]^{j_{0}}
         \ar[rrd]_{j_{1_A}}\\
         &&{\bcal U}(A,A)(1_A,1_A)
         }
       $$
       Which means that 
       $${\cal J}_{A_{00}}: I \to {\bcal U}(A,A)(1_{A},1_{A}) = j_{1_{A}}.$$

  %{3cat}

  % {3cat}
    \clearpage
          \newpage
        \MySection{The 3--category of enriched 2--categories}
  
  As in \cite{forcey1}, I now describe the (strict) 3--category  ${\cal V}$--2--Cat (or ${\cal V}$--Cat--Cat)
       whose objects are (strict, small)  ${\cal V}$--2--categories.
         
         \begin{definition} \label{2:enriched:funct}
         For two ${\cal V}$--2--categories ${\bcal U}$ and ${\bcal W}$ 
         a ${\cal V}${\it --2--functor} $T:{\bcal U} \to {\bcal W}$ is a function on objects 
         $\left|{\bcal U}\right| \to \left|{\bcal W}\right|$ and a family of ${\cal V}$--functors
         $T_{UU'}:{\bcal U}(U,U') \to {\bcal W}(TU,TU').$ These latter obey commutativity of the usual diagrams.
         \begin{enumerate}
          \item For $U,U',U'' \in \left|{\bcal U}\right|$
          $$
            \xymatrix@C=65pt@R=65pt{
            &\bullet
            \ar[rr]^{{\cal M}_{UU'U''}}
            \ar[d]^{T_{U'U''} \otimes^{(1)}_1 T_{UU'}}
            &&\bullet
            \ar[d]^{T_{UU''}}&\\
            &\bullet
            \ar[rr]_{{\cal M}_{(TU)(TU')(TU'')}}
            &&\bullet
            }
          $$
          \item
          $$
            \xymatrix{
            &&\bullet
            \ar[dd]^{T_{UU}}\\
            {\cal I}
            \ar[rru]^{{\cal J}_U}
            \ar[rrd]_{{\cal J}_{TU}}\\
            &&\bullet
            }
          $$
          \end{enumerate}
         \end{definition}
          For objects this means that $T_{U'U''}(f)T_{UU'}(g) = T_{UU''}(fg)$ and $T_{UU}(1_U) = 1_{TU}.$
          The reader should unpack both diagrams into terms of hom--object morphisms and ${\cal V}$--functoriality. 
          The fact that the hom--object morphisms are actually hom--category ${\cal V}$--functors corresponds to the
          need for ${\cal V}$--2--functors to preserve all the structure that exists, including the vertical composition.

         \begin{definition} \label{2:enr:nat:trans}
         A ${\cal V}${\it --2--natural transformation} $\alpha:T\to S:{\bcal U}\to {\bcal W}$ is a function sending
         each $U \in  \left|{\bcal U}\right|$ to a 
         ${\cal V}$--functor $\alpha_{U}: {\cal I}\to {\bcal W}(TU,SU)$ in such a way that we have commutativity of
         
         $$
           \xymatrix{
           &{\cal I} \otimes^{(1)}_1 {\bcal U}(U,U')
           \ar[rr]^-{ \alpha_{U'} \otimes^{(1)}_1 T_{UU'}}
           &&{\bcal W}(TU',SU') \otimes^{(1)}_1 {\bcal W}(TU,TU')
           \ar[rd]^-{{\cal M}}
         \\
           {\bcal U}(U,U')
           \ar[ru]^{=}
         \ar[rd]_{=}
           &&&&{\bcal W}(TU,SU')
         \\
           &{\bcal U}(U,U') \otimes^{(1)}_1 {\cal I}
           \ar[rr]_-{S_{UU'} \otimes^{(1)}_1 \alpha_{U}}
           &&{\bcal W}(SU,SU') \otimes^{(1)}_1 {\bcal W}(TU,SU)
           \ar[ru]^-{{\cal M}}
           }
         $$
         \end{definition}
        Unpacking this a bit, we see that $\alpha_U$ is an object $q = \alpha_U(0)$ in 
        the ${\cal V}$--category ${\bcal W}(TU,SU)$ and 
        a morphism $\alpha_{U_{00}}: I \to {\bcal W}(TU,SU)(q,q).$ By the ${\cal V}$--functoriality of $\alpha_U$
        we see that $\alpha_{U_{00}} = j_q.$ The axiom then states that $q'T_{UU'}(f) = S_{UU'}(f)q$ for all $f$,
        and that 
        $${\cal M}_{(TU)(TU')(SU')_{(q',T_{UU'}(f))(q',T_{UU'}(g))}} \circ (j_{q'} \otimes_2 T_{UU'_{fg}})
        ={\cal M}_{(TU)(SU)(SU')_{(S_{UU'}(f),q)(S_{UU'}(g),q)}} \circ (S_{UU'_{fg}}\otimes_2 j_q )$$
        
       The heuristic picture (following the pattern set in the definition of a ${\cal V}$--2--category) is as follows:
       $$
          \xymatrix@R-=3pt{
          &
          \\
           {TU}
          \ar@/^1pc/[rr]^{Tf}
          \ar@/_1pc/[rr]_{Tg}
          &&{TU'}
          \ar[rr]^{q'}
          &&{SU'}
          \\
         &\\
         }
         =
         \xymatrix@R-=3pt{
         &&&\\
         {TU}
         \ar[rr]^{q}
         &&{SU}
         \ar@/^1pc/[rr]^{Sf}
         \ar@/_1pc/[rr]_{Sg}
         &&{SU'}
         \\
         &&&\\
         }
    $$

        %Left and right whiskering of ${\cal V}$--2--functors onto ${\cal V}$--2--natural transformations are given by precisely the same 
        %descriptions as in the low dimensional case, with $I$ replaced by ${\cal I}$, etc.

           \begin{definition} \label{modification}
           Given two ${\cal V}$--2--natural transformations a 
         ${\cal V}${\it --modification} between them $\mu:\theta\to \phi :T \to S:{\bcal U}\to {\bcal W}$
         is a function that sends each object $U \in \left|{\bcal U}\right|$ to a morphism
         $\mu_U:I \to {\bcal W}(TU,SU)(\theta_U(0),\phi_U(0))$ in such a way that the following diagram commutes.
         (Let $\theta_U(0) = q$, $\phi_U(0) = \hat{q}$, $\theta_{U'}(0) = q'$ and $\phi_{U'}(0) = \hat{q'}.$)
         
         \noindent
	 		          \begin{center}
	          \resizebox{6.5in}{!}{ %!!!!
	          %%%%\resizebox{7in}{!}{
          $$
             \xymatrix@C=-65pt@R=25pt{
             &
             &{\bcal W}(TU',SU')(q',\hat{q'}) \otimes_2 {\bcal W}(TU,TU')(T_{UU'}(f),T_{UU'}(g))
             \ar[rd]^-{{\cal M}}
           \\
           &
           I \otimes_2 {\bcal U}(U,U')(f,g)
             \ar[ru]_{\text{ }\mu_{U'} \otimes_2 T_{UU'_{fg}}}
           &&{\bcal W}(TU,SU')(q'T_{UU'}(f),\hat{q'}T_{UU'}(g))
           \ar@{=}@/^2pc/[dd]
           \\
           {\bcal U}(U,U')(f,g)
             \ar[ru]^{=}
           \ar[rd]_{=}
           \\
           &
           {\bcal U}(U,U')(f,g) \otimes_2 I
             \ar[rd]^{\text{ }S_{UU'_{fg}} \otimes_2 \mu_{U}}
           &&{\bcal W}(TU,SU')(S_{UU'}(f)q,S_{UU'}(g)\hat{q})
           \\
             &
             &{\bcal W}(SU,SU')(S_{UU'}(f),S_{UU'}(g)) \otimes_2 {\bcal W}(TU,SU)(q,\hat{q})
             \ar[ru]^-{{\cal M}}
             }
           $$
           }
	   	                    \end{center}
	       \end{definition}             
         Notice that since $\theta_{U_{00}} = j_{\theta_U(0)}$ for all ${\cal V}$--2--natural transformations $\theta$ we have 
         that the morphism $\mu_U$ seen as a ``family'' consisting of a single morphism 
         (corresponding to $0 \in \left|{\cal I}\right|$)
         constitutes a ${\cal V}$--natural transformation from $\theta_U$ to $\phi_U.$ Occasionally I
         reflect this by denoting $\mu_U$ as ${\mu_U}_0.$
           
         The heuristic picture here is:
         $$
            \xymatrix@R-=3pt{
            &&&\\
            {TU}
            \ar@/^1pc/[rr]^{Tf}
            \ar@/_1pc/[rr]_{Tg}
            &&{TU'}
            \ar@/^1pc/[rr]^{q'}
            \ar@/_1pc/[rr]_{\hat{q'}}
            &&{SU'}
            \\
            &&&\\
            }
           =
           \xymatrix@R-=3pt{
           &&&\\
           {TU}
           \ar@/^1pc/[rr]^{q}
           \ar@/_1pc/[rr]_{\hat{q}}
           &&{SU}
           \ar@/^1pc/[rr]^{Sf}
           \ar@/_1pc/[rr]_{Sg}
           &&{SU'}
           \\
           &&&\\
           }
    $$

     \begin{theorem}
     ${\cal V}$--2--categories, ${\cal V}$--2--functors, ${\cal V}$--2--natural transformations
     and ${\cal V}$--modifications form a 3--category called ${\cal V}$--2--Cat. 
     \end{theorem}
     
     \begin{proof}\text{ }(Part 1.)
     Recall that a 3--category is a category enriched over 2--Cat. This is expanded in terms of axioms in \cite{Borc}.
     Our objects are ${\cal V}$--2--categories. There are two parts of the proof. In part 1 we show that for every pair ${\bcal U}, {\bcal W}$ of 
     ${\cal V}$--2--categories we have a 2--category made up of ${\cal V}$--2--functors, ${\cal V}$--2--natural transformations
     and ${\cal V}$--modifications. For now then  ${\cal V}$--2--functors are the 0--cells, ${\cal V}$--2--natural transformations
     are the 1--cells, and ${\cal V}$--modifications are the 2--cells as in the following picture.
     
     $$
        \xymatrix@R-=5pt{
        {\bcal U}
        \ar[dd]_{F}
        &\ar@{=>}[dd]^{\mu}
        &{\bcal U}
        \ar[dd]^{G}
        \\
         {}
        \ar@/^2pc/[rr]^{\gamma}
        \ar@/_2pc/[rr]_{\rho}
        &&{}
        \\
       {\bcal W}&&{\bcal W}\\
    }
    $$
    %\clearpage
    Throughout I will use the following notation:
    Composition along a ${\cal V}$--2--natural transformation will be indicated with ``$\circ.$'' Composition
    along a ${\cal V}$--2--functor will be indicated with ``$*.$'' Composition along a ${\cal V}$--2--category
    will be indicated by juxtaposition.
     
      Composition of ${\cal V}$--2--natural transformations $\gamma:T\to S$ and $\beta:S\to R$ 
      along a ${\cal V}$--2--functor $S$ is given by
           
           $$
           \xymatrix{(\beta * \gamma)_U} =
           \xymatrix{
           {\cal I} = {\cal I} \otimes^{(1)}_1 {\cal I}
           \ar[d]^{\beta_U \otimes^{(1)}_1 \gamma_U}
           \\
           {\bcal W}(SU,RU) \otimes^{(1)}_1 {\bcal W}(TU,SU)
           \ar[d]^{{\cal M}}
           \\
           {\bcal W}(TU,RU)
           }
           $$
           
      Let $\beta_U(0) = \hat{q}$ and $\gamma_U(0) = \check{q}.$
      By expanding the definition of this composition we see that $(\beta * \gamma)_U(0) = \hat{q}\check{q}$ 
      and that $(\beta * \gamma)_{U_{00}} = j_{\hat{q}\check{q}}$
      
      Since it prefigures a similar proof for ${\cal V}$--modifications,
      I include the proof that this composition forms a valid ${\cal V}$--2--natural transformation even though it follows
      closely the analogous proof for ${\cal V}$--natural transformations as in \cite{EK1}.
      For $\beta * \gamma$ to be a ${\cal V}$--2--natural transformation the exterior of the following diagram must commute.
      
      \noindent
            		          \begin{center}
	          \resizebox{6in}{!}{ %!!!!
	          %%%%\resizebox{6.5in}{!}{
      \begin{sideways}
      %\begin{footnotesize}
      $$
      \xymatrix
      @R=21pt@C=-15pt
      {
      &&&&({\bcal W}(SU',RU') \otimes^{(1)}_1 {\bcal W}(TU',SU'))\otimes^{(1)}_1 {\bcal W}(TU,TU')
      \ar[ddrr]|{{\cal M}\otimes^{(1)}_1 1}
      \ar[ddd]^{\alpha^{(1)1}}
      \\\\
      &&({\cal I} \otimes^{(1)}_1 {\cal I}) \otimes^{(1)}_1 {\bcal U}(U,U')
      \ar[uurr]|{(\beta_{U'} \otimes^{(1)}_1 \gamma_{U'})\otimes^{(1)}_1 T_{UU'}}
      &&&&{\bcal W}(TU',RU')\otimes^{(1)}_1 {\bcal W}(TU,TU')
      \ar[ddddddr]^{{\cal M}}
      \\
      &&&
      \bullet 
          %{\bcal W}(SU',RU') \otimes^{(1)}_1 ({\cal I} \otimes^{(1)}_1 {\bcal U}(U,U'))
      \ar[r]^{1 \otimes^{(1)}_1 (\gamma_{U'}\otimes^{(1)}_1 T_{UU'})}
      &
      \bullet 
          %{\bcal W}(SU',RU') \otimes^{(1)}_1 ({\bcal W}(TU',SU')\otimes^{(1)}_1 {\bcal W}(TU,TU'))
      \ar[ddr]|{1\otimes^{(1)}_1 {\cal M}}
      \\
      &{\cal I} \otimes^{(1)}_1 {\bcal U}(U,U')
      \ar[uur]^{=}
      \ar[dr]|{\beta_{U'} \otimes^{(1)}_1 1}
      \\
      &&\bullet
      %{\bcal W}(SU',RU') \otimes^{(1)}_1 {\bcal U}(U,U')
      \ar[uur]^{=}
      \ar[ddr]^{=}
      &&&{\bcal W}(SU',RU') \otimes^{(1)}_1 {\bcal W}(TU,SU')
      \ar[dddrr]|{{\cal M}}
      \\\\
      &&&\bullet 
      %{\bcal W}(SU',RU') \otimes^{(1)}_1 ({\bcal U}(U,U') \otimes^{(1)}_1 {\cal I})
      \ar[r]^{1\otimes^{(1)}_1 (S_{UU'}\otimes^{(1)}_1 \gamma_U)}
      &\bullet 
      %{\bcal W}(SU',RU') \otimes^{(1)}_1 ({\bcal W}(SU,SU')\otimes^{(1)}_1 {\bcal W}(TU,SU))
      \ar[uur]|{1\otimes^{(1)}_1 {\cal M}}
      \\
      {\bcal U}(U,U')
      \ar[uuuur]^{=}
      \ar[ddddr]^{=}
      &&&&&&&{\bcal W}(TU,RU')
      \\
      &&&\bullet 
      \ar[r]^{(\beta_{U'}\otimes^{(1)}_1 S_{UU'})\otimes^{(1)}_1 1}
      &\bullet
      \ar[ddr]|{{\cal M}\otimes^{(1)}_1 1}
      \ar[uu]_{\alpha^{(1)1}}
      \\\\
      &&\bullet
      \ar[uur]^{=}
      \ar[ddr]^{=}
      &&&{\bcal W}(SU,RU') \otimes^{(1)}_1 {\bcal W}(TU,SU)
      \ar[uuurr]|{{\cal M}}
      \\
      &{\bcal U}(U,U') \otimes^{(1)}_1 {\cal I}
      \ar[ur]|{1\otimes^{(1)}_1 \gamma_U}
      \ar[ddr]^{=}
      \\
      &&&\bullet
      \ar[r]^{(R_{UU'}\otimes^{(1)}_1 \beta_{U})\otimes^{(1)}_1 1}
      &\bullet
      \ar[uur]|{{\cal M}\otimes^{(1)}_1 1}
      \ar[ddd]^{\alpha^{(1)1}}
      \\
      &&{\bcal U}(U,U')\otimes^{(1)}_1 ({\cal I} \otimes^{(1)}_1 {\cal I})
      \ar[ddrr]|{R_{UU'}\otimes^{(1)}_1 (\beta_{U}\otimes^{(1)}_1 \gamma_U)}
      &&&&{\bcal W}(RU,RU')\otimes^{(1)}_1 {\bcal W}(TU,RU)
      \ar[uuuuuur]_{{\cal M}}
      \\\\
      &&&&{\bcal W}(RU,RU')\otimes^{(1)}_1 ({\bcal W}(SU,RU)\otimes^{(1)}_1 {\bcal W}(TU,SU))
      \ar[uurr]|{1\otimes^{(1)}_1 {\cal M}}
      }
      $$
      %\end{footnotesize}
      
      \end{sideways}
       }
      	                    \end{center}
	                    
	                    %\clearpage
       %\newpage
       The arrows marked with an ``='' all occur as copies of ${\cal I}$ are tensored to the object at the arrow's source. 
       The 3 leftmost regions commute by the naturality of $\alpha^{(1)1}=\alpha^2.$ The 2 embedded central hexagons commute by the
       definition of ${\cal V}$--2--natural transformations for $\gamma$ and $\beta.$ The three pentagons on the right
       are copies of the pentagon axiom for the composition ${\cal M}.$ The associativity of this composition also follows directly from
       the latter axiom. 
       
         The identities for this composition are ${\cal V}$--2--natural transformations $\mathbf{1}_T:T\to T$ where 
        $(\mathbf{1}_T)_U = {\cal J}_{TU}.$ That this describes a 2-sided identity for the composition above is easily checked using the 
        unit axioms for a ${\cal V}$--2--category.

                 The composition of two 
                        ${\cal V}$--modifications along a ${\cal V}$--2--transformation is given by
                        the composition of the underlying ${\cal V}$--natural
           transformations. So given ${\cal V}$--2--natural
                        transformations $\alpha, \beta$ and $\sigma:F \to G:{\bcal
           U}\to {\bcal W}$, and ${\cal V}$--modifications
                        $\mu:\alpha \to \beta$ and $\nu:\beta \to \sigma$ as in the following picture

                        $$
  		        
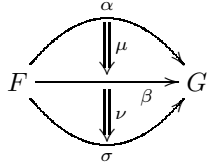
\xymatrix@R-=16pt{
  		        &\ar@{=>}[d]^{\mu}\\
  		        F
  		        \ar@/^2pc/[rr]^{\alpha}
  		        \ar[rr]_>>>>>{\beta}
  		        \ar@/_2pc/[rr]_{\sigma}
  		        &\ar@{=>}[d]^{\nu}
  		        &G\\
  		        &\\
  		        }
    $$
    where
           $\alpha_U(0)=q,$ $\beta_U(0) = \hat{q}$ and 
                        $\sigma_U(0) = \check{q},$ we have 
     $$\xymatrix{(\nu \circ \mu)_U} =
                    \xymatrix{
                    I = I \otimes_1 I
                    \ar[d]^{\nu_U \otimes_1 \mu_U}
                    \\
                    {\bcal W}(FU,GU)(\hat{q},\check{q}) \otimes_1 {\bcal
           W}(FU,GU)(q,\hat{q})
                    \ar[d]^{M}
                    \\
                    {\bcal W}(FU,GU)(q,\check{q})
                    }
                    $$
                    
           We see that this composition is associative by the associativity pentagon
           for $M.$ We also see that the result of a 
           composition is a ${\cal V}$--natural transformation as well. It needs to
           be checked that the result of a composition
           is a valid ${\cal V}$--modification. This is seen by showing that the
           exterior of the following diagram commutes.

  The first bullet in the following diagram is ${\bcal U}(U,U')(f,g).$
  Other objects include:
  $$A= (I\otimes_1 I) \otimes_2 {\bcal U}(U,U')(f,g)\text{ }\text{ ; }  B= (I\otimes_1 I) \otimes_2 ({\bcal U}(U,U')(f,g)\otimes_1 I)\text{ ; }$$
  $$C= (I\otimes_1 I) \otimes_2 ({\bcal U}(U,U')(f,g)\otimes_1 {\bcal U}(U,U')(f,f))\text{ }\text{ ; }  D= (I\otimes_2 {\bcal U}(U,U')(f,g)) \otimes_1 (I\otimes_2 {\bcal U}(U,U')(f,f))\text{ ; }$$
  $$E= {\bcal U}(U,U')(f,g)\otimes_1 I\text{ }\text{ ; }  H= {\bcal U}(U,U')(f,g)\otimes_1 {\bcal U}(U,U')(f,f)\text{ ; }$$
  $$K= {\bcal U}(U,U')(f,g)\otimes_2 (I\otimes_1 I)\text{ }\text{ ; }  L= ({\bcal U}(U,U')(f,g)\otimes_1 I)\otimes_2 (I\otimes_1 I)\text{ ; }$$
  $$N= ({\bcal U}(U,U')(f,g)\otimes_2 I) \otimes_1 ({\bcal U}(U,U')(f,f)\otimes_2 I)\text{ }\text{ ; }$$
  $$P= ({\bcal U}(U,U')(f,g)\otimes_2 {\bcal U}(U,U')(f,f)) \otimes_1 ({\bcal U}(U,U')(f,f)\otimes_2 I).$$
  \clearpage
   \noindent
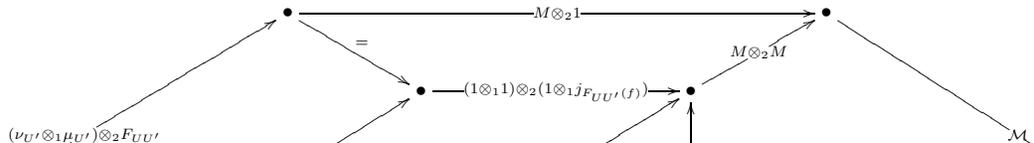
   \begin{center}
	          \resizebox{6in}{!}{%!!!!
	          %%%%\resizebox{6.5in}{!}{
   \begin{sideways}
    $$
    \xymatrix@C=49pt@R=26pt{
    &&&&\bullet
    \ar[dr]^{=}
    \ar[rrrr]|{M\otimes_2 1}
    &&&&\bullet
    \ar[dddrrr]|{{\cal M}}
    \\
    &&&&&\bullet
    \ar[rr]|{(1\otimes_1 1)\otimes_2 (1\otimes_1 j_{F_{UU'}(f)})}
    &&\bullet
    \ar[ur]|{M\otimes_2 M}
    \\
    &&&\text{\scriptsize {[a]}}
    &&\text{\scriptsize {[b]}}
    &&&\text{\scriptsize {[d]}}
    \\
    &A
    \ar[uuurrr]|{(\nu_{U'}\otimes_1 \mu_{U'})\otimes_2 F_{UU'}}
    \ar[rr]^{=}
    &&B
    \ar[uurr]|{(\nu_{U'}\otimes_1 \mu_{U'})\otimes_2 (F_{UU'}\otimes_1 1)}
    \ar[rr]|{(1\otimes_1 1)\otimes_2 (1\otimes_1 j_f)}
    &&C
    \ar[uurr]|{(\nu_{U'}\otimes_1 \mu_{U'})\otimes_2 (F_{UU'}\otimes_1 F_{UU'})}
    &\text{\scriptsize {[c]}}
    &
    &&&&\bullet
    \ar@{=}@/^3pc/[dddddddd]
    \\
    &&&&&&&\bullet
    \ar[uuu]_{\eta^{1,2}}
    \ar[drr]|{{\cal M}\otimes_1 {\cal M}}
    \\
    &\text{\scriptsize {[e]}}
    &&\text{\scriptsize {[f]}}
    &&D
    \ar[uu]^{\eta^{1,2}}
    \ar[urr]|{(\nu_{U'}\otimes_2 F_{UU'})\otimes_1 ( \mu_{U'}\otimes_2 F_{UU'})}
    &&&&\bullet
    \ar[uurr]|{M}
    \ar@{=}@/^2pc/[dddd]
    \\\\
    \bullet
    %{\bcal U}(U,U')(f,g)
    \ar[uuuur]^{=}
    \ar[rr]^{=}
    \ar[ddddr]^{=}
    &&E
    \ar[uuuur]^{=}
    \ar[r]|{1\otimes_1 j_f}
    \ar[ddddr]^{=}
    &H
    \ar[uurr]^{=}
    \ar[ddrr]^{=}
    \ar[uuuurr]^{=}
    \ar[ddddrr]^{=}
    &&&\text{\scriptsize {[g]}}
    \\\\
    &\text{\scriptsize {[h]}}
    &&\text{\scriptsize {[i]}}
    &&N
    \ar[dd]_{\eta^{1,2}}
    \ar[drr]|{(G_{UU'}\otimes_2 \nu_{U})\otimes_1 (G_{UU'}\otimes_2 \mu_{U})}
    &&&&\bullet
    \ar[ddrr]|{M}
    \\
    &&&&&&&\bullet
      \ar[ddd]^{\eta^{1,2}}
    \ar[urr]|{{\cal M}\otimes_1 {\cal M}}
    \\
    &K
    \ar[dddrrr]|{G_{UU'}\otimes_2 (\nu_U\otimes_1 \mu_U)}
    \ar[rr]^{=}
    &&L
    \ar[ddrr]|{(G_{UU'}\otimes_1 1)\otimes_2 (\nu_U\otimes_1 \mu_U)}
    \ar[rr]|{(1\otimes_1 j_f)\otimes_2 (1\otimes_1 1)}
    &&P
    \ar[ddrr]|{(G_{UU'}\otimes_1 G_{UU'})\otimes_2 (\nu_U\otimes_1 \mu_U)}
    &\text{\scriptsize {[l]}}
    &
    &&&&\bullet
    \\
    &&&\text{\scriptsize {[j]}}
    &&\text{\scriptsize {[k]}}
    &&&\text{\scriptsize {[m]}}
    \\
    &&&&&\bullet
    \ar[rr]|{(1\otimes_1 j_{G_{UU'}(f)})\otimes_2 (1\otimes_1 1)}
    &&\bullet
    \ar[dr]|{M\otimes_2 M}
    \\
    &&&&\bullet
    \ar[ur]^{=}
    \ar[rrrr]|{1\otimes_2 M}
    &&&&\bullet
    \ar[uuurrr]|{{\cal M}}
    }
    $$
    \end{sideways}
  }
    	                    \end{center}

  \clearpage
  %\newpage
  The arrows marked with an ``='' all occur as copies of $I$ are tensored to the object at the arrow's source. Therefore
  the quadrilateral regions [a],[e],[f],[h],[i] and [j] all commute trivially. The uppermost and lowermost quadrilaterals
  commute by the property of composing with units in an enriched category. The two triangles commute by the external
  unit condition for iterated monoidal categories. Regions [b] and [k] commute by respect of units by enriched functors.
  Regions [c] and [l] commute by naturality of $\eta.$ Region [g] commutes by the definition of ${\cal V}$--modification
  for $\nu$ and $\mu.$ Regions [d] and [m] commute by the ${\cal V}$--functoriality of ${\cal M}.$ The following 
  heuristic
  diagram for this proof is quite instructive. (See the pattern set for these diagrams in the definition 
  of a ${\cal V}$--2--category.)
  $$
    \xymatrix@R-=16pt{
    &\ar@{=>}[d]^{1_{Ff}}
    &&
    \\
    FU
    \ar@/^2pc/[rr]^{Ff}
    \ar@{..>}[rr]_{Ff}
    \ar@/_2pc/[rr]_{Fg}
    &
    &FU'
    \ar@/^2pc/[rr]^{q'}
    \ar[rr]^{\hat{q'}}
    \ar@/_2pc/[rr]^{\check{q'}}
    &
    &GU'\\
    &&&&&
    }
    =
     \xymatrix@R-=16pt{
      &
      &&\ar@{=>}[d]^{1_{Gf}}
      \\
      FU
      \ar@/^2pc/[rr]^{q}
      \ar[rr]^{\hat{q}}
      \ar@/_2pc/[rr]^{\check{q}}
      &
      &GU
      \ar@/^2pc/[rr]^{Gf}
      \ar@{..>}[rr]_{Gf}
      \ar@/_2pc/[rr]_{Gg}
      &
      &GU'\\
      &&&&&
      }
    $$

               Thus identities  $\mathbf{1}_{\alpha}$ for this
               composition are families of ${\cal V}$--natural equivalences. Since $\alpha_U$ is a 
         ${\cal V}$--functor from ${\cal I}$ to ${\bcal W}(TU,SU)$ this means 
         specifically that $((\mathbf{1}_{\alpha})_U)_0 = j_{\alpha_{U}(0)} = j_q.$ Recall that here the ``family''
         has only one member, corresponding to the single object in ${\cal I}.$ That this describes a 2-sided identity for the composition above is easily checked using the 
        unit axioms for a ${\cal V}$--category.
  
  In order to define composition of all allowable pasting diagrams in the 2-category, we need only to define the 
    composition described by left and right whiskering diagrams (as partial functors) 
    and check that these can be combined into a
    well--defined horizontal composition. The first 
    picture shows a 1-cell (that is a ${\cal V}$--2--natural transformation between ${\cal V}$--2--functors 
    $F,G:{\bcal U}\to {\bcal W}$) following a 2-cell (a ${\cal V}$--modification).
    These are composed to form a new 2-cell as follows
    
    $$
    \xymatrix@R-=3pt{
    &\ar@{=>}[dd]^{\mu}\\
    F
    \ar@/^1pc/[rr]^{\psi}
    \ar@/_1pc/[rr]_{\beta}
    &&G
    \ar[rr]^{\gamma}
    &&H
    \\
    &\\
    } \text{   is composed to become } \xymatrix@R-=3pt{
    &\ar@{=>}[dd]^{\gamma *\mu}\\
    F
    \ar@/^1pc/[rr]^{\gamma *\psi}
    \ar@/_1pc/[rr]_{\gamma *\beta}
    &&H
    \\
    &\\
    }
    $$
     where $\gamma *\psi$ and $\gamma *\beta$ are described above, and $\gamma *\mu$
    has components given by the following composition: (Let $\psi_U(0) = q$,
    $\beta_U(0) = \hat{q}$ and 
    $\gamma_U(0) = \check{q}.$ 
    % Let $\psi_{U'}(0) = q'$, $\beta_{U'}(0) = \hat{q'}$ and $\gamma_{U'}(0) = \check{q'}.$ 
    Note that $j_{\check{q}} = \gamma_{U_{00}}.$)
    $$ \xymatrix{(\gamma *\mu)_{U}} = \xymatrix{
                    I = I \otimes_2 I
                    \ar[d]^{j_{\check{q}} \otimes_2 \mu_U}
                    \\
                    {\bcal W}(GU,HU)(\check{q},\check{q}) \otimes_2 {\bcal W}(FU,GU)(q,\hat{q})
                    \ar[d]^{{\cal M}}
                    \\
                    {\bcal W}(FU,HU)(\check{q}q,\check{q}\hat{q})
                    }
                    $$
                    \noindent
      \begin{center}
        \begin{list}{}{}            
        For this composition to yield a valid ${\cal V}$--modification 
        the exterior of the following diagram must commute.
    
      \noindent
      \begin{center}
	          \resizebox{6in}{!}{ %!!!!
	          %%%%\resizebox{6.5in}{!}{
	      
    \begin{sideways}\nopagebreak    
            %\begin{footnotesize}
        $$
        \xymatrix
        @R=27pt@C=-64pt
        {
        &&&&({\bcal W}(GU',HU')(\check{q'},\check{q'}) \otimes_2 {\bcal W}(FU',GU')(q',\hat{q'}))\otimes_2 {\bcal W}(FU,FU')(Ff,Fg)
        \ar[ddrr]|{{\cal M}\otimes_2 1}
        \ar[ddd]^{\alpha^2}
        \\\\
        &&(I \otimes_2 I) \otimes_2 {\bcal U}(U,U')(f,g)
        \ar[uurr]|{(j_{\check{q'}} \otimes_2 \mu_{U'})\otimes_2 F_{UU'}}
        &&&&{\bcal W}(FU',HU')(\check{q'}q',\check{q'}\hat{q'})\otimes_2 {\bcal W}(FU,FU')(Ff,Fg)
        \ar[ddddr]^{{\cal M}}
        \\
        &&&
        \bullet 
            %{\bcal W}(GU',HU')(\check{q'},\check{q'}) \otimes_2 (I \otimes_2 {\bcal U}(U,U')(f,g))
        \ar[r]^{1 \otimes_2 (\mu_{U'}\otimes_2 F_{UU'})}
        &
        \bullet 
            %{\bcal W}(GU',HU')(\check{q'},\check{q'}) \otimes_2 ({\bcal W}(FU',GU')(q',\hat{q'})\otimes_2 {\bcal W}(FU,FU')(Ff,Fg))
        \ar[dr]|{1\otimes_2 {\cal M}}
        \\
        &I \otimes_2 {\bcal U}(U,U')(f,g)
        \ar[uur]^{=}
        \ar[dr]|{j_{\check{q'}} \otimes_2 1}
        &&&&{\bcal W}(GU',HU')(\check{q'},\check{q'}) \otimes_2 {\bcal W}(FU,GU')(q'Ff,\hat{q'}Fg)
        \ar[ddrr]|{{\cal M}}
        \ar@{=}[dd]
        \\
        &&\bullet
        %{\bcal W}(GU',HU')(\check{q'},\check{q'}) \otimes_2 {\bcal U}(U,U')(f,g)
        \ar@/_1pc/[uur]_{=}
        \ar@/^1pc/[ddr]^{=}
        &&&
        \\
        &&&&&{\bcal W}(GU',HU')(\check{q'},\check{q'}) \otimes_2 {\bcal W}(FU,GU')(Gfq,Gg\hat{q})
        \ar[ddrr]|{{\cal M}}
        &&{\bcal W}(FU,HU')(\check{q'}q'Ff,\check{q'}\hat{q'}Fg) 
        \ar@{=}[dd]
        \\
        &&&\bullet 
        %{\bcal W}(GU',HU')(\check{q'},\check{q'}) \otimes_2 ({\bcal U}(U,U')(f,g) \otimes_2 I)
        \ar[r]^{1\otimes_2 (G_{UU'}\otimes_2 \mu_U)}
        &\bullet 
        %{\bcal W}(GU',HU')(\check{q'},\check{q'}) \otimes_2 ({\bcal W}(GU,GU')(Gf,Gg)\otimes_2 {\bcal W}(FU,GU)(q,\hat{q}))
        \ar[ur]|{1\otimes_2 {\cal M}}
        \\
        {\bcal U}(U,U')(f,g)
        \ar[uuuur]^{=}
        \ar[ddddr]_{=}
        &&&&&&&{\bcal W}(FU,HU')(\check{q'}Gfq,\check{q'}Gg\hat{q})
        \\
        &&&\bullet 
        \ar[r]^{(j_{\check{q'}}\otimes_2 G_{UU'})\otimes_2 1}
        &\bullet
        \ar[dr]|{{\cal M}\otimes_2 1}
        \ar[uu]_{\alpha^2}
        \\
        &&&&&{\bcal W}(GU,HU')(\check{q'}Gf,\check{q'}Gg) \otimes_2 {\bcal W}(FU,GU)(q,\hat{q})
        \ar[uurr]|{{\cal M}}
        \ar@{=}[dd]
        &&{\bcal W}(FU,HU')(Hf\check{q}q,Hg\check{q}\hat{q}) 
        \ar@{=}[uu]
        \\
        &&\bullet
        \ar@/_1pc/[uur]_{=}
        \ar@/^1pc/[ddr]^{=}
        &&&
        \\
        &{\bcal U}(U,U')(f,g) \otimes_2 I
        \ar[ur]|{1\otimes_2 \mu_U}
        \ar[ddr]_{=}
        &&&&{\bcal W}(GU,HU')(Hf\check{q},Hg\check{q}) \otimes_2 {\bcal W}(FU,GU)(q,\hat{q})
        \ar[uurr]|{{\cal M}}
        \\
        &&&\bullet
        \ar[r]^{(H_{UU'}\otimes_2 j_{\check{q}})\otimes_2 1}
        &\bullet
        \ar[ur]|{{\cal M}\otimes_2 1}
        \ar[ddd]^{\alpha^2}
        \\
        &&{\bcal U}(U,U')(f,g)\otimes_2 (I \otimes_2 I)
        \ar[ddrr]|{H_{UU'}\otimes_2 (j_{\check{q}}\otimes_2 \mu_U)}
        &&&&{\bcal W}(HU,HU')(Hf,Hg)\otimes_2 {\bcal W}(FU,HU)(\check{q}q,\check{q}\hat{q})
        \ar[uuuur]_{{\cal M}}
        \\\\
        &&&&{\bcal W}(HU,HU')(Hf,Hg)\otimes_2 ({\bcal W}(GU,HU)(\check{q},\check{q})\otimes_2 {\bcal W}(FU,GU)(q,\hat{q}))
        \ar[uurr]|{1\otimes_2 {\cal M}}
        }
        $$\nopagebreak
        %\end{footnotesize}
         \end{sideways}
         }
		                    \end{center}
		                    \end{list}
		                    \end{center}
    \clearpage
    %\newpage
    The arrows marked with an ``='' all occur as copies of $I$ are tensored to the object at the arrow's source. 
         The 3 leftmost regions commute by the naturality of $\alpha^2.$ The 2 embedded central ``hexagons'' commute by the
         definition of ${\cal V}$--modifications for $\mu$ and $\mathbf{1}_{\gamma}.$ The three pentagons on the right
         are copies of the pentagon axiom for the composition ${\cal M}.$ 
         
    The second picture shows a 2--cell following a 1--cell. These are composed as follows
    $$
    \xymatrix@R-=3pt{
    &&&\ar@{=>}[dd]^{\mu}\\
    E
    \ar[rr]^{\rho}
    &&F
    \ar@/^1pc/[rr]^{\psi}
    \ar@/_1pc/[rr]_{\beta}
    &&G
    \\
    &&&\\
    } \text{   is composed to become } \xymatrix@R-=3pt{
    &\ar@{=>}[dd]^{\mu *\rho}\\
    E
    \ar@/^1pc/[rr]^{\psi *\rho}
    \ar@/_1pc/[rr]_{\beta *\rho}
    &&G
    \\
    &\\
    }
    $$
    
    where $\mu *\rho$ has components given by the following composition: (Let $\rho_U(0) = \overline{q}.$
    Note that $j_{\overline{q}} = \rho_{U_{00}}.$)
    $$\xymatrix{(\mu *\rho)_{U}} = \xymatrix{
                    I = I \otimes_2 I
                    \ar[d]^{ \mu_U\otimes_2 j_{\overline{q}}}
                    \\
                    {\bcal W}(FU,GU)(q,\hat{q}) \otimes_2 {\bcal W}(EU,FU)(\overline{q},\overline{q})
                    \ar[d]^{{\cal M}}
                    \\
                    {\bcal W}(FU,HU)(q\overline{q},\hat{q}\overline{q})
                    }
                    $$
                    
      For this composition to yield a valid ${\cal V}$--modification the exterior of the following diagram must commute.               
      
      \noindent
          		          \begin{center}
	          \resizebox{6in}{!}{ %!!!!
	          %%%%\resizebox{6.5in}{!}{
 \begin{sideways}
               %\begin{footnotesize}
          $$
          \xymatrix
          @R=27pt@C=-64pt
          {
          &&&&({\bcal W}(FU',GU')(q',\hat{q'}) \otimes_2 {\bcal W}(EU',FU')(\overline{q'},\overline{q'}))\otimes_2 {\bcal W}(EU,EU')(Ef,Eg)
          \ar[ddrr]|{{\cal M}\otimes_2 1}
          \ar[ddd]^{\alpha^2}
          \\\\
          &&(I \otimes_2 I) \otimes_2 {\bcal U}(U,U')(f,g)
          \ar[uurr]|{(\mu_{U'} \otimes_2 j_{\overline{q'}})\otimes_2 E_{UU'}}
          &&&&{\bcal W}(EU',GU')(q'\overline{q'},\hat{q'}\overline{q'})\otimes_2 {\bcal W}(EU,EU')(Ef,Eg)
          \ar[ddddr]^{{\cal M}}
          \\
          &&&
          \bullet 
              %{\bcal W}(FU',GU')(q',\hat{q'}) \otimes_2 (I \otimes_2 {\bcal U}(U,U')(f,g))
          \ar[r]^{1 \otimes_2 (j_{\overline{q'}}\otimes_2 E_{UU'})}
          &
          \bullet 
              %{\bcal W}(FU',GU')(q',\hat{q'}) \otimes_2 ({\bcal W}(EU',FU')(\overline{q'},\overline{q'})\otimes_2 {\bcal W}(EU,EU')(Ef,Eg))
          \ar[dr]|{1\otimes_2 {\cal M}}
          \\
          &I \otimes_2 {\bcal U}(U,U')(f,g)
          \ar[uur]^{=}
          \ar[dr]|{\mu_{U'} \otimes_2 1}
          &&&&{\bcal W}(FU',GU')(q',\hat{q'}) \otimes_2 {\bcal W}(EU,FU')(\overline{q'}Ef,\overline{q'}Eg)
          \ar[ddrr]|{{\cal M}}
          \ar@{=}[dd]
          \\
          &&\bullet
          %{\bcal W}(FU',GU')(q',\hat{q'}) \otimes_2 {\bcal U}(U,U')(f,g)
          \ar@/_1pc/[uur]_{=}
          \ar@/^1pc/[ddr]^{=}
          &&&
          \\
          &&&&&{\bcal W}(FU',GU')(q',\hat{q'}) \otimes_2 {\bcal W}(EU,FU')(Ff\overline{q},Fg\overline{q})
          \ar[ddrr]|{{\cal M}}
          &&{\bcal W}(EU,GU')(q'\overline{q'}Ef,\hat{q'}\overline{q'}Eg) 
          \ar@{=}[dd]
          \\
          &&&\bullet 
          %{\bcal W}(FU',GU')(q',\hat{q'}) \otimes_2 ({\bcal U}(U,U')(f,g) \otimes_2 I)
          \ar[r]^{1\otimes_2 (F_{UU'}\otimes_2 j_{\overline{q}})}
          &\bullet 
          %{\bcal W}(FU',GU')(q',\hat{q'}) \otimes_2 ({\bcal W}(FU,FU')(Ff,Fg)\otimes_2 {\bcal W}(EU,FU)(\overline{q},\overline{q}))
          \ar[ur]|{1\otimes_2 {\cal M}}
          \\
          {\bcal U}(U,U')(f,g)
          \ar[uuuur]^{=}
          \ar[ddddr]_{=}
          &&&&&&&{\bcal W}(EU,GU')({q'}Ff\overline{q},\hat{q'}Fg\overline{q})
          \\
          &&&\bullet 
          \ar[r]^{(\mu_{U'}\otimes_2 F_{UU'})\otimes_2 1}
          &\bullet
          \ar[dr]|{{\cal M}\otimes_2 1}
          \ar[uu]_{\alpha^2}
          \\
          &&&&&{\bcal W}(FU,GU')({q'}Ff,\hat{q'}Fg) \otimes_2 {\bcal W}(EU,FU)(\overline{q},\overline{q})
          \ar[uurr]|{{\cal M}}
          \ar@{=}[dd]
          &&{\bcal W}(EU,GU')(Gfq\overline{q},Gg\hat{q}\overline{q}) 
          \ar@{=}[uu]
          \\
          &&\bullet
          \ar@/_1pc/[uur]_{=}
          \ar@/^1pc/[ddr]^{=}
          &&&
          \\
          &{\bcal U}(U,U')(f,g) \otimes_2 I
          \ar[ur]|{1\otimes_2 j_{\overline{q}}}
          \ar[ddr]_{=}
          &&&&{\bcal W}(FU,GU')(Gf{q},Gg\hat{q}) \otimes_2 {\bcal W}(EU,FU)(\overline{q},\overline{q})
          \ar[uurr]|{{\cal M}}
          \\
          &&&\bullet
          \ar[r]^{(G_{UU'}\otimes_2 \mu_{U})\otimes_2 1}
          &\bullet
          \ar[ur]|{{\cal M}\otimes_2 1}
          \ar[ddd]^{\alpha^2}
          \\
          &&{\bcal U}(U,U')(f,g)\otimes_2 (I \otimes_2 I)
          \ar[ddrr]|{G_{UU'}\otimes_2 (\mu_U \otimes_2 j_{\overline{q}})}
          &&&&{\bcal W}(GU,GU')(Gf,Gg)\otimes_2 {\bcal W}(EU,GU)(q\overline{q},\hat{q}\overline{q})
          \ar[uuuur]_{{\cal M}}
          \\\\
          &&&&{\bcal W}(GU,GU')(Gf,Gg)\otimes_2 ({\bcal W}(FU,GU)(q,\hat{q})\otimes_2 {\bcal W}(EU,FU)(\overline{q},\overline{q}))
          \ar[uurr]|{1\otimes_2 {\cal M}}
          }
          $$
          %\end{footnotesize}
          \end{sideways}
          }
  \end{center}          
  
      \clearpage
      %\newpage
      The arrows marked with an ``='' all occur as copies of $I$ are tensored to the object at the arrow's source. 
           The 3 leftmost regions commute by the naturality of $\alpha^2.$ The 2 embedded central ``hexagons'' commute by the
           definition of ${\cal V}$--modifications for $\mu$ and $\mathbf{1}_{\rho}.$ The three pentagons on the right
         are copies of the pentagon axiom for the composition ${\cal M}.$

    What we have developed here are the partial functors of the composition morphism implicit in enriching
    over $\mathbf{Cat}$. The said composition morphism is a functor of two variables. That the partial functors 
    can be combined to make the functor of two variables is implied by the 
    commutativity of a diagram that describes
    the two ways of combining them (see \cite{MacLane}).
    One thing that needs to be checked is that composing horizontally adjacent 
    2--cells is well--defined. We also need to check that the partial functors are indeed functorial. 
    This is shown by checking that the whiskering distributes over the vertical composition, and 
    checking that whiskering is the same as horizontally composing with an 
    identity 2--cell. (The latter is actually showing more than  
    that whiskering onto an identity 2--cell is the same as 
    horizontally composing two identity 2--cells, which in turn is more than what we really need: i.e. whiskering onto an identity
    2--cell gives an identity 2--cell for the composed 1--cells. It is often however, just as convenient to prove.) I
    start with the first axiom of functoriality.
    
    First we need to check that the whiskering distributes, 
    i.e. that $(\rho *\nu) \circ (\rho *\mu) = \rho *(\nu \circ \mu)$ and that $(\nu*\xi)\circ(\mu*\xi)=(\nu\circ\mu)*\xi$ 
    as in the following picture.
    $$
      
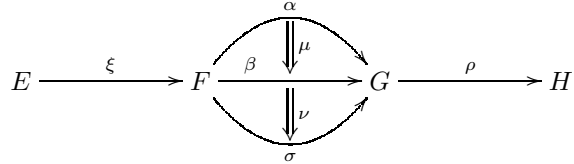
\xymatrix@R-=16pt{
      &&&\ar@{=>}[d]^{\mu}
      &&
      \\
      E
      \ar[rr]^{\xi}
      &&F
      \ar@/^2pc/[rr]^{\alpha}
      \ar[rr]^<<<<<{\beta}
      \ar@/_2pc/[rr]_{\sigma}
      &\ar@{=>}[d]^{\nu}
      &G
      \ar[rr]^{\rho}
      &
      &H
      \\
      &&&&&&&&
      }
    $$
    This requires the exteriors of the following two diagrams to commute
    (Let $\xi_U(0)=\overline{q}$, $\rho_U(0) = \overline{\overline{q}}$, $\alpha_U(0)=q,$ $\beta_U(0) = \hat{q}$ and 
                        $\sigma_U(0) = \check{q},$)

      \noindent
          		          \begin{center}
	          \resizebox{6in}{!}{ %!!!!
                 %%%%\resizebox{6.5in}{!}{
  \begin{sideways} 
  \begin{footnotesize}
            $$
            \xymatrix@R=25pt@C=-65pt{ %!!!!
            %%%% \xymatrix@R=35pt@C=-65pt{
            &&&&{\bcal W}(FU,HU)(\overline{\overline{q}}\hat{q},\overline{\overline{q}}\check{q})\otimes_1 {\bcal W}(FU,HU)(\overline{\overline{q}}{q},\overline{\overline{q}}\hat{q})
            \ar@<+.5ex>@/^4pc/[ddddddrr]^{M} 
            \\
            &&&{}\save[]*\txt{\text{ }\\$({\bcal W}(GU,HU)(\overline{\overline{q}},\overline{\overline{q}}) \otimes_2 {\bcal W}(FU,GU)(\hat{q},\check{q}))\otimes_1 ({\bcal W}(GU,HU)(\overline{\overline{q}},\overline{\overline{q}}) \otimes_2 {\bcal W}(FU,GU)({q},\hat{q}))$\\\text{ }}
            \ar[ur]^{{\cal M}\otimes_1 {\cal M}}
            %\ar[ddr]^{\eta}
            \restore
            \\
            &&&&&
            \\
            &&&&{}\save[]*\txt{\text{ }\\$({\bcal W}(GU,HU)(\overline{\overline{q}},\overline{\overline{q}}) \otimes_1 {\bcal W}(GU,HU)(\overline{\overline{q}},\overline{\overline{q}}))\otimes_2 ({\bcal W}(FU,GU)(\hat{q},\check{q}) \otimes_1 {\bcal W}(FU,GU)({q},\hat{q}))$\\\text{ }}
            \ar[dddr]^{M\otimes_2 M}
            \ar@{<-}[uul]^{\eta}
            \restore
            \\\\\\
            &\text{ }\text{ }\text{ }\text{ }\text{ }\text{ }\text{ }\text{ }\text{ }\text{ }\text{ }\text{ }\text{ }\text{ }\text{ }\text{ }\text{ }
            &&I=I\otimes_1 I = (I\otimes_2 I)\otimes_1 (I\otimes_2 I)\text{---------}{}^{\eta}\to (I\otimes_1 I)\otimes_2 (I\otimes_1 I)=I\otimes_2 (I\otimes_1 I)
            \ar[rr]^-{j_{\overline{\overline{q}}}\otimes_2 M\circ(\nu_U\otimes_1 \mu_U)}
            \ar@<-17ex>[uuu]|{(j_{\overline{\overline{q}}} \otimes_1 j_{\overline{\overline{q}}})\otimes_2 (\nu_U \otimes_1 \mu_U)}
            \ar@<12ex>[uuuuu]|{(j_{\overline{\overline{q}}} \otimes_2 \nu_U)\otimes_1 (j_{\overline{\overline{q}}} \otimes_2 \mu_U)}
            &&{\bcal W}(GU,HU)(\overline{\overline{q}},\overline{\overline{q}}) \otimes_2 {\bcal W}(FU,GU)(q,\check{q})\text{--------}{}^{\cal M}\to {\bcal W}(FU,HU)(\overline{\overline{q}}q,\overline{\overline{q}}\check{q})
            &\text{ }\text{ }\text{ }\text{ }\text{ }\text{ }\text{ }\text{ }\text{ }\text{ }\text{ }\text{ }\text{ }\text{ }\text{ }\text{ }\text{ }
            \\\\
            &\text{ }\text{ }\text{ }\text{ }\text{ }\text{ }\text{ }\text{ }\text{ }\text{ }\text{ }\text{ }\text{ }\text{ }\text{ }\text{ }\text{ }
              &&I=I\otimes_1 I = (I\otimes_2 I)\otimes_1 (I\otimes_2 I)\text{---------}{}^{\eta}\to (I\otimes_1 I)\otimes_2 (I\otimes_1 I)=I\otimes_2 (I\otimes_1 I)
              \ar[rr]^-{M\circ(\nu_U\otimes_1 \mu_U)\otimes_2 j_{\overline{q}}}
              \ar@<17ex>[ddd]|{( \nu_U\otimes_1 \mu_U)\otimes_2 (j_{\overline{q}} \otimes_1 j_{\overline{q}})}
              \ar@<-12ex>[ddddd]|{( \nu_U\otimes_2 j_{\overline{q}})\otimes_1 (\mu_U \otimes_2 j_{\overline{q}})}
              &&{\bcal W}(FU,GU)({q},\check{q}) \otimes_2 {\bcal W}(EU,FU)(\overline{q},\overline{q}\text{--------}{}^{\cal M}\to {\bcal W}(EU,GU)({q}\overline{q},\check{q}\overline{q}) 
              &\text{ }\text{ }\text{ }\text{ }\text{ }\text{ }\text{ }\text{ }\text{ }\text{ }\text{ }\text{ }\text{ }\text{ }\text{ }\text{ }\text{ }
              \\\\\\
              &&&&{}\save[]*\txt{\text{ }\\$({\bcal W}(FU,GU)(\hat{q},\check{q}) \otimes_1 {\bcal W}(FU,GU)({q},\hat{q}))\otimes_2 ({\bcal W}(EU,FU)(\overline{q},\overline{q}) \otimes_1 {\bcal W}(EU,FU)(\overline{q},\overline{q}))$\\\text{ }}
              \ar[uuur]^{M\otimes_2 M} \restore
              \\
              &&&&&
              \\
              &&&{}\save[]*\txt{\text{ }\\$({\bcal W}(FU,GU)(\hat{q},\check{q}) \otimes_2 {\bcal W}(EU,FU)(\overline{q},\overline{q}))\otimes_1 ({\bcal W}(FU,GU)({q},\hat{q}) \otimes_2 {\bcal W}(EU,FU)(\overline{q},\overline{q}))$\\\text{ }}
              \ar[uur]^{\eta}
              \ar[dr]^{{\cal M}\otimes_1 {\cal M}} \restore
              \\
              &&&&{\bcal W}(EU,GU)(\hat{q}\overline{q},\check{q}\overline{q}) \otimes_1 {\bcal W}(EU,GU)({q}\overline{q},\hat{q}\overline{q}) 
              \ar@<-.5ex>@/_4pc/[uuuuuurr]^{M}
              }
          $$
            \end{footnotesize}
             \end{sideways}
    }
  \end{center}
   
   %\clearpage
    %\newpage
    These commute since the interior regions all commute. The leftmost quadrilaterals commute by naturality of $\eta.$ The 
      central triangular regions commute by the unit axioms of ${\cal V}$--categories. The pentagonal regions
    commute by the ${\cal V}$--functoriality of ${\cal M}.$
    
    This commutativity has verified that the partial functors of the horizontal composition functor (whiskers) in fact do respect 
    the composition in their domain.

    We still need the two ways of composing the below cells using whiskers to be equivalent:
    
    $$
    \xymatrix@R-=3pt{
    &\ar@{=>}[dd]^{\mu}&&\ar@{=>}[dd]^{\nu}\\
    F
    \ar@/^1pc/[rr]^{\psi}
    \ar@/_1pc/[rr]_{\beta}
    &&G
    \ar@/^1pc/[rr]^{\gamma}
    \ar@/_1pc/[rr]_{\sigma}
    &&H
    \\
    &&&\\
    }
    $$
    
    That is, we need: 
    $$\nu * \mu = (\sigma *\mu) \circ (\nu *\psi) = (\nu *\beta) \circ (\gamma *\mu).$$
    In terms of the above definitions, the exterior of the following diagram must commute
    (Let $\rho_U(0) = \overline{\overline{q}}.$)

   \noindent
       		          \begin{center}
	          \resizebox{5.5in}{!}{ %!!!!
	          %%%% \resizebox{6.35in}{!}{
	          \begin{sideways}
      \begin{small}
      $$
      \xymatrix@C=-65pt{ %!!!!
      %%%% \xymatrix@R=40pt@C=-65pt{
      &&&&{\bcal W}(FU,HU)(\overline{\overline{q}}q,\overline{\overline{q}}\hat{q})\otimes_1 {\bcal W}(FU,HU)(\check{q}{q},\overline{\overline{q}}{q})
      \ar@<.5ex>@/^3pc/[ddddddrr]^{M}
      \\
      &&&{}\save[]*\txt{\text{ }\\$({\bcal W}(GU,HU)(\overline{\overline{q}},\overline{\overline{q}}) \otimes_2 {\bcal W}(FU,GU)(q,\hat{q}))\otimes_1 ({\bcal W}(GU,HU)(\check{q},\overline{\overline{q}}) \otimes_2 {\bcal W}(FU,GU)({q},{q}))$\\\text{ }}
      \ar[ur]^{{\cal M}\otimes_1 {\cal M}}
      %\ar[ddr]^{\eta}
      \restore
      \\
      &&&&&
      \\
      &&&&{}\save[]*\txt{\text{ }\\$({\bcal W}(GU,HU)(\overline{\overline{q}},\overline{\overline{q}}) \otimes_1 {\bcal W}(GU,HU)(\check{q},\overline{\overline{q}}))\otimes_2 ({\bcal W}(FU,GU)(q,\hat{q}) \otimes_1 {\bcal W}(FU,GU)({q},{q}))$\\\text{ }}
      \ar[dddr]^{M\otimes_2 M}
      \ar@{<-}[uul]^{\eta}
      \restore
      \\\\\\
      &\text{ }\text{ }\text{ }\text{ }\text{ }\text{ }\text{ }\text{ }\text{ }\text{ }\text{ }\text{ }\text{ }\text{ }\text{ }\text{ }\text{ }
      &&I=I\otimes_1 I = (I\otimes_2 I)\otimes_1 (I\otimes_2 I)\text{---------}{}^{\eta}\to (I\otimes_1 I)\otimes_2 (I\otimes_1 I)=I\otimes_2 I
      \ar@<-17ex>[uuu]|{(j_{\overline{\overline{q}}} \otimes_1 \nu_U)\otimes_2 (\mu_U \otimes_1 j_{q})}
      \ar@<17ex>[ddd]|{( \nu_U\otimes_1 j_{\check{q}})\otimes_2 (j_{\hat{q}} \otimes_1 \mu_U)}
      \ar[rr]^-{\nu_U\otimes_2 \mu_U}
      \ar@<12ex>[uuuuu]|{(j_{\overline{\overline{q}}} \otimes_2 \mu_U)\otimes_1 ( \nu_U\otimes_2 j_{q})}
      \ar@<-12ex>[ddddd]|{( \nu_U\otimes_2 j_{\hat{q}})\otimes_1 (j_{\check{q}} \otimes_2 \mu_U)}
      &&{\bcal W}(GU,HU)(\check{q},\overline{\overline{q}}) \otimes_2 {\bcal W}(FU,GU)(q,\hat{q})\text{--------}{}^{\cal M}\to {\bcal W}(FU,HU)(\check{q}q,\overline{\overline{q}}\hat{q})
      &\text{ }\text{ }\text{ }\text{ }\text{ }\text{ }\text{ }\text{ }\text{ }\text{ }\text{ }\text{ }\text{ }\text{ }\text{ }\text{ }\text{ }
      \\\\\\
      &&&&{}\save[]*\txt{\text{ }\\$({\bcal W}(GU,HU)(\check{q},\overline{\overline{q}}) \otimes_1 {\bcal W}(GU,HU)(\check{q},\check{q}))\otimes_2 ({\bcal W}(FU,GU)(\hat{q},\hat{q}) \otimes_1 {\bcal W}(FU,GU)({q},\hat{q}))$\\\text{ }}
      \ar[uuur]^{M\otimes_2 M} \restore
      \\
      &&&&&
      \\
      &&&{}\save[]*\txt{\text{ }\\$({\bcal W}(GU,HU)(\check{q},\overline{\overline{q}}) \otimes_2 {\bcal W}(FU,GU)(\hat{q},\hat{q}))\otimes_1 ({\bcal W}(GU,HU)(\check{q},\check{q}) \otimes_2 {\bcal W}(FU,GU)({q},\hat{q}))$\\\text{ }}
      \ar[uur]^{\eta}
      \ar[dr]^{{\cal M}\otimes_1 {\cal M}} \restore
      \\
      &&&&{\bcal W}(FU,HU)(\check{q}\hat{q},\overline{\overline{q}}\hat{q}) \otimes_1 {\bcal W}(FU,HU)(\check{q}{q},\check{q}\hat{q}) 
      \ar@<-.5ex>@/_3pc/[uuuuuurr]^{M}
      }
      $$
      \end{small}
        \end{sideways}
    }
  \end{center}
  
  %\clearpage
    %\newpage
    This commutes since the interior regions all commute. The leftmost quadrilaterals commute by naturality of $\eta.$ The 
    central triangular regions commute by the unit axioms of ${\cal V}$--categories. The upper and lower pentagonal regions
    commute by the ${\cal V}$--functoriality of ${\cal M}.$ This composition gives a valid ${\cal V}$--modification
    since the whiskered pieces are valid and since the composition along a ${\cal V}$--2--natural transformation gives a 
    valid ${\cal V}$--modification. The central leg of the above diagram gives a more direct description of the composition
    of ${\cal V}$--2--modifications along a ${\cal V}$--2--functor.  From this description it is automatic that whiskering a
    ${\cal V}$--2--natural transformation to a ${\cal V}$--modification along a ${\cal V}$--2--functor is the same as 
    composing along that ${\cal V}$--2--functor with an identity ${\cal V}$--modification corresponding to the whisker. We
    also see from this description that the associativity of this composition follows immediately from the associativity axiom of
    a strict ${\cal V}$--2--category.
    
     Now we can show the functoriality of the entire composition functor. (The general
  proof regarding partial functors of a functor of two variables is in \cite{MacLane}.)
    This states that, in the following picture, 
    $({\nu_2}*{\nu_1})\circ ({\mu_2}*{\mu_1}) = ({\nu_2}\circ {\mu_2})*({\nu_1}\circ {\mu_1}).$
    
    $$
    
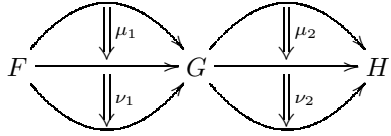
\xymatrix@R-=16pt{
    &\ar@{=>}[d]^{\mu_1}
    &&\ar@{=>}[d]^{\mu_2}
    \\
    F
    \ar@/^2pc/[rr]
    \ar[rr]
    \ar@/_2pc/[rr]
    &\ar@{=>}[d]^{\nu_1}
    &G
    \ar@/^2pc/[rr]
    \ar[rr]
    \ar@/_2pc/[rr]
    &\ar@{=>}[d]^{\nu_2}
    &H
    \\
    &&&&&
    }
    $$
    
    This is shown by using what we have established. The following series of pictures serve to illustrate the proof.
    
          \noindent
              		          \begin{center}
    	          \resizebox{6in}{!}
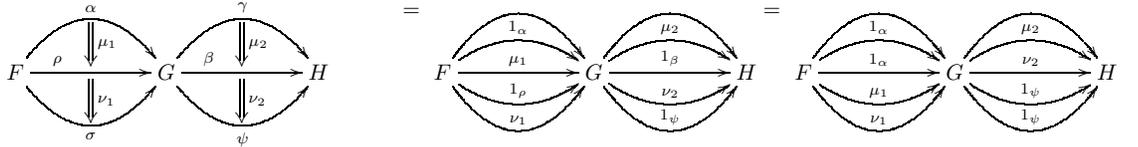
{ %!!!!
                 %%%%\resizebox{6.5in}{!}{
    $$
      \xymatrix@R-=16pt{
      &\ar@{=>}[d]^{\mu_1}
      &&\ar@{=>}[d]^{\mu_2}
      \\
      F
      \ar@/^2pc/[rr]^{\alpha}
      \ar[rr]^<<<<<{\rho}
      \ar@/_2pc/[rr]_{\sigma}
      &\ar@{=>}[d]^{\nu_1}
      &G
      \ar@/^2pc/[rr]^{\gamma}
      \ar[rr]^<<<<<{\beta}
      \ar@/_2pc/[rr]_{\psi}
      &\ar@{=>}[d]^{\nu_2}
      &H
      \\
      &&&&&
      }
    =
    \xymatrix@R-=16pt{
        &
        &&
        \\
        F
        \ar@/^2.2pc/[rr]
        \ar@/^1.2pc/[rr]^{1_{\alpha}}
        \ar[rr]^{\mu_1}
        \ar@/_1.2pc/[rr]^{1_{\rho}}
        \ar@/_2.2pc/[rr]^{\nu_1}
        &
        &G
        \ar@/^2.2pc/[rr]
        \ar@/^1.2pc/[rr]^{\mu_2}
        \ar[rr]^{1_{\beta}}
        \ar@/_1.2pc/[rr]^{\nu_2}
        \ar@/_2.2pc/[rr]^{1_{\psi}}
        &
        &H
        \\
        &&
      }
  =
  \xymatrix@R-=16pt{
        &
        &&
        \\
        F
        \ar@/^2.2pc/[rr]
        \ar@/^1.2pc/[rr]^{1_{\alpha}}
        \ar[rr]^{1_{\alpha}}
        \ar@/_1.2pc/[rr]^{\mu_1}
        \ar@/_2.2pc/[rr]^{\nu_1}
        &
        &G
        \ar@/^2.2pc/[rr]
        \ar@/^1.2pc/[rr]^{\mu_2}
        \ar[rr]^{\nu_2}
        \ar@/_1.2pc/[rr]^{1_{\psi}}
        \ar@/_2.2pc/[rr]^{1_{\psi}}
        &
        &H
        \\
        &&
      }
  $$
  }
  \end{center}
  
  Or symbolically:
  $$({\nu_2}*{\nu_1})\circ ({\mu_2}*{\mu_1})$$ 
  $$= (({\psi}*\nu_1)\circ(\nu_2*{\rho}))\circ (({\beta}*\mu_1)\circ(\mu_2*{\alpha}))$$
  $$= ({\psi}*\nu_1)\circ((\nu_2*{\rho})\circ ({\beta}*\mu_1))\circ(\mu_2*{\alpha})$$
  $$= ({\psi}*\nu_1)\circ(({\psi}*\mu_1)\circ (\nu_2*{\alpha}))\circ(\mu_2*{\alpha})$$
  $$= (({\psi}*\nu_1)\circ({\psi}*\mu_1))\circ ((\nu_2*{\alpha})\circ(\mu_2*{\alpha}))$$
  $$= ({\psi}*(\nu_1\circ\mu_1))\circ ((\nu_2\circ\mu_2)*{\alpha})$$
  $$= ({\nu_2}\circ {\mu_2})*({\nu_1}\circ {\mu_1}) $$
  
  The exchange identity is precisely the functoriality (respect of vertical composition) of the functor of 
  two variables that describes the horizontal composition.
  We also have the respect of units by the horizontal composition simply by using the exchange identity above with
  units in the lower two 2--cells.
  Of course the 
  central leg version of the composition can be directly verified to be functorial. The roundabout route is nice since it 
  covers lots of pasting diagrams, and is a good model for future such verifications.
  
  The unit for composing ${\cal V}$--modifications along a ${\cal V}$--functor is the identity ${\cal V}$--modification
    $\mathbf{1}_{1_T}$ where $1_T:T\to T$ is the identity ${\cal V}$--2--natural transformation for a ${\cal V}$--2--functor $T.$
    Since ${(1_T)}_U = {\cal J}_{TU}$ then ${(1_T)}_U(0)={\cal J}_{TU}(0)= 1_{TU}.$ 
    Thus ${({(\mathbf{1}_{1_T})}_U)}_0 = j_{1_{TU}}={\cal J}_{{TU}_{00}}.$ That $\mathbf{1}_{1_T}$ is a 2--sided unit for the
    composition is seen by the unit axioms of  strict ${\cal V}$--2--category.

\end{proof}

  %{part2}

  %\clearpage
  %\newpage
   \begin{proof}\text{ }(Part 2.)
   
   In part 2 of the proof we describe how for each triple of ${\cal V}$--2--categories
   we have a 2--functor of two variables that serves to compose morphisms along a 
     common ${\cal V}$--2--category as in the following picture.
     
     $$
       
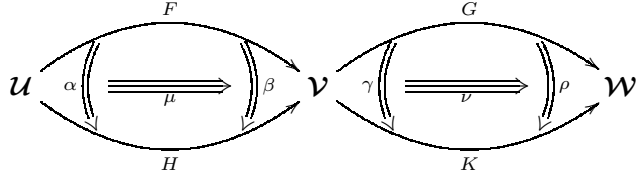
\xymatrix@R=10pt@C=20pt{
       &\ar@/_.5pc/@{=>}[dd]_{\alpha}
       &&\ar@/^.5pc/@{=>}[dd]^{\beta}
       &&\ar@/_.5pc/@{=>}[dd]_{\gamma}
       &&\ar@/^.5pc/@{=>}[dd]^{\rho}
       \\
       {\bcal U}
       \ar@/^2pc/[rrrr]^F
       \ar@/_2pc/[rrrr]_H
       &
       \ar@3{->}[rr]_{\mu}
       &&&{\bcal V}
       \ar@/^2pc/[rrrr]^G
       \ar@/_2pc/[rrrr]_K
       &
       \ar@3{->}[rr]_{\nu}
       &&&{\bcal W}\\
       &&&&&&&&&&&
       }
    $$
     
     At each stage of description we also need to check that the composition
     along a common ${\cal V}$--2--category is associative and respects all units, as well as making sure that
     for morphisms the composition is functorial. This latter property always exhibits itself as an exchange identity. 
  
  Composition of ${\cal V}$--2--functors is just composition of the object functions and 
  composition of the hom--category ${\cal V}$--functors, with appropriate subscripts. Thus
  $(ST)_{UU'}(f) = S_{TUTU'}(T_{UU'}(f)).$ Then it is straightforward to verify that the axioms are obeyed,
  as in $$(ST)_{U'U''}(f)(ST)_{UU'}(g)$$ 
  $$= S_{TU'TU''}(T_{U'U''}(f))S_{TUTU'}(T_{UU'}(g))$$ 
  $$= S_{TUTU''}(T_{U'U''}(f)T_{UU'}(g))$$
  $$= S_{TUTU''}(T_{UU''}(fg))$$
  $$= (ST)_{UU''}(fg).$$
  That this composition is associative follows from the associativity of composition of the underlying functions and
  ${\cal V}$--functors. The 2--sided identity for this composition $1_{\bcal U}$ is made of the 
  identity function (on objects) and
  identity ${\cal V}$--functors (for hom--categories.) 
  
  Next we define the composition of ${\cal V}$--2--natural transformations along a ${\cal V}$--2--category. This is accomplished by first
  describing the trivial cases--whiskering a ${\cal V}$--2--functor to a ${\cal V}$--2--natural transformation along a ${\cal V}$--2--category.
  (By proceeding in terms of whiskers I get the opportunity to both discuss all the possible pasting in the 3--category
  and to exhibit the sub--2--categories implicit in its structure.)
  
  The first 
    picture shows a 1-cell (${\cal V}$--2--functor) following a 2-cell (${\cal V}$--2--natural transformation).
    These are composed to form a new 2-cell as follows
    
    $$
    \xymatrix@R-=3pt{
    &\ar@{=>}[dd]^{\alpha}\\
    {\bcal U}
    \ar@/^1pc/[rr]^F
    \ar@/_1pc/[rr]_H
    &&{\bcal V}
    \ar[rr]^G
    &&{\bcal W}
    \\
    &\\
    } \text{   is composed to become } \xymatrix@R-=3pt{
    &\ar@{=>}[dd]^{G\alpha}\\
    {\bcal U}
    \ar@/^1pc/[rr]^{GF}
    \ar@/_1pc/[rr]_{GH}
    &&{\bcal W}
    \\
    &\\
    }
    $$
    
    where $G\alpha$
    has components given by
    $$\xymatrix{(G\alpha)_{U}} = \xymatrix{
                                    {\cal I}  
                                    \ar[d]_{\alpha_{U}}
                                    \\{\bcal V}(FU,HU)
                                    \ar[d]_{G_{FU,HU}}
                                    \\{\bcal W}(GFU,GHU)}$$
                                    
    Notice that by definition and axioms of enrichment, 
    letting $\alpha_U(0)=q$, we have $(G\alpha)_U(0) = Gq$ and $(G\alpha)_{U_{00}} = j_{Gq}.$
     
     This whiskering gives a valid ${\cal V}$--2--natural transformation by the following commuting diagram.
    %\clearpage 
    
    \noindent
        		          \begin{center}
	          \resizebox{6.5in}{!}
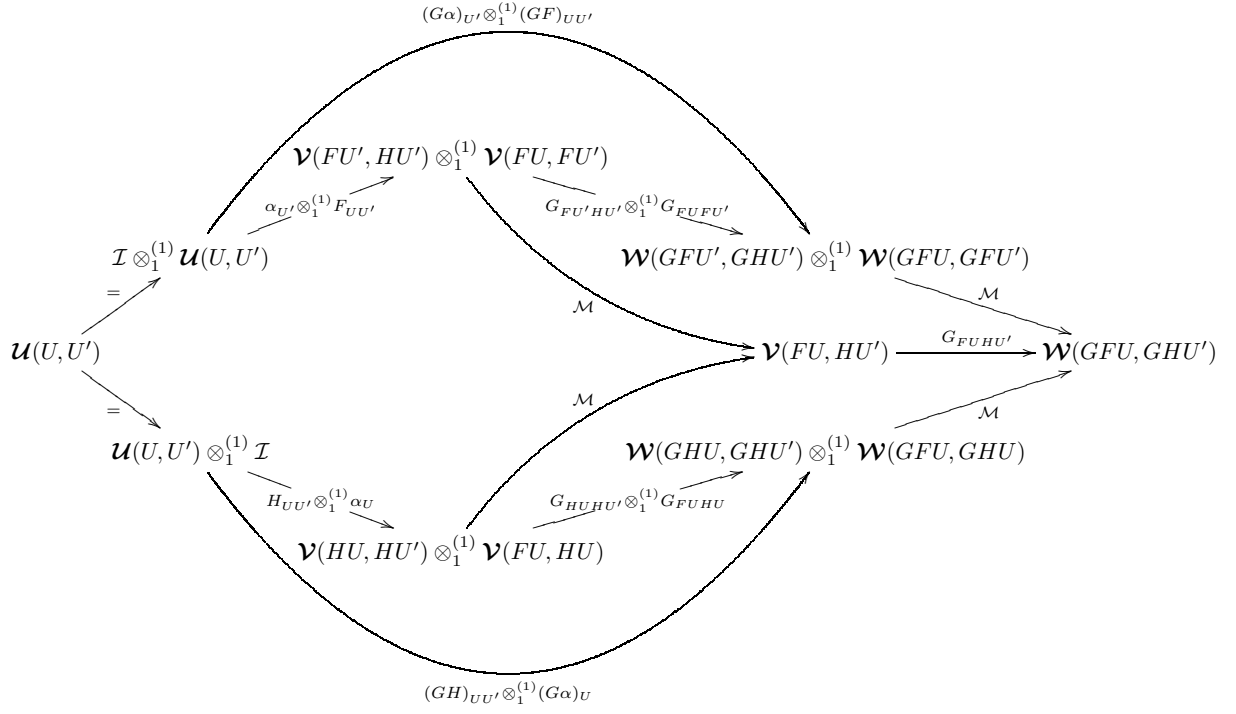
{ %!!!!
	          %%%%\resizebox{7.35in}{!}{
    $$
             \xymatrix@C=-1pt{
             &&&{\bcal V}(FU',HU') \otimes^{(1)}_1 {\bcal V}(FU,FU')
             \ar@/_2pc/[rdd]_-{{\cal M}}
             \ar[rd]|{G_{FU'HU'}\otimes^{(1)}_1 G_{FUFU'}}
             \\
             &{\cal I} \otimes^{(1)}_1 {\bcal U}(U,U')
             \ar[urr]|-{ \alpha_{U'} \otimes^{(1)}_1 F_{UU'}}
             \ar@/^8pc/[rrr]^{(G\alpha)_{U'}\otimes^{(1)}_1 (GF)_{UU'}}
             &&&{\bcal W}(GFU',GHU') \otimes^{(1)}_1 {\bcal W}(GFU,GFU')
             \ar[rd]^-{{\cal M}}
           \\
             {\bcal U}(U,U')
             \ar[ru]^{=}
           \ar[rd]_{=}
             &&&&{\bcal V}(FU,HU')
             \ar[r]^{G_{FUHU'}}
             &{\bcal W}(GFU,GHU')
           \\
             &{\bcal U}(U,U') \otimes^{(1)}_1 {\cal I}
             \ar[drr]|-{H_{UU'} \otimes^{(1)}_1 \alpha_{U}}
             \ar@/_8pc/[rrr]_{(GH)_{UU'}\otimes^{(1)}_1 (G\alpha)_{U}}
             &&
             &{\bcal W}(GHU,GHU') \otimes^{(1)}_1 {\bcal W}(GFU,GHU)
             \ar[ru]_-{{\cal M}}
             \\
             &&&{\bcal V}(HU,HU') \otimes^{(1)}_1 {\bcal V}(FU,HU)
             \ar@/^2pc/[ruu]^-{{\cal M}}
             \ar[ru]|{G_{HUHU'}\otimes^{(1)}_1 G_{FUHU}}
             }
         $$
         }
  \end{center}
    
    The central region expresses the ${\cal V}$--2--naturality of $\alpha.$ The two rightmost regions commute by the definition
    of ${\cal V}$--2--functor.
    
    The second picture shows a 2-cell following a 1-cell. These are composed as follows
    $$
    
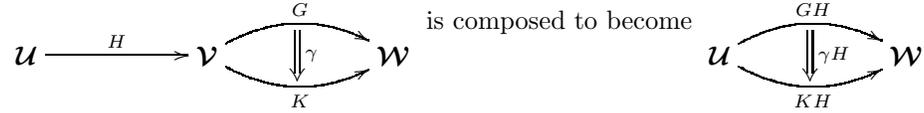
\xymatrix@R-=3pt{
    &&&\ar@{=>}[dd]^{\gamma}\\
    {\bcal U}
    \ar[rr]^H
    &&{\bcal V}
    \ar@/^1pc/[rr]^G
    \ar@/_1pc/[rr]_K
    &&{\bcal W}
    \\
    &&&\\
    } \text{   is composed to become } \xymatrix@R-=3pt{
    &\ar@{=>}[dd]^{\gamma H}\\
    {\bcal U}
    \ar@/^1pc/[rr]^{GH}
    \ar@/_1pc/[rr]_{KH}
    &&{\bcal W}
    \\
    &\\
    }
    $$
    
    where $\gamma H$ has components given by 
    $(\gamma H)_{U} = \gamma_{HU}$. 
  
  This whiskering gives a valid ${\cal V}$--2--natural transformation by the following commuting diagram.
  %\clearpage
  
  \noindent
      		          \begin{center}
	          \resizebox{6.5in}{!}{ %!!!!
	          %%%%\resizebox{7.35in}{!}{
  $$
           \xymatrix@C=-1pt{
           &&{\cal I} \otimes^{(1)}_1 {\bcal V}(HU,HU')
           \ar[rrd]|-{ \gamma_{HU'} \otimes^{(1)}_1 G_{HUHU'}}
           \\
           &{\cal I} \otimes^{(1)}_1 {\bcal U}(U,U') \text{ }\text{ }\text{ }\text{ }\text{ }\text{ }\text{ }\text{ }\text{ }\text{ }\text{ }\text{ }\text{ }\text{ }\text{ }\text{ }
           \ar[ru]|{1\otimes^{(1)}_1 H_{UU'}}
           \ar@/^7.2pc/[rrr]^{(\gamma H)_{U'}\otimes^{(1)}_1 (GH)_{UU'}}
           &
           &&{\bcal W}(GHU',KHU') \otimes^{(1)}_1 {\bcal W}(GHU,GHU')
           \ar[rd]^-{{\cal M}}
         \\
         {\bcal U}(U,U')
         \ar[r]^{H_{UU'}}
         \ar[ur]^{=}
         \ar[dr]_{=}
         &{\bcal V}(HU,HU')
         \ar@/_1pc/[ruu]_{=}
         \ar@/^1pc/[rdd]^{=}
           &&&&{\bcal W}(GHU,KHU')
         \\
           &{\bcal U}(U,U') \otimes^{(1)}_1 {\cal I} \text{ }\text{ }\text{ }\text{ }\text{ }\text{ }\text{ }\text{ }\text{ }\text{ }\text{ }\text{ }\text{ }\text{ }\text{ }\text{ }
           \ar[rd]|{H_{UU'}\otimes^{(1)}_1 1}
           \ar@/_7.2pc/[rrr]_{(KH)_{UU'}\otimes^{(1)}_1 (\gamma H)_{U}}
           &
           &&{\bcal W}(KHU,KHU') \otimes^{(1)}_1 {\bcal W}(GHU,KHU)
           \ar[ru]^-{{\cal M}}
         \\
         &&{\bcal V}(HU,HU') \otimes^{(1)}_1 {\cal I}
           \ar[rru]|-{K_{HUHU'} \otimes^{(1)}_1 \gamma_{HU}}  
           }
         $$
         }
  \end{center}
  
  The central region expresses the ${\cal V}$--2--naturality of $\gamma.$
    
    To show the exchange identity here we proceed by checking the usual agreement and functoriality of partials.
    First I will check that the partial functors described by whiskering are
    indeed functorial. These proofs continue to parallel the lower dimensional case. First we check that the right 
    whiskering distributes, 
      i.e. that $(Z \alpha) * (Z \beta) = Z (\alpha * \beta)$ as in the following picture. (Recall that ``$*$'' denotes the
      composition along ${\cal V}$--2--functors as in the first part of the proof.)
      $$
        \xymatrix@R-=16pt{
        &\ar@{=>}[d]^{\beta}
        &&
        \\
        {\bcal U}
        \ar@/^2pc/[rr]^{T}
        \ar[rr]^<<<<<{S}
        \ar@/_2pc/[rr]_{R}
        &\ar@{=>}[d]^{\alpha}
        &{\bcal V}
        \ar[rr]^{Z}
        &
        &{\bcal W}
        \\
        &&&&&
        }
      $$
    The two sides of the proposed equality form the legs of the following diagram, which commutes due to 
    the definition of the 
    ${\cal V}$--2--functoriality of Z:
    $$
             \xymatrix@C=-5pt{
             &{\cal I} = {\cal I} \otimes^{(1)}_1 {\cal I}
             \ar[d]^{\alpha_U \otimes^{(1)}_1 \beta_U}
             \\
             &{\bcal V}(SU,RU) \otimes^{(1)}_1 {\bcal V}(TU,SU)
             \ar[dl]^{{\cal M}}
             \ar[dr]^{Z\otimes^{(1)}_1 Z}
             \\
             {\bcal V}(TU,RU)
             \ar[dr]^{Z}
             &&{\bcal W}(ZSU,ZRU) \otimes^{(1)}_1 {\bcal V}(ZTU,ZSU)
             \ar[dl]^{{\cal M}}
             \\
             &{\bcal W}(ZTU,ZRU) 
             }
           $$
  %\clearpage
  For the same requirement on the other partial functor we have the picture
  $$
    \xymatrix@R-=16pt{
    &
    &&\ar@{=>}[d]^{\gamma}
    \\
    {\bcal U}
    \ar[rr]^{T}
    &
    &{\bcal V}
    \ar@/^2pc/[rr]
    \ar[rr]
    \ar@/_2pc/[rr]
    &\ar@{=>}[d]^{\delta}
    &{\bcal W}
    \\
    &&&&&
    }
    $$
  From the definitions is it immediate that $((\delta * \gamma)T)_U = (\delta * \gamma)_{TU} = (\delta T * \gamma T)_U.$
  
  Now we can compose ${\cal V}$--2--natural transformations along a ${\cal V}$--2--category, as in the following
  picture.
  $$
    \xymatrix@R-=3pt{
    &\ar@{=>}[dd]^{\alpha}&&\ar@{=>}[dd]^{\gamma}\\
    {\bcal U}
    \ar@/^1pc/[rr]^F
    \ar@/_1pc/[rr]_H
    &&{\bcal V}
    \ar@/^1pc/[rr]^G
    \ar@/_1pc/[rr]_K
    &&{\bcal W}
    \\
    &&&\\
    }
    $$
  
  As usual there are two ways to 
  do so that need to be reconciled. They both consist of defining the composition along the ${\cal V}$--2--category 
  in terms of a composition along a common ${\cal V}$--2--functor as in part 1 of the proof. Thus since the whiskered
  pieces are valid ${\cal V}$--2--natural transformations, by a previous diagram their composition will be as well.
  The first way of composing is given by:
     $$\xymatrix{\gamma\alpha} =
           \xymatrix{
           {\cal I} = {\cal I} \otimes^{(1)}_1 {\cal I}
           \ar[d]^{(\gamma H)_U \otimes^{(1)}_1 (G\alpha)_U}
           \\
           {\bcal W}(GHU,KHU) \otimes^{(1)}_1 {\bcal W}(GFU,GHU)
           \ar[d]^{{\cal M}}
           \\
           {\bcal W}(GFU,KHU)
           }
           $$
  
  The second is given by
  $$\xymatrix{\gamma\alpha} =
           \xymatrix{
           {\cal I} = {\cal I} \otimes^{(1)}_1 {\cal I}
           \ar[d]^{(K\alpha)_U \otimes^{(1)}_1 (\gamma F)_U}
           \\
           {\bcal W}(KFU,KHU) \otimes^{(1)}_1 {\bcal W}(GFU,KFU)
           \ar[d]^{{\cal M}}
           \\
           {\bcal W}(GFU,KHU)
           }
           $$
  Letting $(\gamma F)_U(0) = \hat{q}$ and $(\gamma H)_U(0) = \hat{q'}$ and recalling that $(G\alpha)_U(0) = Gq$ we have
  $(\gamma\alpha)_U(0) = \hat{q'}Gq = Kq\hat{q} $ and by ${\cal V}$--functoriality of ${\cal M}$ that 
  $(\gamma\alpha)_{U_{00}} = j_{\hat{q'}Gq}.$

   That the two ways of composing are actually the same is based on the ${\cal V}$--2--naturality of $\gamma,$ the definition of which 
  makes up the central region of the following commuting diagram. The other regions commute trivially.
  %\clearpage
  
  \noindent
      		          \begin{center}
	          \resizebox{6.55in}{!}{ %!!!!
	          %%%%\resizebox{7.45in}{!}{
  $$
           \xymatrix@C=-2pt{
           &&&{\cal I} \otimes^{(1)}_1 {\bcal V}(FU,HU)
           \ar[rrd]|-{ \gamma_{HU} \otimes^{(1)}_1 G_{FUHU}}
           \\
           &\text{ }\text{ }\text{ }\text{ }\text{ }\text{ }\text{ }\text{ }\text{ }\text{ }\text{ }\text{ }\text{ }\text{ }
           &&&&{\bcal W}(GHU,KHU) \otimes^{(1)}_1 {\bcal W}(GFU,GHU)
           \ar[rd]^-{{\cal M}}
         \\
           {\cal I} \otimes^{(1)}_1 {\cal I}
  	          \ar[r]^-{=}
  	          \ar[uurrr]^{1\otimes^{(1)}_1 \alpha_U}
           \ar[ddrrr]_{\alpha_U\otimes^{(1)}_1 1}
           \ar@/^7.5pc/[rrrrru]^{(\gamma H)_U \otimes^{(1)}_1 (G\alpha)_U}
           \ar@/_7.5pc/[rrrrrd]_{(K\alpha)_U \otimes^{(1)}_1 (\gamma F)_U}
           &{\cal I}
           \ar[r]^-{\alpha_U}
           &{\bcal V}(FU,HU)
           \ar[ruu]^{=}
         \ar[rdd]_{=}
           &&&&{\bcal W}(GFU,KHU)
         \\
           &\text{ }\text{ }\text{ }\text{ }\text{ }\text{ }\text{ }\text{ }\text{ }\text{ }\text{ }\text{ }\text{ }\text{ }
           &&&&{\bcal W}(KFU,KHU) \otimes^{(1)}_1 {\bcal W}(GFU,KFU)
           \ar[ru]^-{{\cal M}}
           \\
           &&&{\bcal V}(FU,HU) \otimes^{(1)}_1 {\cal I}
           \ar[rru]|-{K_{FUHU} \otimes^{(1)}_1 \gamma_{FU}}
           }
         $$
         }
  \end{center}

  Now whiskering 
    a 1--cell Q on the right (or left) of a 2--cell $\alpha: T \to S$  should be the same as horizontally 
    composing ${\textbf 1}_{Q}$ on the respective side of $\alpha$.
    Pictorially for the right-hand whiskering:
    $$
    \xymatrix@R-=3pt{
    &\ar@{=>}[dd]^{\alpha}\\
    {\bcal U}
    \ar@/^1pc/[rr]^T
    \ar@/_1pc/[rr]_S
    &&{\bcal V}
    \ar[rr]^Q
    &&{\bcal W}
    \\
    &\\
    } = \xymatrix@R-=3pt{
    &\ar@{=>}[dd]^{\alpha}&&\ar@{=>}[dd]^{{\textbf 1}_{Q}}\\
    {\bcal U}
    \ar@/^1pc/[rr]^T
    \ar@/_1pc/[rr]_S
    &&{\bcal V}
    \ar@/^1pc/[rr]^Q
    \ar@/_1pc/[rr]_Q
    &&{\bcal W}
    \\
    &&&\\
    }
    $$
    To see this equality we need check only one way of composing 
    ${\textbf 1}_{Q}\alpha$ since we have shown it to be well defined -- i.e. we check that $Q\alpha =  
    {\textbf 1}_{Q}\alpha = Q\alpha * {\textbf 1}_{Q}T$. This is true 
    immediately from the relationship of ${\cal M}$ and ${\cal J}$. Now pictorially for the left-hand whiskering:
    
    $$
    \xymatrix@R-=3pt{
    &&&\ar@{=>}[dd]^{\alpha}\\
    {\cal D}
    \ar[rr]^P
    &&{\cal A}
    \ar@/^1pc/[rr]^T
    \ar@/_1pc/[rr]_S
    &&{\cal B}
    \\
    &&&\\
    } = \xymatrix@R-=3pt{
    &\ar@{=>}[dd]^{{\textbf 1}_{P}}&&\ar@{=>}[dd]^{\alpha}\\
    {\cal D}
    \ar@/^1pc/[rr]^P
    \ar@/_1pc/[rr]_P
    &&{\cal A}
    \ar@/^1pc/[rr]^T
    \ar@/_1pc/[rr]_S
    &&{\cal B}
    \\
    &&&\\
    }
    $$
    That $\alpha P = \alpha{\textbf 1}_{P} = S{\textbf 1}_{P} * \alpha P$ also is shown by using the relationship
    of ${\cal M}$ and ${\cal J}$ and by the ${\cal V}$--2--functoriality of $S.$
  
  Associativity of this composition follows from the associativity of composing ${\cal V}$--2--natural transformations
  along a ${\cal V}$--2--functor. It also requires the functoriality of the partial functors. In the following picture
  
  $$
         \xymatrix@R-=3pt{
         &\ar@{=>}[dd]^{\alpha}&&\ar@{=>}[dd]^{\gamma}&&\ar@{=>}[dd]^{\beta}\\
         {\bcal U}
         \ar@/^1pc/[rr]^F
         \ar@/_1pc/[rr]_H
         &&{\bcal V}
         \ar@/^1pc/[rr]^G
         \ar@/_1pc/[rr]_K
         &&{\bcal W}
         \ar@/^1pc/[rr]^P
         \ar@/_1pc/[rr]_Q
         &&{\bcal X}
         \\
         &&&&&&\\
         }      
         $$
  we have $$\beta(\gamma\alpha) = Q(K\alpha * \gamma F)* \beta GF$$
  $$= (QK\alpha * Q \gamma F) * \beta GF$$
  $$= QK\alpha * (Q \gamma F * \beta GF)$$
  $$= QK\alpha * (Q \gamma  * \beta G)F = (\beta\gamma)\alpha$$
  %\clearpage
  where the assumed associativities of whiskers are easily verified.
  
  The unit for composing ${\cal V}$--2--natural transformations along a ${\cal V}$--2--category is the 
  identity ${\cal V}$--2--natural transformation
    $\mathbf{1}_{1_{\bcal U}}$ where $1_{\bcal U}:{\bcal U}\to {\bcal U}$ is the identity ${\cal V}$--2--functor.
    Note that $(\mathbf{1}_{1_{\bcal U}})_U(0)= {\cal J}_{1_{\bcal U}U}(0) = {\cal J}_U(0) = 1_U$ and that
    ${(\mathbf{1}_{1_{\bcal U}})_U}_{00} = {\cal J}_{U_{00}} = j_{1_U}.$ Since the composition is based on that 
    of  ${\cal V}$--2--natural transformations along a ${\cal V}$--2--functor, to see that $\mathbf{1}_{1_{\bcal U}}$
    is a 2--sided unit all we need to check is that for any $\alpha:T\to S:{\bcal U}\to{\bcal W}$ we have that
    $1_{\bcal W}\alpha = \alpha = \alpha 1_{\bcal U}.$ This is clear from the definitions of whiskering above.
  
  Next we consider compositions involving ${\cal V}$--modifications along a ${\cal V}$--2--category.
  I start by defining whiskering of ${\cal V}$--2--functors and then use that definition to define whiskering of
  ${\cal V}$--2--natural transformations. 
  First the right whiskering of a ${\cal V}$--2--functor onto a 
  ${\cal V}$--modification as in the picture:
  $$
       \xymatrix@R=10pt@C=20pt{
       &\ar@/_.5pc/@{=>}[dd]_{\alpha}
       &&\ar@/^.5pc/@{=>}[dd]^{\beta}
       &&
       &
       &
       \\
       {\bcal U}
       \ar@/^2pc/[rrrr]^F
       \ar@/_2pc/[rrrr]_H
       &
       \ar@3{->}[rr]_{\mu}
       &&&{\bcal V}
       \ar[rrrr]_K
       &
       &&&{\bcal W}\\
       &&&&&&&&
       }
       \mapsto
       \xymatrix@R=10pt@C=20pt{
            &\ar@/_.5pc/@{=>}[dd]_{K\alpha}
            &&\ar@/^.5pc/@{=>}[dd]^{K\beta}
            &
            \\
            {\bcal U}
            \ar@/^2pc/[rrrr]^{KF}
            \ar@/_2pc/[rrrr]_{KH}
            &
            \ar@3{->}[rr]_{K\mu}
            &&&{\bcal W}
            \\
       &&&&
       }
    $$
  where:
  $$\xymatrix{(K\mu)_{U}} = \xymatrix{
                                     I  
                                    \ar[d]_{\mu_{U}}
                                    \\{\bcal V}(FU,HU)(q,\check{q})
                                    \ar[d]_{K_{FU,HU}}
                                    \\{\bcal W}(KFU,KHU)(Kq,K\check{q})}$$
  Where we let $\alpha_U(0)= q$ and $\beta_{U}(0)= \check{q}.$
  That this forms a valid ${\cal V}$--modification can be seen upon inspecting the following diagram. Its commutativity
  relies on the fact that $\mu$ is a ${\cal V}$--modification and on the ${\cal V}$--2--functoriality of $K.$
  
  \noindent
        		          \begin{center}
	          \resizebox{6.55in}{!}{ %!!!!
	          %%%%\resizebox{7.45in}{!}{
  $$
             \xymatrix@C=-73pt@R=35pt{
             &&&{\bcal W}(KFU',KHU')(Kq',K\check{q'}) \otimes_2 {\bcal W}(KFU,KFU')(KFf,KFg)
             \ar[dr]^-{{\cal M}}
             \\
             &
             &{\bcal V}(FU',HU')(q',\check{q'}) \otimes_2 {\bcal V}(FU,FU')(Ff,Fg)
             \ar[rd]^-{{\cal M}}
             \ar[ru]^{K\otimes_2 K}
             &&{\bcal W}(KFU,KHU')(Kq'KFf,K\check{q'}KFg)
             \ar@{=}[d]
           \\
           &
           I \otimes_2 {\bcal U}(U,U')(f,g)
             \ar[ru]_{\text{ }\mu_{U'} \otimes_2 F_{UU'_{fg}}}
             \ar `l[uu] `[ruu][ruu]^{(K\mu)_{U'}\otimes_2 KF}
           &&{\bcal V}(FU,HU')(q'Ff,\check{q'}Fg)
           \ar@{=}@/^1pc/[dd]
           \ar[r]^-{K}
           &{\bcal W}(KFU,KHU')(Kq'Ff,K\check{q'}Fg)
           \ar@{=}@/^1pc/[dd]
           \\
           {\bcal U}(U,U')(f,g)
             \ar[ru]^{=}
           \ar[rd]_{=}
           \\
           &
           {\bcal U}(U,U')(f,g) \otimes_2 I
             \ar[rd]^{\text{ }H_{UU'_{fg}} \otimes_2 \mu_{U}}
             \ar `l[dd] `[rdd][rdd]_{KH\otimes_2 (K\mu)_{U}}
           &&{\bcal V}(FU,HU')(Hfq,Hg\check{q})
           \ar[r]^-{K}
           &{\bcal W}(KFU,KHU')(KHfq,KHg\check{q})
           \\
             &
             &{\bcal V}(HU,HU')(Hf,Hg) \otimes_2 {\bcal V}(FU,HU)(q,\check{q})
             \ar[ru]^-{{\cal M}}
             \ar[rd]_{K\otimes_2 K}
             &&{\bcal W}(KFU,KHU')(KHfKq,KHgK\check{q})
             \ar@{=}[u]
             \\
             &&&{\bcal W}(KHU,KHU')(KHf,KHg) \otimes_2 {\bcal W}(KFU,KHU)(Kq,K\check{q})
             \ar[ur]_-{{\cal M}}
             }
           $$
           }
  \end{center}
  %\clearpage
  
  Secondly the left whiskering of a ${\cal V}$--2--functor onto a 
  ${\cal V}$--modification as in the picture:
  $$
       \xymatrix@R=10pt@C=20pt{
       &
       &
       &
       &&\ar@/_.5pc/@{=>}[dd]_{\gamma}
       &&\ar@/^.5pc/@{=>}[dd]^{\rho}
       \\
       {\bcal U}
       \ar[rrrr]^F
       &
       &&&{\bcal V}
       \ar@/^2pc/[rrrr]^G
       \ar@/_2pc/[rrrr]_K
       &
       \ar@3{->}[rr]_{\nu}
       &&&{\bcal W}\\
       &&&&&&&&
       }
       \mapsto
            \xymatrix@R=10pt@C=20pt{
                 &\ar@/_.5pc/@{=>}[dd]_{\gamma F}
                 &&\ar@/^.5pc/@{=>}[dd]^{\rho F}
                 &
                 \\
                 {\bcal U}
                 \ar@/^2pc/[rrrr]^{GF}
                 \ar@/_2pc/[rrrr]_{KF}
                 &
                 \ar@3{->}[rr]_{\nu F}
                 &&&{\bcal W}
                 \\
            &&&&
       }
    $$
  is given by $(\nu F)_U = \nu_{FU}.$
  
  This one is a valid ${\cal V}$--modification because of the following diagram, where we use the fact that $\nu$ is 
  a ${\cal V}$--modification.
  
\noindent
      		          \begin{center}
	          \resizebox{6.55in}{!}{ %!!!!
	          %%%%\resizebox{7.45in}{!}{
  $$
             \xymatrix@C=-33pt@R=40pt{
             &
             &&{\bcal W}(GFU',KFU')(\hat{q'},\overline{q'}) \otimes_2 {\bcal W}(GFU,GFU)(GFf,GFg)
             \ar[d]^-{{\cal M}}
           \\
           &
           I\otimes_2 {\bcal U}(U,U')(f,g)
           \ar[r]^-{1\otimes_2 F_{UU'_{fg}}}
           \ar@/^1.4pc/[rru]|-{(\nu F)_U' \otimes_2 GF}
           &
           I \otimes_2 {\bcal V}(FU,FU')(Ff,Fg)
             \ar[ru]|-{\text{ }\nu_{FU'} \otimes_2 G_{FUFU'_{FfFg}}}
           &\text{ }\text{ }\text{ }\text{ }{\bcal W}(GFU,KFU')(\hat{q'}GFf,\overline{q'}GFg)
           \ar@{=}[dd]
           &
           \\
           \text{ }\text{ }\text{ }
           \ar@<4pt>[ru]^{=}
           \ar@<-4pt>[rd]^{=}
           &{\bcal U}(U,U')(f,g)\text{---}{}^{F_{UU'_{fg}}}\to{\bcal V}(FU,FU')(Ff,Fg)
           \ar[ru]_{=}
           \ar[rd]^{=}
           &
           \\
           &
           {\bcal U}(U,U')(f,g)\otimes_2 I
  	 \ar[r]^-{F_{UU'_{fg}}\otimes_2 1}
           \ar@/_1.4pc/[rrd]|-{KF \otimes_2 (\nu F)_U}
           &
           {\bcal V}(FU,FU')(Ff,Fg) \otimes_2 I
             \ar[rd]|-{\text{ }K_{FUFU'_{FfFg}} \otimes_2 \nu_{FU}}
           &\text{ }\text{ }\text{ }\text{ }{\bcal W}(GFU,KFU)(KFf\hat{q},KFg\overline{q})
           &
           \\
             &&
             &{\bcal W}(KFU,KFU')(KFf,KFg) \otimes_2 {\bcal W}(GFU,KFU)(\hat{q},\overline{q})
             \ar[u]^-{{\cal M}}
             }
           $$
           }
  \end{center}
  
  The functoriality of these partials is shown just as for the whiskering of ${\cal V}$--2--functors onto 
  ${\cal V}$--2--natural transformations. 
  Consider a ${\cal V}$--modification $\xi:\phi\to\psi:T\to F:{\bcal U}\to {\bcal V}.$ 
  For right whiskering we have that $K(\mu*\xi)= K\mu*K\xi$ by
  the ${\cal V}$--2--functoriality of K. For a ${\cal V}$--2--functor $S:{\bcal X}\to {\bcal U}$
  we have $(\mu*\xi)S = \mu S * \xi S$ since $((\mu*\xi)S)_{X} = (\mu*\xi)_{SX} = (\mu S* \xi S)_{X}.$
  
  In the next step we 
  basically see the generalizations of these last two compositions.
  Next I define the right whiskering of a ${\cal V}$--2--natural transformation onto a 
  ${\cal V}$--modification as in the picture:
  
  $$
       \xymatrix@R=10pt@C=20pt{
       &\ar@/_.5pc/@{=>}[dd]_{\alpha}
       &&\ar@/^.5pc/@{=>}[dd]^{\beta}
       &&
       &\ar@{=>}[dd]^{\rho}
       &
       \\
       {\bcal U}
       \ar@/^2pc/[rrrr]^F
       \ar@/_2pc/[rrrr]_H
       &
       \ar@3{->}[rr]_{\mu}
       &&&{\bcal V}
       \ar@/^2pc/[rrrr]^G
       \ar@/_2pc/[rrrr]_K
       &
       &&&{\bcal W}\\
       &&&&&&&&
       }
       \mapsto
            \xymatrix@R=10pt@C=20pt{
                 &\ar@/_.5pc/@{=>}[dd]_{\rho\alpha}
                 &&\ar@/^.5pc/@{=>}[dd]^{\rho\beta}
                 &
                 \\
                 {\bcal U}
                 \ar@/^2pc/[rrrr]^{GF}
                 \ar@/_2pc/[rrrr]_{KH}
                 &
                 \ar@3{->}[rr]_{\rho\mu}
                 &&&{\bcal W}
                 \\
            &&&&
       }
    $$

  The ${\cal V}$--modification $\rho\mu:\rho\alpha \to \rho\beta$ can be defined in two ways.
  Let $\alpha_U(0)= q$, $\beta_{U}(0)= \check{q}$, $\rho_{FU}(0) = \overline{q}$, $\beta_{U}(0)= \check{q'}$ and $\rho_{HU}(0) = \overline{q'}.$
  The first way of composing is given by:
     $$\xymatrix{\rho\mu_U} =
           \xymatrix{
           I = I \otimes_2 I
           \ar[d]^{(K\mu)_U \otimes_2 (\rho F)_{U_{00}}}
           \\
           {\bcal W}(KFU,KHU)(K{q},K\check{q}) \otimes_2 {\bcal W}(GFU,KFU)(\overline{q},\overline{q})
           \ar[d]^{{\cal M}}
           \\
           {\bcal W}(GFU,KHU)(K{q}\overline{q},K\check{q}\overline{q})
           }
           $$
  %\clearpage
  The second is given by
  $$\xymatrix{\rho\mu_U} =
           \xymatrix{
           I = I \otimes_2 I
           \ar[d]^{(\rho H)_{U_{00}} \otimes_2 (G\mu)_U}
           \\
           {\bcal W}(GHU,KHU)(\overline{q'},\overline{q'}) \otimes_2 {\bcal W}(GFU,GHU)(G{q},G\check{q})
           \ar[d]^{{\cal M}}
           \\
           {\bcal W}(GFU,KHU)(\overline{q'}G{q},\overline{q'}G\check{q})
           }
           $$
  
  That the two ways agree is given by the following commuting diagram, which depends on the fact that $\rho$ is a 
  ${\cal V}$--2--natural transformation.
  
         \noindent
	       		          \begin{center}
	          \resizebox{6.55in}{!}{ %!!!!
	          %%%%\resizebox{7.45in}{!}{
         $$
             \xymatrix@C=-35pt@R=25pt{
             &
             &&{\bcal W}(GHU,KHU)(\overline{q'},\overline{q'}) \otimes_2 {\bcal W}(GFU,GHU)(G{q},G\check{q})
             \ar[d]^-{{\cal M}}
           \\
           &
           I\otimes_2 I
           \ar[r]_-{1\otimes_2 \mu_{U}}
           \ar@/^1.4pc/[rru]|-{(\rho H)_U \otimes_2 (G\mu)_U}
           &
           I \otimes_2 {\bcal V}(FU,HU)(q,\check{q})
             \ar[ru]|-{\text{ }\rho_{HU} \otimes_2 G_{FUHU_{q\check{q}}}}
           &{\bcal W}(GFU,KHU)(\overline{q'}G{q},\overline{q'}G\check{q})
           \ar@{=}[dd]
           &
           \\
           \ar[ru]^{=}
           \ar[rd]^{=}
           &\text{ }\text{ }\text{ }\text{ }\text{ }\text{ }\text{ }\text{ }I\text{---}{}^{\mu_{U}}\to{\bcal V}(FU,HU)(q,\check{q})\text{ }\text{ }\text{ }\text{ }\text{ }\text{ }\text{ }\text{ }\text{ }\text{ }\text{ }\text{ }\text{ }
             \ar[ru]_{=}
           \ar[rd]^{=}
           \\
           &
           I\otimes_2 I
  	 \ar[r]^-{\mu_{U}\otimes_2 1}
           \ar@/_1.4pc/[rrd]|-{(K\mu)_U \otimes_2 (\rho F)_U}
           &
           {\bcal V}(FU,HU)(q,\check{q}) \otimes_2 I
             \ar[rd]|-{\text{ }K_{FUHU_{q\check{q}}} \otimes_2 \rho_{FU}}
           &{\bcal W}(GFU,KHU)(K{q}\overline{q},K\check{q}\overline{q})
           &
           \\
             &&
             &{\bcal W}(KFU,KHU)(K{q},K\check{q}) \otimes_2 {\bcal W}(GFU,KFU)(\overline{q},\overline{q})
             \ar[u]^-{{\cal M}}
             }
           $$
           }
  \end{center}
  
  That this composition yields a ${\cal V}$--modification is easily seen when we note that it is by definition
  the same as composing certain ${\cal V}$--modifications along a common ${\cal V}$--2--functor. For example the
  above composition is of the ${\cal V}$--modifications $K\mu$ and $1_{\rho F}$ along the ${\cal V}$--2--functor
  $KF.$
  
  Now we can define the left whiskering of a ${\cal V}$--2--natural transformation onto a 
  ${\cal V}$--modification as in the picture:
  $$
       
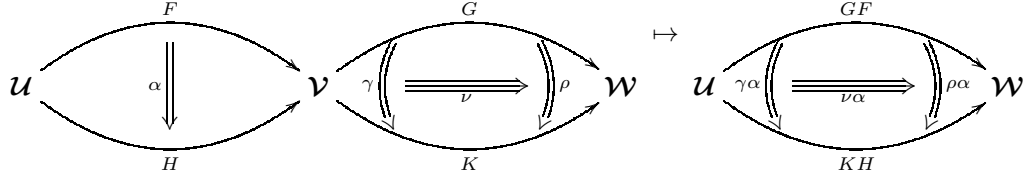
\xymatrix@R=10pt@C=20pt{
       &
       &\ar@{=>}[dd]_{\alpha}
       &
       &&\ar@/_.5pc/@{=>}[dd]_{\gamma}
       &&\ar@/^.5pc/@{=>}[dd]^{\rho}
       \\
       {\bcal U}
       \ar@/^2pc/[rrrr]^F
       \ar@/_2pc/[rrrr]_H
       &
       &&&{\bcal V}
       \ar@/^2pc/[rrrr]^G
       \ar@/_2pc/[rrrr]_K
       &
       \ar@3{->}[rr]_{\nu}
       &&&{\bcal W}\\
       &&&&&&&&
       }
       \mapsto
                 \xymatrix@R=10pt@C=20pt{
                      &\ar@/_.5pc/@{=>}[dd]_{\gamma\alpha}
                      &&\ar@/^.5pc/@{=>}[dd]^{\rho\alpha}
                      &
                      \\
                      {\bcal U}
                      \ar@/^2pc/[rrrr]^{GF}
                      \ar@/_2pc/[rrrr]_{KH}
                      &
                      \ar@3{->}[rr]_{\nu\alpha}
                      &&&{\bcal W}
                      \\
                 &&&&
       }
    $$
  The ${\cal V}$--modification $\nu\alpha:\gamma\alpha \to \rho\alpha$ can be defined in two ways.
  Let $\alpha_U(0)= q$, $\gamma_{FU}(0)= \hat{q}$, $\rho_{FU}(0) = \overline{q}$, $\gamma_{HU}(0)= \hat{q'}$ and $\rho_{HU}(0) = \overline{q'}.$
  The first way of composing is given by:
     $$\xymatrix{\nu\alpha_U} =
           \xymatrix{
           I = I \otimes_2 I
           \ar[d]^{(\nu H)_U \otimes_2 (G\alpha)_{U_{00}}}
           \\
           {\bcal W}(GHU,KHU)(\hat{q'},\overline{q'}) \otimes_2 {\bcal W}(GFU,GHU)(Gq,Gq)
           \ar[d]^{{\cal M}}
           \\
           {\bcal W}(GFU,KHU)(\hat{q'}Gq,\overline{q'}Gq)
           }
           $$
  %\clearpage
  The second is given by
  $$\xymatrix{\nu\alpha_U} =
           \xymatrix{
           I = I \otimes_2 I
           \ar[d]^{(K\alpha)_{U_{00}} \otimes_2 (\nu F)_U}
           \\
           {\bcal W}(KFU,KHU)(Kq,Kq) \otimes_2 {\bcal W}(GFU,KFU)(\hat{q},\overline{q})
           \ar[d]^{{\cal M}}
           \\
           {\bcal W}(GFU,KHU)(Kq\hat{q},Kq\overline{q})
           }
           $$
  
  That the two ways agree is given by the following commuting diagram, which depends on the fact that $\nu$ is a 
  ${\cal V}$--modification.
  
          \noindent
	        		          \begin{center}
	          \resizebox{6.55in}{!}{ %!!!!
	          %%%%\resizebox{7.45in}{!}{
          $$
             \xymatrix@C=-35pt@R=30pt{
             &
             &&{\bcal W}(GHU,KHU)(\hat{q'},\overline{q'}) \otimes_2 {\bcal W}(GFU,GHU)(Gq,Gq)
             \ar[d]^-{{\cal M}}
           \\
           &
           I\otimes_2 I
           \ar[r]_-{1\otimes_2 \alpha_{U_{00}}}
           \ar@/^1.4pc/[rru]|-{(\nu H)_U \otimes_2 (G\alpha)_U}
           &
           I \otimes_2 {\bcal V}(FU,HU)(q,q)
             \ar[ru]|-{\text{ }\nu_{HU} \otimes_2 G_{FUHU_{qq}}}
           &{\bcal W}(GFU,KHU)(\hat{q'}Gq,\overline{q'}Gq)
           \ar@{=}[dd]
           &
           \\
           %I
           %\ar[r]^{\alpha_{U_{00}}}
           \ar[ru]^{=}
           \ar[rd]^{=}
           &\text{ }\text{ }\text{ }\text{ }\text{ }\text{ }\text{ }\text{ }I\text{---}{}^{\alpha_{U_{00}}}\to{\bcal V}(FU,HU)(q,q)\text{ }\text{ }\text{ }\text{ }\text{ }\text{ }\text{ }\text{ }\text{ }\text{ }\text{ }\text{ }\text{ }
             \ar[ru]_{=}
           \ar[rd]^{=}
           \\
           &
           I\otimes_2 I
  	 \ar[r]^-{\alpha_{U_{00}}\otimes_2 1}
           \ar@/_1.4pc/[rrd]|-{(K\alpha)_U \otimes_2 (\nu F)_U}
           &
           {\bcal V}(FU,HU)(q,q) \otimes_2 I
             \ar[rd]|-{\text{ }K_{FUHU_{qq}} \otimes_2 \nu_{FU}}
           &{\bcal W}(GFU,KHU)(Kq\hat{q},Kq\overline{q})
           &
           \\
             &&
             &{\bcal W}(KFU,KHU)(Kq,Kq) \otimes_2 {\bcal W}(GFU,KFU)(\hat{q},\overline{q})
             \ar[u]^-{{\cal M}}
             }
           $$
           }
  \end{center}
  
  Again this composition yields a ${\cal V}$--modification since it is by definition
  the same as composing certain ${\cal V}$--modifications along a common ${\cal V}$--2--functor.
  
  Necessary for the functoriality of the partials given by the above left and right whiskering is that we have
  $\rho(\omega \circ \mu) = \rho\omega \circ \rho\mu$ and $(tau \circ \nu)\alpha = (\tau\alpha \circ \nu\alpha)$
  as in the following
  pictures:
   \noindent
  	        		          \begin{center}
  	          \resizebox{6.55in}{!}{ %!!!!
	          %%%%\resizebox{7.45in}{!}{
  $$
       \xymatrix@R=10pt@C=20pt{
       &\ar@/_.5pc/@{=>}[dd]_{\alpha}
       &\ar@{=>}[dd]^<<<{\beta}
       &\ar@/^.5pc/@{=>}[dd]^{\delta}
       &&
       &\ar@{=>}[dd]^{\rho}&
       \\
       {\bcal U}
       \ar@/^2pc/[rrrr]^F
       \ar@/_2pc/[rrrr]_H
       &
       \ar@3{->}[r]_{\mu}
       &\ar@3{->}[r]_{\omega}
       &&{\bcal V}
       \ar@/^2pc/[rrrr]^G
       \ar@/_2pc/[rrrr]_K
       &&&&{\bcal W}\\
       &&&&&&&&
       }
       \text{ and }
       \xymatrix@R=10pt@C=20pt{
            &
            &\ar@{=>}[dd]_{\alpha}
            &
            &&\ar@/_.5pc/@{=>}[dd]_{\gamma}
            &\ar@{=>}[dd]^<<<{\rho}
            &\ar@/^.5pc/@{=>}[dd]^{\sigma}
            \\
            {\bcal U}
            \ar@/^2pc/[rrrr]^F
            \ar@/_2pc/[rrrr]_H
            &          
            &&&{\bcal V}
            \ar@/^2pc/[rrrr]^G
            \ar@/_2pc/[rrrr]_K
            &
            \ar@3{->}[r]_{\nu}
            &\ar@3{->}[r]_{\tau}
            &&{\bcal W}\\
            &&&&&&&&
            }
    $$
    }
    \end{center}
  These requirements are met since the exteriors of the following two diagrams commute, respectively.
  Here let $\delta_U(0)= \hat{\hat{q}}$ and $\tau_U(0) = \overline{\overline{q}}.$

  \noindent
      		          \begin{center}
	          \resizebox{6.25in}{!}{  
  \begin{sideways}
        \begin{footnotesize}
          $$
          \xymatrix@R=30pt@C=-8pt{
          &&&&{\bcal W}(GFU,KHU)(K\check{q}\overline{q},K\hat{\hat{q}}\overline{q})\otimes_1 {\bcal W}(GFU,KHU)(Kq\overline{q},K\check{q}\overline{q})
          \ar@<5.5ex>@/^15pc/[ddddddrrrr]^{M}
          \\
          &&&{}\save[]*\txt{\text{ }\\$({\bcal W}(KFU,KHU)(K\check{q},K\hat{\hat{q}}) \otimes_2 {\bcal W}(GFU,KFU)(\overline{q},\overline{q}))\otimes_1 ({\bcal W}(KFU,KHU)(Kq,K\check{q}) \otimes_2 {\bcal W}(GFU,KFU)(\overline{q},\overline{q}))$\\\text{ }\\\text{ }}
          \ar[ur]^{{\cal M}\otimes_1 {\cal M}}
          %\ar[ddr]^{\eta}
          \restore
          \\
          &&&&&
          \\
          &&&&{}\save[]*\txt{\text{ }\\$({\bcal W}(KFU,KHU)(K\check{q},K\hat{\hat{q}}) \otimes_1 {\bcal W}(KFU,KHU)(Kq,K\check{q}))\otimes_2 ({\bcal W}(GFU,KFU)(\overline{q},\overline{q}) \otimes_1 {\bcal W}(GFU,KFU)(\overline{q},\overline{q}))$\\\text{ }\\\text{ }}
          \ar[dddr]^{M\otimes_2 M}
          \ar@{<-}[uul]^{\eta}
          \restore
          \\\\\\
          &\text{ }\text{ }\text{ }\text{ }\text{ }\text{ }\text{ }\text{ }\text{ }\text{ }\text{ }\text{ }\text{ }\text{ }\text{ }\text{ }\text{ }
          &&I=I\otimes_1 I = (I\otimes_2 I)\otimes_1 (I\otimes_2 I)\text{---------}{}^{\eta}\to (I\otimes_1 I)\otimes_2 (I\otimes_1 I)=I\otimes_2 I
          %\ar@/_1pc/[rr]^-{K(\omega\circ\mu)\otimes_2 \rho F}
          \ar@<-17ex>[uuu]|{(K\omega \otimes_2 \rho F)\otimes_1 (K\mu \otimes_2 \rho F)}
          \ar@<12ex>[uuuuu]|{(K\omega \otimes_1 K\mu)\otimes_2 (\rho F \otimes_1 \rho F)}
          &&{}\save[]*\txt{\text{ }\\\text{ }\\${\bcal W}(KFU,KHU)(Kq,K\hat{\hat{q}}) \otimes_2 {\bcal W}(GFU,KFU)(\overline{q},\overline{q})\text{--}{}^{\cal M}\to {\bcal W}(GFU,KHU)(Kq\overline{q},K\hat{\hat{q}}\overline{q})\text{ }\text{ }\text{ }\text{ }\text{ }\text{ }$\\\text{ }}
          \ar@{<-}[ll]^-{K(\omega\circ\mu)\otimes_2 \rho F}
          \restore
          &\text{ }\text{ }\text{ }\text{ }\text{ }\text{ }\text{ }\text{ }&\text{ }\text{ }\text{ }\text{ }\text{ }\text{ }\text{ }&\text{ }\text{ }
          \\\\
          &\text{ }\text{ }\text{ }\text{ }\text{ }\text{ }\text{ }\text{ }\text{ }\text{ }\text{ }\text{ }\text{ }\text{ }\text{ }\text{ }\text{ }
            &&I=I\otimes_1 I = (I\otimes_2 I)\otimes_1 (I\otimes_2 I)\text{---------}{}^{\eta}\to (I\otimes_1 I)\otimes_2 (I\otimes_1 I)=I\otimes_2 I
            %\ar@/^1pc/[rr]_-{K\alpha\otimes_2 (\tau\circ\nu)F}
            \ar@<17ex>[ddd]|{( K\alpha\otimes_1 K\alpha)\otimes_2 (\tau F \otimes_1 \nu F)}
            \ar@<-12ex>[ddddd]|{( K\alpha\otimes_2 \tau F)\otimes_1 (K\alpha \otimes_2 \nu F)}
            &&{}\save[]*\txt{\text{ }\\${\bcal W}(KFU,KHU)(Kq,Kq) \otimes_2 {\bcal W}(GFU,KFU)(\hat{q},\bar{\bar{q}})\text{--}{}^{\cal M}\to {\bcal W}(GFU,KHU)(Kq\hat{q},Kq\bar{\bar{q}})\text{ }\text{ }\text{ }\text{ }\text{ }\text{ }$\\\text{ }\\\text{ }}
            \ar@{<-}[ll]^-{K\alpha\otimes_2 (\tau\circ\nu)F}
            \restore  
            &\text{ }\text{ }\text{ }\text{ }\text{ }\text{ }\text{ }\text{ }&\text{ }\text{ }\text{ }\text{ }\text{ }\text{ }\text{ }&\text{ }\text{ }
            \\\\\\
            &&&&{}\save[]*\txt{\text{ }\\\text{ }\\$({\bcal W}(KFU,KHU)(Kq,Kq) \otimes_1 {\bcal W}(KFU,KHU)(Kq,Kq))\otimes_2 ({\bcal W}(GFU,KFU)(\bar{q},\bar{\bar{q}}) \otimes_1 {\bcal W}(GFU,KFU)(\hat{q},\bar{q}))$\\\text{ }}
            \ar[uuur]^{M\otimes_2 M}
            \ar@{<-}[ddl]^{\eta}
            \restore
            \\
            &&&&&
           \\
            &&&{}\save[]*\txt{\text{ }\\\text{ }\\$({\bcal W}(KFU,KHU)(Kq,Kq) \otimes_2 {\bcal W}(GFU,KFU)(\bar{q},\bar{\bar{q}}))\otimes_1 ({\bcal W}(KFU,KHU)(Kq,Kq) \otimes_2 {\bcal W}(GFU,KFU)(\hat{q},\bar{q}))$\\\text{ }}
            %\ar[uur]^{\eta}
            \ar[dr]^{{\cal M}\otimes_1 {\cal M}}
            \restore
            \\
            &&&&{\bcal W}(GFU,KHU)(Kq\bar{q},Kq\bar{\bar{q}}) \otimes_1 {\bcal W}(GFU,KHU)(Kq\hat{q},Kq\bar{q})
            \ar@<-5.5ex>@/_15pc/[uuuuuurrrr]_{M}
            }
        $$
          \end{footnotesize}
          \end{sideways}
  }
  \end{center}
   
   %\clearpage
  %\newpage
  
  It is straightforward to check that if an identity ${\cal V}$--2--natural transformation for a given ${\cal V}$--2--functor
  is whiskered onto the left or 
  right of a ${\cal V}$--modification the definitions give exactly the respective whiskering of the 
  ${\cal V}$--2--functor itself. Thus the following definition of the horizontal composition of ${\cal V}$--modifications 
  along a ${\cal V}$--2--category, given in terms of composing along a common ${\cal V}$--2--natural transformation
  could be written less generally but equivalently in terms of composing along a common 
  ${\cal V}$--2--functor. Either way the result is a valid ${\cal V}$--modification based on an earlier proof.
  The equivalence will actually be a corollary of the proof of the well defined nature of the composition.
  Now we are considering the full picture:
  $$
       
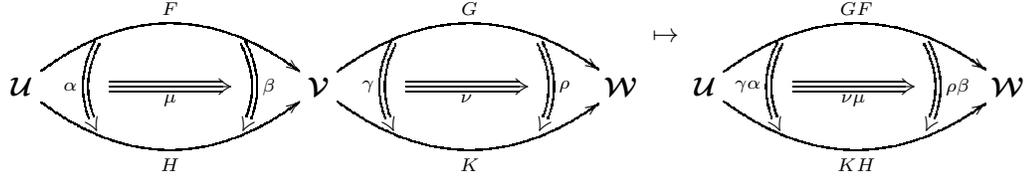
\xymatrix@R=10pt@C=20pt{
       &\ar@/_.5pc/@{=>}[dd]_{\alpha}
       &&\ar@/^.5pc/@{=>}[dd]^{\beta}
       &&\ar@/_.5pc/@{=>}[dd]_{\gamma}
       &&\ar@/^.5pc/@{=>}[dd]^{\rho}
       \\
       {\bcal U}
       \ar@/^2pc/[rrrr]^F
       \ar@/_2pc/[rrrr]_H
       &
       \ar@3{->}[rr]_{\mu}
       &&&{\bcal V}
       \ar@/^2pc/[rrrr]^G
       \ar@/_2pc/[rrrr]_K
       &
       \ar@3{->}[rr]_{\nu}
       &&&{\bcal W}\\
       &&&&&&&&
       }
       \mapsto
                 \xymatrix@R=10pt@C=20pt{
                      &\ar@/_.5pc/@{=>}[dd]_{\gamma\alpha}
                      &&\ar@/^.5pc/@{=>}[dd]^{\rho\beta}
                      &
                      \\
                      {\bcal U}
                      \ar@/^2pc/[rrrr]^{GF}
                      \ar@/_2pc/[rrrr]_{KH}
                      &
                      \ar@3{->}[rr]_{\nu\mu}
                      &&&{\bcal W}
                      \\
                 &&&&
       }
    $$
  and the two ways of defining $\nu\mu$ in terms of composing along a common ${\cal V}$--2--natural transformation
  are as follow. 
  
  The first is:
  $$\xymatrix{(\nu\mu)_U = (\rho\mu \circ \nu\alpha)_U} =$$

\noindent
        		          \begin{center}
	          \resizebox{6.55in}{!}{ %!!!!
	          %%%%\resizebox{7.45in}{!}{
             $$
           \xymatrix{
           I = I \otimes_1 I = (I \otimes_2 I) \otimes_1 (I \otimes_2 I)
           \ar[d]^{((K\mu)_U \otimes_2 (\rho F)_{U_{00}})\otimes_1 ((K\alpha)_{U_{00}} \otimes_2 (\nu F)_U)}
           \\
           {}\save[]*\txt{\text{ }\\\text{ }\\$({\bcal W}(KFU,KHU)(Kq,K\check{q}) \otimes_2 {\bcal W}(GFU,KFU)(\overline{q},\overline{q})\otimes_1({\bcal W}(KFU,KHU)(Kq,Kq) \otimes_2 {\bcal W}(GFU,KFU)(\hat{q},\overline{q}))$\\\text{ }}
           \ar[d]^{{\cal M}\otimes_1 {\cal M}} \restore
           \\
           {\bcal W}(GFU,KHU)(Kq\overline{q},K\check{q}\overline{q})\otimes_1 {\bcal W}(GFU,KHU)(Kq\hat{q},Kq\overline{q})
           \ar[d]^{M}
           \\
           {\bcal W}(GFU,KHU)(Kq\hat{q},K\check{q}\overline{q})
           }
           $$
           }
  \end{center}
           
  and the second:
  $$\xymatrix{(\nu\mu)_U = (\nu\beta \circ \gamma\mu)_U} =$$ 
            
            \noindent
        		          \begin{center}
	          \resizebox{6.55in}{!}{ %!!!!
	          %%%%\resizebox{7.45in}{!}{
	          $$
           \xymatrix{
           I = I \otimes_1 I = (I \otimes_2 I) \otimes_1 (I \otimes_2 I)
           \ar[d]^{((\nu H)_U \otimes_2 (G\beta)_{U_{00}})\otimes_1 ((\gamma H)_{U_{00}} \otimes_2 (G\mu)_U)}
           \\
           {}\save[]*\txt{\text{ }\\\text{ }\\$({\bcal W}(GHU,KHU)(\hat{q'},\overline{q'}) \otimes_2 {\bcal W}(GFU,GHU)(G\check{q},G\check{q}))\otimes_1({\bcal W}(GHU,KHU)(\hat{q'},\hat{q'}) \otimes_2 {\bcal W}(GFU,GHU)(G{q},G\check{q}))$\\\text{ }}
           \ar[d]^{{\cal M}\otimes_1 {\cal M}} \restore
           \\
           {\bcal W}(GFU,KHU)(\hat{q'}G\check{q},\overline{q'}G\check{q})\otimes_1 {\bcal W}(GFU,KHU)(\hat{q'}G{q},\hat{q'}G\check{q})
           \ar[d]^{M}
           \\
           {\bcal W}(GFU,KHU)(\hat{q'}G{q},\overline{q'}G\check{q})
           }
           $$
           }
  \end{center}
  
  Note that in both of the preceding two definitions we have made two choices between equivalent ways of representing
  component ${\cal V}$--modifications. The preceding two definitions are equivalent based on the following commutative
  diagram.
  \clearpage
  
  \noindent
  \begin{center}
  \begin{list}{}{}
  \noindent
  \centering
	          \resizebox{6.25in}{!}{ %!!!!
	          %%%% \resizebox{7.35in}{!}{
  $
  ({\bcal W}(GHU,KHU)(\hat{q'},\overline{q'}) \otimes_2 {\bcal W}(GFU,GHU)(G\check{q},G\check{q}))\otimes_1({\bcal W}(GHU,KHU)(\hat{q'},\hat{q'}) \otimes_2 {\bcal W}(GFU,GHU)(G{q},G\check{q})) 
  $
  }

   \noindent
  \centering
	          \resizebox{6.25in}{!}{ %!!!!
	          %%%%\resizebox{7.4in}{!}{
  $$%_{{\displaystyle \downarrow}}
  \xymatrix@C=-15pt@R=12pt{
  &&&&&&\bullet%{{{{\displaystyle \ddots}_{{}_{\displaystyle \searrow}}}_{}}_{}}_{\text{ }}%^{\displaystyle \downarrow}
  \ar[rrdddddd]^{{\cal M}\otimes_1 {\cal M}}
  \ar[ldddd]^{\eta}
  \\\\
  &&&(I\otimes_2 I)\otimes_1 (I\otimes_2 I)
  \ar[rrruu]|{((\nu H)_U \otimes_2 (G\beta)_{U_{00}})\otimes_1 ((\gamma H)_{U_{00}} \otimes_2 (G\mu)_U)}
  \ar[rddd]^{\eta}
  \\\\
  &&&&&\bullet
  \ar[rdddd]^{M\otimes_2 M}
  \\
  &&&&(I\otimes_1 I)\otimes_2 (I\otimes_1 I)
  \ar[ru]%^{((\nu H)_U \otimes_1 (\gamma H)_{U_{00}})\otimes_2 ((G\beta)_{U_{00}} \otimes_1 (G\mu)_U)}
  \\
  &&&&&&&&\bullet%{\bcal W}(GFU,KHU)(\hat{q'}G\check{q},\overline{q'}G\check{q})\otimes_1 {\bcal W}(GFU,KHU)(\hat{q'}G{q},\hat{q'}G\check{q})
  \ar[rdddd]^{M}
  \\\\
  &I\otimes_1 I
  \ar[rruuuuuu]^{=}
  &&&&&{\bcal W}(GHU,KHU)(\hat{q'},\overline{q'})\otimes_2 {\bcal W}(GFU,GHU)(Gq,G\check{q})
  \ar[rrrdd]^{{\cal M}}
  \ar@<.48pt>[rrrdd]
  \ar@<-.48pt>[rrrdd]
  \\
  &&I\otimes_2 I
  \ar[rruuuu]^{=}
  \ar@/^.5pc/[rrrru]^-{(\nu H)_{U}\otimes_2 (G\mu)_{U}\text{ }\text{ }}
  \ar@<.48pt>@/^.5pc/[rrrru]
  \ar@<-.48pt>@/^.5pc/[rrrru]
  \ar[rrd]|{1\otimes_2 \mu_{U}}
  \\
  &&&&I\otimes_2 {\bcal V}(FU,HU)(q,\check{q})
  \ar[rruu]|{\nu_{HU}\otimes_2 G_{FUHU}}
  &&&&&{\bcal W}(GFU,KHU)(\hat{q'}Gq,\overline{q'}G\check{q})
  \ar@{=}@/^1pc/[dddd]
  \ar@<.32pt>@{=}@/^1pc/[dddd]
  \ar@<-.25pt>@{=}@/^1pc/[dddd]
  \\\\
  I
  \ar[ruuuu]^{=}
  \ar[rruuu]_{=}
  \ar@<.48pt>[rruuu]
  \ar@<-.48pt>[rruuu]
  \ar[rrr]^-{\mu_U}
  \ar[rdddd]_{=}
  \ar[rrddd]^{=}
  \ar@<.48pt>[rrddd]
  \ar@<-.48pt>[rrddd]
  &&&{\bcal V}(FU,HU)(q,\check{q})
  \ar[ruu]^{=}
  \ar[rdd]_{=}
  \\\\
  &&&&{\bcal V}(FU,HU)(q,\check{q})\otimes_2 I
  \ar[rrdd]|{K_{FUHU}\otimes_2 \nu_{FU}}
  &&&&&{\bcal W}(GFU,KHU)(Kq\hat{q},K\check{q\overline{q}})
  \\
  &&I\otimes_2 I
  \ar[rrdddd]^{=}
  \ar@/_.5pc/[rrrrd]_-{(K\mu)_{U}\otimes_2 (\nu F)_{U}\text{ }\text{ }}
  \ar@<.48pt>@/_.5pc/[rrrrd]
  \ar@<-.48pt>@/_.5pc/[rrrrd]
  \ar[rru]|{\mu_{U}\otimes_2 1}
  \\
  &I\otimes_1 I
  \ar[rrdddddd]^{=}
  &&&&&{\bcal W}(KFU,KHU)(Kq,K\check{q})\otimes_2 {\bcal W}(GFU,KFU)(\hat{q},\overline{q})
  \ar[rrruu]_{{\cal M}}
  \ar@<.48pt>[rrruu]
  \ar@<-.48pt>[rrruu]
  \\\\
  &&&&&&&&\bullet%{\bcal W}(GFU,KHU)(Kq\overline{q},K\check{q}\overline{q})\otimes_1 {\bcal W}(GFU,KHU)(Kq\hat{q},Kq\overline{q})
  \ar[ruuuu]^{M}
  \\
  &&&&(I\otimes_1 I)\otimes_2 (I\otimes_1 I)
  \ar[rd]%^{((K\mu)_U \otimes_1 (K\alpha)_{U_{00}})\otimes_2 ((\rho F)_{U_{00}} \otimes_1 (\nu F)_U)}
  \\
  &&&&&\bullet
  \ar[ruuuu]^{M\otimes_2 M}
  \\\\
  &&&(I\otimes_2 I)\otimes_1 (I\otimes_2 I)
  \ar[rrrdd]|{((K\mu)_U \otimes_2 (\rho F)_{U_{00}})\otimes_1 ((K\alpha)_{U_{00}} \otimes_2 (\nu F)_U)}
  \ar[ruuu]^{\eta}
  \\\\
  &&&&&&\bullet%{\text{ }}_{{}_{{}_{{\displaystyle \nwarrow}_{\displaystyle \ddots}}}}%\vdots_{\displaystyle \uparrow}
  \ar[rruuuuuu]^{{\cal M}\otimes_1 {\cal M}}
  \ar[luuuu]^{\eta}
  }
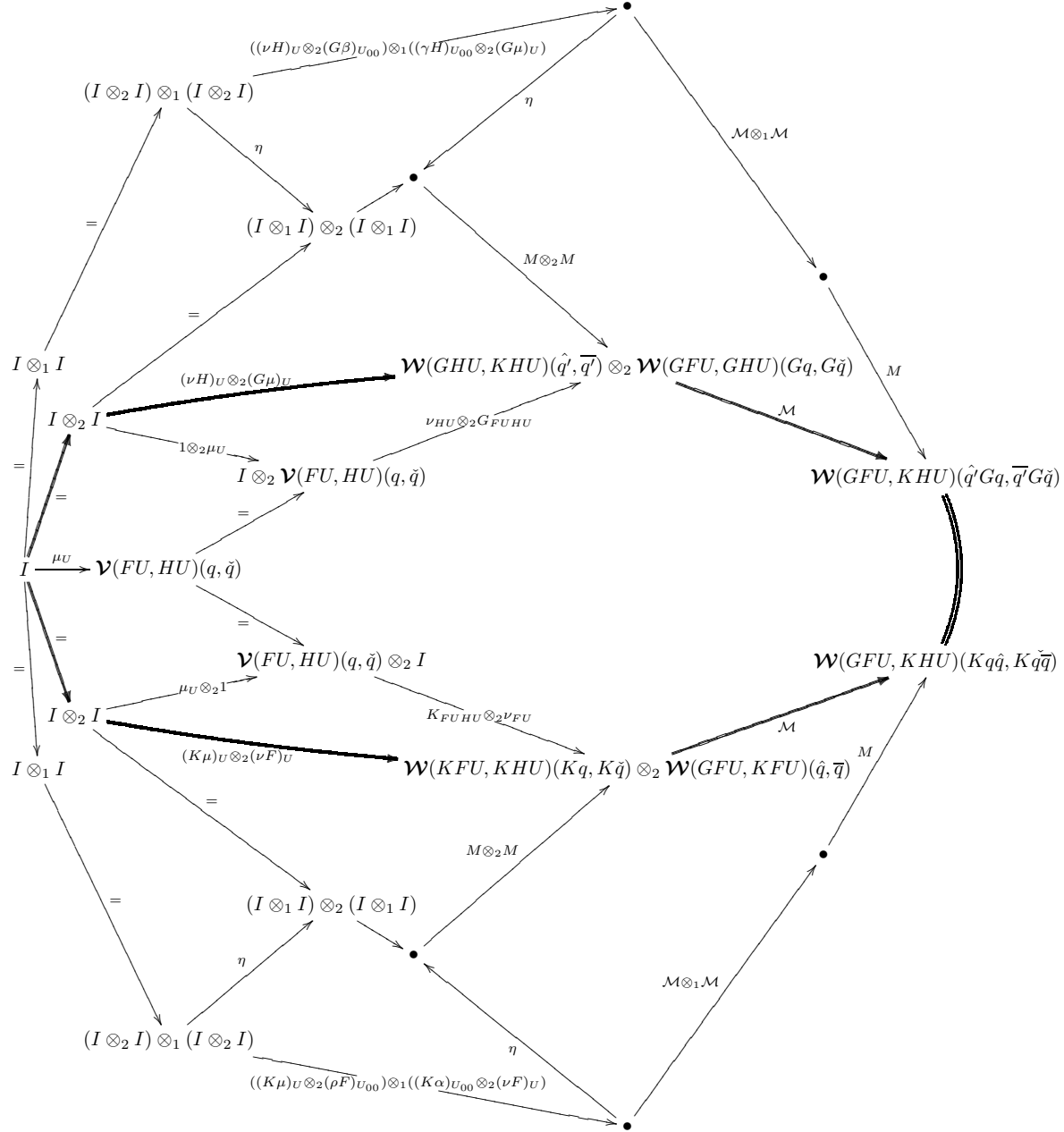
  $$%^{{\displaystyle \uparrow}}
  }
    
    \noindent
    \centering
	          \resizebox{6.25in}{!}{ %!!!!
	          %%%%\resizebox{7.35in}{!}{
  $
  ({\bcal W}(KFU,KHU)(Kq,K\check{q}) \otimes_2 {\bcal W}(GFU,KFU)(\overline{q},\overline{q})\otimes_1({\bcal W}(KFU,KHU)(Kq,Kq) \otimes_2 {\bcal W}(GFU,KFU)(\hat{q},\overline{q}))
  $
  }
   \end{list}
   \end{center}
  %{part3}

  \clearpage
  %\newpage

  The exterior commutes since all the interior regions commute. The top and bottom bullets are labeled by the
  text at the top and bottom of the diagram. The other bullets and
  unlabeled arrows should be easily filled in, noting that
  the uppermost and lowest quadrilaterals commute by the naturality of $\eta.$ 
  The arrows marked with an ``='' all occur as copies of $I$ are tensored to the object at the arrow's source. Therefore
  the western regions with the initial $I$ as a vertex all commute trivially. The large central region expresses the
  fact that $\nu$ is a ${\cal V}$--modification. The pentagonal regions on the right commute by
  the ${\cal V}$--functoriality of ${\cal M}.$ The remaining interior regions commute by definition and by the 
  axioms of a ${\cal V}$--category.
   
   The thick arrows in the central portion of the above diagram outline the definition of composing
  ${\cal V}$--modifications along a ${\cal V}$--2--category in terms of composing along a 
  common ${\cal V}$--2--functor. Thus this diagram also demonstrates that the two ways of doing so are
  equivalent to each other and to the method which uses composition along a 
  common ${\cal V}$--2--natural transformation.
  
  Next we continue checking functoriality of partials. As usual I check the stronger property that the 
  composition defined by those partials gives the whiskering itself as a composition with a unit. First we check that
  composing in the following two pictures yields the same ${\cal V}$--modification.
 
 \noindent
         		          \begin{center}
 	          \resizebox{6.55in}{!}
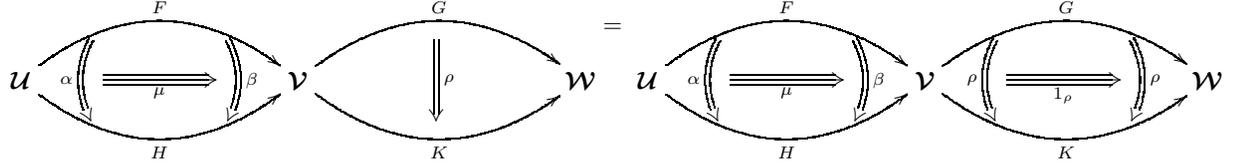
{ %!!!!
	          %%%%\resizebox{7.45in}{!}{
  $$
       \xymatrix@R=10pt@C=20pt{
       &\ar@/_.5pc/@{=>}[dd]_{\alpha}
       &&\ar@/^.5pc/@{=>}[dd]^{\beta}
       &&
       &\ar@{=>}[dd]^{\rho}
       &
       \\
       {\bcal U}
       \ar@/^2pc/[rrrr]^F
       \ar@/_2pc/[rrrr]_H
       &
       \ar@3{->}[rr]_{\mu}
       &&&{\bcal V}
       \ar@/^2pc/[rrrr]^G
       \ar@/_2pc/[rrrr]_K
       &
       &&&{\bcal W}\\
       &&&&&&&&
       } =
       \xymatrix@R=10pt@C=20pt{
            &\ar@/_.5pc/@{=>}[dd]_{\alpha}
            &&\ar@/^.5pc/@{=>}[dd]^{\beta}
            &&\ar@/_.5pc/@{=>}[dd]_{\rho}
            &&\ar@/^.5pc/@{=>}[dd]^{\rho}
            \\
            {\bcal U}
            \ar@/^2pc/[rrrr]^F
            \ar@/_2pc/[rrrr]_H
            &
            \ar@3{->}[rr]_{\mu}
            &&&{\bcal V}
            \ar@/^2pc/[rrrr]^G
            \ar@/_2pc/[rrrr]_K
            &
            \ar@3{->}[rr]_{1_{\rho}}
            &&&{\bcal W}\\
            &&&&&&&&
            }
    $$
 }
 \end{center}
 
 Using the definition of the composition of ${\cal V}$--modifications in terms of composing along a common
  ${\cal V}$--2--functor it is easy to see that this equality follows from the fact that 
  $(\rho F)_{U_{00}}= (\rho)_{FU_{00}} = j_{(\rho_{FU}(0))} = (1_{\rho})_{FU} = (1_{\rho}F)_U.$
  
  As noted earlier the compositions in the first two of the
  following pictures yield the same ${\cal V}$--modifications as well,
  due to the fact that $(1_G F)_{U_{00}} = (1_G)_{FU_{00}} = {\cal J}_{GFU_{00}}.$ Thus by the above equality
  all three are equivalent:
  
  $$
       
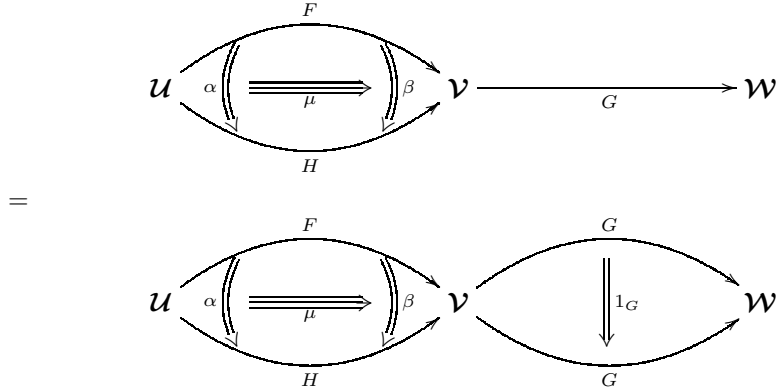
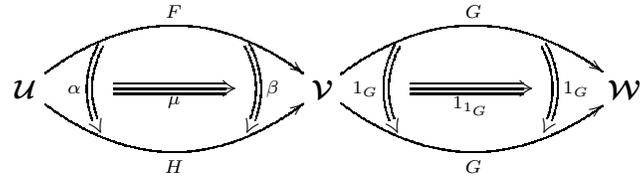
\xymatrix@R=10pt@C=20pt{
       &\ar@/_.5pc/@{=>}[dd]_{\alpha}
       &&\ar@/^.5pc/@{=>}[dd]^{\beta}
       &&
       &
       &
       \\
       {\bcal U}
       \ar@/^2pc/[rrrr]^F
       \ar@/_2pc/[rrrr]_H
       &
       \ar@3{->}[rr]_{\mu}
       &&&{\bcal V}
       \ar[rrrr]_G
       &
       &&&{\bcal W}\\
       &&&&&&&&
       }
       $$=$$
       \xymatrix@R=10pt@C=20pt{
       &\ar@/_.5pc/@{=>}[dd]_{\alpha}
       &&\ar@/^.5pc/@{=>}[dd]^{\beta}
       &&
       &\ar@{=>}[dd]^{1_G}
       &
       \\
       {\bcal U}
       \ar@/^2pc/[rrrr]^F
       \ar@/_2pc/[rrrr]_H
       &
       \ar@3{->}[rr]_{\mu}
       &&&{\bcal V}
       \ar@/^2pc/[rrrr]^G
       \ar@/_2pc/[rrrr]_G
       &
       &&&{\bcal W}\\
       &&&&&&&&
       }$$ = $$
       \xymatrix@R=10pt@C=20pt{
            &\ar@/_.5pc/@{=>}[dd]_{\alpha}
            &&\ar@/^.5pc/@{=>}[dd]^{\beta}
            &&\ar@/_.5pc/@{=>}[dd]_{1_G}
            &&\ar@/^.5pc/@{=>}[dd]^{1_G}
            \\
            {\bcal U}
            \ar@/^2pc/[rrrr]^F
            \ar@/_2pc/[rrrr]_H
            &
            \ar@3{->}[rr]_{\mu}
            &&&{\bcal V}
            \ar@/^2pc/[rrrr]^G
            \ar@/_2pc/[rrrr]_G
            &
            \ar@3{->}[rr]_{1_{1_G}}
            &&&{\bcal W}\\
            &&&&&&&&
            }
    $$
  \clearpage
  On the other side we need to check that the following compositions are equivalent:
  \noindent
          		          \begin{center}
  	          \resizebox{6.55in}{!}{ %!!!!
	          %%%%\resizebox{7.45in}{!}{
  $$
       \xymatrix@R=10pt@C=20pt{
       &
       &\ar@{=>}[dd]_{\alpha}
       &
       &&\ar@/_.5pc/@{=>}[dd]_{\gamma}
       &&\ar@/^.5pc/@{=>}[dd]^{\rho}
       \\
       {\bcal U}
       \ar@/^2pc/[rrrr]^F
       \ar@/_2pc/[rrrr]_H
       &
       &&&{\bcal V}
       \ar@/^2pc/[rrrr]^G
       \ar@/_2pc/[rrrr]_K
       &
       \ar@3{->}[rr]_{\nu}
       &&&{\bcal W}\\
       &&&&&&&&
       }=
       \xymatrix@R=10pt@C=20pt{
            &\ar@/_.5pc/@{=>}[dd]_{\alpha}
            &&\ar@/^.5pc/@{=>}[dd]^{\alpha}
            &&\ar@/_.5pc/@{=>}[dd]_{\gamma}
            &&\ar@/^.5pc/@{=>}[dd]^{\rho}
            \\
            {\bcal U}
            \ar@/^2pc/[rrrr]^F
            \ar@/_2pc/[rrrr]_H
            &
            \ar@3{->}[rr]_{1_{\alpha}}
            &&&{\bcal V}
            \ar@/^2pc/[rrrr]^G
            \ar@/_2pc/[rrrr]_K
            &
            \ar@3{->}[rr]_{\nu}
            &&&{\bcal W}\\
            &&&&&&&&
            }
    $$
  }
  \end{center}
  
  This follows from the definition of ${\cal V}$--2--functor (relation to unit axiom) since 
  $(G\alpha)_{U_{00}} = j_{G(\alpha_U(0))} = (G1_{\alpha})_U.$ Furthermore, since
  $(G1_F)_{U_{00}} = {\cal J}_{GFU_{00}},$ we have the equality:
  
  $$
       \xymatrix@R=10pt@C=20pt{
       &
       &
       &
       &&\ar@/_.5pc/@{=>}[dd]_{\gamma}
       &&\ar@/^.5pc/@{=>}[dd]^{\rho}
       \\
       {\bcal U}
       \ar[rrrr]^F
       &
       &&&{\bcal V}
       \ar@/^2pc/[rrrr]^G
       \ar@/_2pc/[rrrr]_K
       &
       \ar@3{->}[rr]_{\nu}
       &&&{\bcal W}\\
       &&&&&&&&
       }$$= $$
       \xymatrix@R=10pt@C=20pt{
            &
            &\ar@{=>}[dd]_{1_F}
            &
            &&\ar@/_.5pc/@{=>}[dd]_{\gamma}
            &&\ar@/^.5pc/@{=>}[dd]^{\rho}
            \\
            {\bcal U}
            \ar@/^2pc/[rrrr]^F
            \ar@/_2pc/[rrrr]_F
            &
            &&&{\bcal V}
            \ar@/^2pc/[rrrr]^G
            \ar@/_2pc/[rrrr]_K
            &
            \ar@3{->}[rr]_{\nu}
            &&&{\bcal W}\\
            &&&&&&&&
            }$$ = $$
            \xymatrix@R=10pt@C=20pt{
                 &\ar@/_.5pc/@{=>}[dd]_{1_F}
                 &&\ar@/^.5pc/@{=>}[dd]^{1_F}
                 &&\ar@/_.5pc/@{=>}[dd]_{\gamma}
                 &&\ar@/^.5pc/@{=>}[dd]^{\rho}
                 \\
                 {\bcal U}
                 \ar@/^2pc/[rrrr]^F
                 \ar@/_2pc/[rrrr]_F
                 &
                 \ar@3{->}[rr]_{1_{1_F}}
                 &&&{\bcal V}
                 \ar@/^2pc/[rrrr]^G
                 \ar@/_2pc/[rrrr]_K
                 &
                 \ar@3{->}[rr]_{\nu}
                 &&&{\bcal W}\\
                 &&&&&&&&
            }
  $$
  \clearpage
  Associativity of this composition follows from the associativity of composing ${\cal V}$--modifications
  along a ${\cal V}$--2--functor. It also requires the functoriality of the partial functors that describe whiskering
  ${\cal V}$--2--functors onto ${\cal V}$--modifications. In the following picture
  
  $$
       
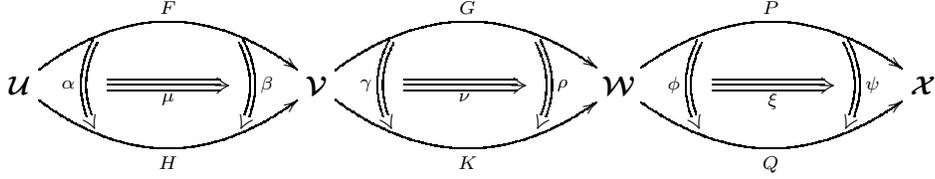
\xymatrix@R=10pt@C=20pt{
       &\ar@/_.5pc/@{=>}[dd]_{\alpha}
       &&\ar@/^.5pc/@{=>}[dd]^{\beta}
       &&\ar@/_.5pc/@{=>}[dd]_{\gamma}
       &&\ar@/^.5pc/@{=>}[dd]^{\rho}
       &&\ar@/_.5pc/@{=>}[dd]_{\phi}
       &&\ar@/^.5pc/@{=>}[dd]^{\psi}
       \\
       {\bcal U}
       \ar@/^2pc/[rrrr]^F
       \ar@/_2pc/[rrrr]_H
       &
       \ar@3{->}[rr]_{\mu}
       &&&{\bcal V}
       \ar@/^2pc/[rrrr]^G
       \ar@/_2pc/[rrrr]_K
       &
       \ar@3{->}[rr]_{\nu}
       &&&{\bcal W}
       \ar@/^2pc/[rrrr]^P
       \ar@/_2pc/[rrrr]_Q
       &
       \ar@3{->}[rr]_{\xi}
       &&&{\bcal X}
       \\
       &&&&&&&&&&&
       }
    $$

  we have $$\xi(\nu\mu) = Q(K\mu * \nu F)* \xi GF$$
  $$= (QK\mu * Q \nu F) * \xi GF$$
  $$= QK\mu * (Q \nu F * \xi GF)$$
  $$= QK\mu * (Q \nu  * \xi G)F = (\xi\nu)\mu$$
  where the assumed associativities of whiskers are easily verified.
  
  Finally the unit for this composition is given by the ${\cal V}$--modification 
  $1_{1_{1_{\bcal U}}}$ as in the following picture:
  
  $$
  
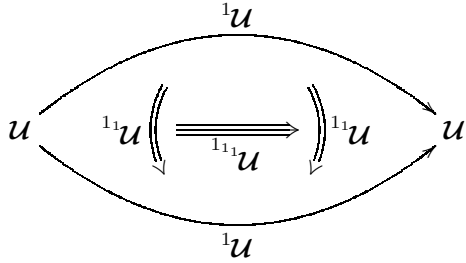
\xymatrix@R=10pt@C=20pt{
                 &&\ar@/_.5pc/@{=>}[dd]_{1_{1_{\bcal U}}}
                 &&\ar@/^.5pc/@{=>}[dd]^{1_{1_{\bcal U}}}
                 \\
                 {\bcal U}
                 \ar@/^3pc/[rrrrrr]^{1_{\bcal U}}
                 \ar@/_3pc/[rrrrrr]_{1_{\bcal U}}
                 &&
                 \ar@3{->}[rr]_{1_{1_{1_{\bcal U}}}}
                 &&&&{\bcal U}
                 \\
                 &&&&&&&&&&
                 }
  $$               
  It is straightforward to check that this is a 2--sided unit for the composition of ${\cal V}$--modifications
  along a common ${\cal V}$--2--category once we recognize that 
  $(1_{1_{1_{\bcal U}}})_U:I\to {\bcal U}(U,U)(1_U,1_U)$ is the morphism $j_{1_U}$ in ${\cal V}.$
  %\clearpage
  \end{proof}
  
  I close with the basic pasting diagram that the above proof has shown to be well-defined. There are 4 exchange 
  identities that this well--definedness depends upon, the requirements for each of which have been met.
  
  $$
  
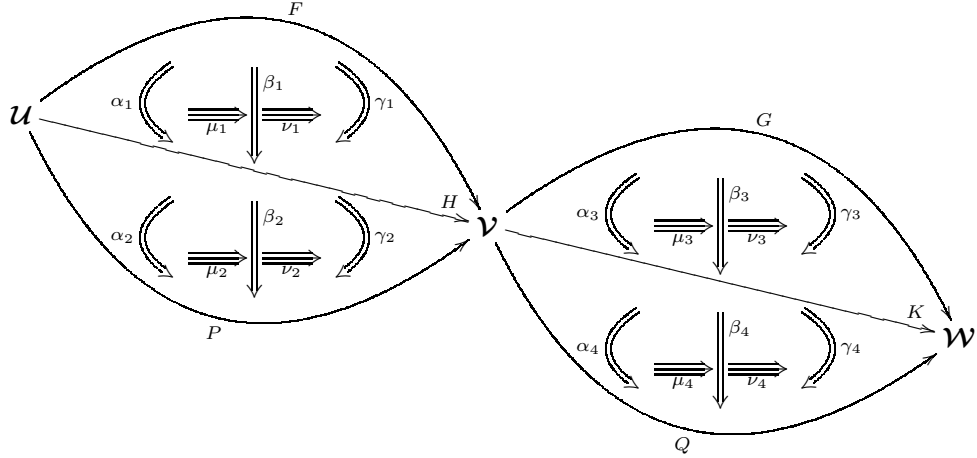
\xymatrix@R=2.8pt@C=22pt{
  && 
  \ar@{=>}@/_1.2pc/[ddd]_{\alpha_1}
  & 
  \ar@{=>}[dddd]^<<<{\beta_1}
  & 
  \ar@{=>}@/^1.2pc/[ddd]^{\gamma_1}
  \\\\
  {\bcal U}
  \ar@/_4.8pc/[rrrrrrdddd]_{P}
  \ar[rrrrrrdddd]^>>>>{H}
  \ar@/^4.8pc/[rrrrrrdddd]^{F}
  && 
  \ar@3{->}[r]_{\mu_1}
  & 
  \ar@3{->}[r]_{\nu_1}
  & 
  \\
  && 
  && 
  \\
  &&& 
  &&&&& 
  \ar@{=>}@/_1.2pc/[ddd]_{\alpha_3}
  & 
  \ar@{=>}[dddd]^<<<{\beta_3}
  & 
  \ar@{=>}@/^1.2pc/[ddd]^{\gamma_3}
  \\
  && 
  \ar@{=>}@/_1.2pc/[ddd]_{\alpha_2}
  & 
  \ar@{=>}[dddd]^<<<{\beta_2}
  & 
  \ar@{=>}@/^1.2pc/[ddd]^{\gamma_2}
  \\
  &&&&&&
  {\bcal V}
  \ar@/_4.8pc/[rrrrrrdddd]_{Q}
  \ar[rrrrrrdddd]^>>>>{K}
  \ar@/^4.8pc/[rrrrrrdddd]^{G}
  && 
  \ar@3{->}[r]_{\mu_3}
  & 
  \ar@3{->}[r]_{\nu_3}
  & 
  \\
  && 
  \ar@3{->}[r]_{\mu_2}
  & 
  \ar@3{->}[r]_{\nu_2}
  & &&&& && 
  \\
  && && &&&&& 
  \\
  &&& 
  &&&&& 
  \ar@{=>}@/_1.2pc/[ddd]_{\alpha_4}
  & 
  \ar@{=>}[dddd]^<<<{\beta_4}
  & 
  \ar@{=>}@/^1.2pc/[ddd]^{\gamma_4}
  \\
  &&&&&&&&&&&&{\bcal W}
  \\
  &&&&&&&& 
  \ar@3{->}[r]_{\mu_4}
  & 
  \ar@3{->}[r]_{\nu_4}
  & 
  &&
  \\
  &&&&&&&& && 
  \\
  &&&&&&&&& 
  \\
  }
  $$
  
  Thus we have:
  $$
  (\alpha_4\alpha_2) * (\alpha_3\alpha_1) = (\alpha_4 * \alpha_3)(\alpha_2*\alpha_1)
  $$
  $$
  (\nu_2\circ\mu_2)*(\nu_1\circ\mu_1)=(\nu_2*\nu_1)\circ(\mu_2*\mu_1)
  $$
  $$
  (\nu_3\nu_1)\circ(\mu_3\mu_1)=(\nu_3\circ\mu_3)(\nu_1\circ\mu_1)
  $$
  $$
  (\mu_4\mu_2)*(\mu_3\mu_1)=(\mu_4*\mu_3)(\mu_2*\mu_1).
$$

  %{bib2}

    % {bib}
    
    %\clearpage
    %\newpage
    
}
\end{document}